\documentclass[12pt]{article}
\usepackage{mathrsfs}
\usepackage{amsthm,amsfonts,amssymb,amsmath}
\usepackage{fullpage}
\usepackage{enumerate}
\usepackage[colorlinks=true]{hyperref}
\usepackage{tikz-cd}
\usepackage[new]{old-arrows}
\usepackage[english]{babel}
\usepackage[autostyle]{csquotes}

\numberwithin{equation}{section}

\newcommand{\an}{{\mathrm{an}}}

\newcommand{\Pic}{\widehat{\mathrm{Pic}}}

\newcommand{\hdeg}{\widehat{\deg}}

\newcommand{\Div}{{\mathrm{Div}}}

\newcommand{\Spec}{\operatorname{Spec}}

\renewcommand{\div}{{\mathrm{div}}}

\newcommand{\ord}{{\mathrm{ord}}}

\newcommand{\supp}{{\mathrm{supp}}}

\newcommand{\norm}[1]{\|{#1}\|}

\newcommand{\ol }{\overline}

\newcommand\hdiv{\widehat{\mathrm{div}}}

\newcommand\Frac{\mathrm{Frac}\,}

\newtheorem{thm}{Theorem}[section]

\newtheorem{theorem}[thm]{Theorem}
\newtheorem{proposition}[thm]{Proposition}
\newtheorem{lemma}[thm]{Lemma}
\newtheorem{corollary}[thm]{Corollary}

\theoremstyle{definition}
\newtheorem{definition}[thm]{Definition}

\theoremstyle{remark}
\newtheorem{remark}[thm]{Remark}

\begin{document}

\title{Numerical triviality of adelic line bundles on quasi-projective varieties}
\author{Lai Shang}
\date{}
\maketitle

\begin{abstract}
    Yuan--Zhang proved the arithmetic Hodge index theorem for adelic line bundles 
on projective varieties over number fields and finitely generated fields. Together 
with Qu--Yin's pseudo-effectivity results, this establishes a numerical 
triviality criterion for adelic line bundles over number fields. In this paper, we further generalize the numerical 
triviality criterion for adelic line bundles to quasi-projective varieties over 
finitely generated fields. We also prove
an arithmetic Hodge index type theorem in the quasi-projective setting which requires stronger conditions. The curve case combines 
admissible metrics and the N\'eron--Tate height pairing with constancy arguments at archimedean and
non-archimedean places. The higher-dimensional height criterion is proved using 
pencil maps constructed by Xie--Yuan. The Hodge index type statement is then reduced to complete intersection curves. 
\end{abstract}
\setcounter{tocdepth}{4}
\tableofcontents

\section{Introduction}
In \cite{YZ13,YZ17}, Yuan and
Zhang established the Hodge index theorem for adelic line bundles on projective varieties,
first over number fields and later over finitely generated fields. This extends the work of Faltings 
and Hriljac \cite{Faltings1984,Hriljac1985} for arithmetic surfaces, who independently
related the arithmetic intersection pairing of degree-zero divisors
with the N\'eron--Tate height pairing on the Jacobian. As for hermitian line bundles, Moriwaki first extended the Hodge index theorem to higher-dimensional 
arithmetic varieties in \cite{Mor96}. 
  
The purpose of this paper 
is to establish a characterization of numerical triviality of adelic 
line bundles on quasi-projective varieties over number 
fields, and more generally, finitely generated 
fields. The case of projective varieties is solved by a 
combination of Yuan--Zhang's work on the arithmetic Hodge index 
theorem introduced above and Qu--Yin's work \cite{QY24} on pseudo-effectivity and essential minimum through their 
arithmetic Demailly approximation. Together, these results lead to the following numerical triviality criterion.
\begin{theorem}[\cite{YZ17,QY24}]
    Let $X$ be a normal geometrically integral projective variety of dimension $d\geq 1$ over a number field $K$. Let $\overline{M}$ be an integrable adelic line bundle on $X$ and $\pi:X\rightarrow \Spec K$ be the structure morphism. The following statements are equivalent:
    \begin{enumerate}[$\mathrm{(}a\mathrm{)}$]
        \item $\overline{M} \cdot \overline{L}_1 \cdots \overline{L}_d = 0$, for any integrable adelic line bundles $\overline{L}_1, \dots, \overline{L}_d$ on $X$.
        \item  The height function $h_{\overline{M}} : X(\overline{K}) \to \mathbb{R}$ is identically 0.  
        \item $\overline{M}\simeq \pi^* \overline{N}$, for an adelic $\mathbb{Q}$-line bundle $\overline{N}$ over $K$ with $\hdeg(\overline{N}) = 0$.
    \end{enumerate}
\end{theorem}
    The implication $(a)\Rightarrow (c)$ is due to \cite[Theorem 2.2]{YZ17}, where we let 
    $\overline{L}_1=\overline{M}$ and obtain the geometric numerical triviality from $(a)$. The implication  
    $(b)\Rightarrow (a)$ holds for any quasi-projective $V$, and is achieved by applying \cite[Theorem A.1]{QY24} to 
    both $-\overline{M}$ and 
    $\overline{M}$ to show the pseudo-effectivity of $\pm \overline{M}$, which implies $(a)$ by taking 
    $\overline{L}_i$ to be nef and expressing arbitrary integrable $\overline{L}_i$ as a difference of nef ones.
    The projection formula proves $(c) \Rightarrow (a)$, while $(c) \Rightarrow (b)$ follows from Lemma \ref{height2}.
    
    In \cite{Car20}, Carney obtained a similar criterion 
    for numerically trivial adelic line bundles on projective varieties over a 
    function field $K=k(B)$, where $B/k$ 
    is an algebraic curve over an algebraically closed field $k$. In that case, the main difference is that there are 
    additional numerically trivial line bundles coming from the constant field $k$ via Chow's $K/k$-trace functor; see \cite[Theorem 3.2]{Car20}.
      
Let $F$ be a finitely generated field over $\mathbb{Q}$ of 
 transcendence degree $r\geq 0$ (number field if $r=0$). When $X$ is a projective variety over $F$, the implication $(a)\Rightarrow(c)$ is 
 already known from \cite[Theorem 1.2]{YZ13}. Here in the finitely generated field case, $(a)$ should be 
 interpreted as $\pi_{*}\langle \overline{M},  \overline{L}_1, \cdots, \overline{L}_d \rangle \equiv_{\mathrm{num}} 0$. Here, the Deligne pairing  
 $\pi_{*}\langle \overline{M},  \overline{L}_1, \cdots, \overline{L}_d \rangle$ is an adelic line bundle on 
 $\Spec F$, which is introduced in \cite[$\S$2.5$\,\&\,\S$4.5]{YZ26}. We will show the numerical 
 triviality criterion for quasi-projective varieties over finitely generated fields in the following theorem:

\begin{theorem}\label{num-triv}
     Let $V$ be a normal geometrically integral quasi-projective variety of dimension $d\ge 1$ over a finitely generated field $F$ of 
     $\mathrm{tr.deg}_{\mathbb{Q}}(F)=r\ge 0$. Let $\pi_V:V\rightarrow \Spec F$ be its structure morphism. Let $\overline{M}$ be an integrable adelic 
     line bundle on $V$ and $\mathfrak{h}_{\overline{M}} : V(\overline{F}) \to \widehat{\mathrm{Pic}}(F/\mathbb{Z})_{\mathrm{int},\mathbb{Q}}$ be the 
     vector-valued height associated with $\overline{M}$. The following are equivalent:
    \begin{enumerate}[$\mathrm{(}a\mathrm{)}$]
        \item $\pi_{V,*}\langle \overline{M},  \overline{L}_1, \cdots, \overline{L}_d \rangle \equiv_{\mathrm{num}} 0$, for any integrable adelic line bundles $\overline{L}_1, \dots, \overline{L}_d$ on $V$.
        \item $\mathfrak{h}_{\overline{M}}(x)\equiv_{\mathrm{num}} 0$ in $\widehat{\mathrm{Pic}}(F/\mathbb{Z})_{\mathrm{int},\mathbb{Q}}$, for all $x\in V(\overline{F})$.  
        \item $\overline{M}\simeq \pi_V^* \overline{N}$ for an adelic $\mathbb{Q}$-line bundle $\overline{N}\in \widehat{\mathrm{Pic}}(F/\mathbb{Z})_{\mathrm{int},\mathbb{Q}}$ with $\overline{N}\equiv_{\mathrm{num}} 0$. 
        \end{enumerate}
\end{theorem}
Here we write $\overline{N}\equiv_{\mathrm{num}} 0$ if $\overline{N}\cdot\overline{H}_1\cdots\overline{H}_r=0$ for arbitrary $r$ 
integrable adelic line bundles $\overline{H}_1,\cdots,\overline{H}_r$ on $\Spec F$. The notion of vector-valued 
height is introduced in \cite[$\S$5.3]{YZ26}, which generalizes the Moriwaki 
height in \cite{Mor00}. We will describe some basic properties in $\S$\ref{vector_heights}.

 Now we recall that for $V/F$ as above and an arbitrary adelic line bundle $\overline{M}$ on $V$, the geometric part $\widetilde{M}$ of $\overline{M}$ is defined to be the image of $\overline{M}$ under the varying-base map
\begin{equation*}
    \begin{aligned}
        \widehat{\mathrm{Pic}}(V/\mathbb{Z}) &\longrightarrow \widehat{\mathrm{Pic}}(V/F), \\
       \overline{M}~~~ &\longmapsto ~~~\widetilde{M}.
    \end{aligned}
\end{equation*}
See \cite[$\S$2.5.5]{YZ26} for details. To show the $(a) \Rightarrow (c)$ part of Theorem \ref{num-triv}, we will propose an arithmetic Hodge index type theorem for quasi-projective 
varieties over finitely generated fields as follows. 
\begin{theorem}[Theorem \ref{Rdr}]\label{chap1-Rdr}
Let $V$ be a normal geometrically integral quasi-projective variety of dimension $d\ge 1$ over a finitely generated field $F$ of 
     $\mathrm{tr.deg}_{\mathbb{Q}}(F)=r\ge 0$. Suppose
$\overline M$ is an integrable adelic line bundle on $V$.
Assume that for every collection
 $\overline A_1,\ldots,\overline A_{d-1} \in \widehat{\mathrm{Pic}}(V/\mathbb{Z}) _{\mathrm{int},\mathbb Q}$, one has
\[
 \widetilde M\cdot
 \widetilde A_1\cdots\widetilde A_{d-1}=0
\]
and
\[
 \pi_{V,*}\bigl\langle
 \overline M,\overline M,
 \overline A_1,\ldots,\overline A_{d-1}
 \bigr\rangle
 \equiv_{\mathrm{num}}0.
\]
When $d=1$, the first equality is interpreted as
$\deg(\widetilde M)=0$. Then
\[
 \overline M\in
 \pi_V^*\widehat{\mathrm{Pic}}(F/\mathbb{Z})_{\mathbb Q}.
\]
\end{theorem}
\begin{remark}
     When $d=1,~r=0$, the numerical 
    triviality of $\pi_{V,*}\langle \overline{M},\overline{M} \rangle$ is equivalent to $\overline{M}^2=0$, by the projection formula in \cite[Lemma 4.4.3(1)]{YZ26}. And for general $d,~r$, the Deligne pairing $ \pi_{V,*}\bigl\langle
 \overline M,\overline M,
 \overline A_1,\ldots,\overline A_{d-1}\rangle\equiv_{\mathrm{num}} 0$ is equivalent to 
$\overline{M}^2\cdot\overline{A}_1\cdots\overline{A}_{d-1}\cdot \pi_V^*\overline{H}_1\cdots \pi_V^*\overline{H}_r=0$ for all integrable adelic line bundles $\overline{A}_i$ on $V$ and 
 $\overline{H}_i$ on $\Spec F$ by the same projection formula.
As adelic line bundles can be approximated by model 
adelic line bundles in the boundary topology, it suffices to let these $\overline{A}_i$ run over all model adelic line bundles on $V$. The implication $(a) \Rightarrow (c)$ as a corollary of this theorem is shown in $\S$\ref{Hodge-implies-a-to-c}.
\end{remark}

\subsection{Proof strategy}

Let $V$ be a 
normal, geometrically integral quasi-projective variety of pure dimension $d\geq 1$ over a finitely generated field $F$ of 
$\mathrm{tr.deg}_{\mathbb{Q}}F=r\ge 0$. 
The proof depends on two families of statements regarding numerical triviality of heights and rigidity, which are 
denoted by $H(d,r)$ and $R(d,r)$ respectively. The assertion $H(d,r)$ claims that if $\mathfrak h_{\overline M}(x)$ is numerically trivial for
every algebraic point $x$ in $V(\ol F)$, then $\overline M$ is a pullback of an adelic line bundle on $\Spec F$. The
assertion $R(d,r)$ exactly states that if $\widetilde{M} \cdot
 \widetilde A_1\cdots\widetilde A_{d-1}=0$ and 
 $$\overline{M}^2\cdot\overline{A}_1\cdots\overline{A}_{d-1}\cdot \pi_V^*\overline{H}_1\cdots \pi_V^*\overline{H}_r=0$$
 for all integrable adelic line bundles $\overline{A}_i$ on $V$ and 
 $\overline{H}_i$ on $\Spec F$, then we have $\overline M\in
 \pi_V^*\widehat{\mathrm{Pic}}(F/\mathbb{Z})_{\mathbb Q}$, for all $V/F$ of $\dim V=d$ and 
 $\mathrm{tr.deg}_{\mathbb{Q}}F= r$. These will be completely proved in \S\ref{Hdr} and \S\ref{chap5-Rdr}. 

We start by proving the curve case $H(1,r)$ by reducing it to proving $R(1,r)$ as in $\S$\ref{Q1r}, then obtain the 
general $H(d,r)$ by induction via the methods in \cite{XY22} which lower the 
dimension of varieties by $1$ and add $1$ to the transcendence degree of the 
base field. And we prove that $H(d,r)$ and $R(1,r)$ together imply that $R(d,r)$ holds in $\S$\ref{subsec-Rdr general}. 

\begin{theorem}[$H(1,r)$ holds]
\label{thm:intro-Q1r}
Let $C$ be a normal geometrically integral quasi-projective curve over
$F$, let $\pi_C\colon C\longrightarrow\operatorname{Spec}F$
be the structure morphism, and let
\[
 \overline M\in
 \widehat{\mathrm{Pic}}(C/\mathbb{Z})
 _{\mathrm{int},\mathbb Q}.
\]
If
\[
 \mathfrak h_{\overline M}(x)\equiv_{\mathrm{num}}0,
 \qquad\text{for every }x\in C(\overline F),
\]
then
\[
 \overline M\in
 \pi_C^*\widehat{\mathrm{Pic}}(F/\mathbb{Z})_{\mathbb Q}.
\]
\end{theorem}

The proof of $H(1,r)$ first converts the pointwise numerical triviality of
vector-valued heights into geometric degree zero and into the numerical
vanishing of the self-pairing, thereby reducing it to $R(1,r)$.

We prove that $H(d,r)$ holds for all $d\geq 1$ and $r\geq 0$, in the following theorem:
\begin{theorem}[Theorem \ref{H(d,r)}] \label{chap1-Hdr}
  Let $F$ be a finitely generated 
field of transcendence degree $r$ over $\mathbb{Q}$, and $V/F$ be a geometrically integral, normal quasi-projective 
variety of dimension $d$. For an integrable adelic $\mathbb{Q}$-line bundle 
$ \overline{M} \in \widehat{\mathrm{Pic}}(V/\mathbb{Z})_{\mathrm{int},\mathbb{Q}} $, if the vector-valued height 
$$
\mathfrak{h}_{\overline{M}}(x)\equiv_{\mathrm{num}} 0 \in \widehat{\mathrm{Pic}}(F/\mathbb{Z})_{\mathbb{Q}},
$$
for every $x\in V(\overline{F})$, then we have
 $$ \overline{M}\in\pi_{V}^{*}\widehat{\mathrm{Pic}}(F/\mathbb{Z})_{\mathbb{Q}}. 
 $$
\end{theorem}

This is obtained by the following induction argument:

\begin{theorem}[Induction]
\label{thm:intro-Q-induction}
For every $d\geq 2$ and $r\geq 0$,
\[
\mathrm{if~} H(d-1,r+1) \mathrm{~holds,~then~}  H(d,r)~\mathrm{holds}.
\]
Consequently, the curve case $H(1,s)$ for all $s\geq0$ implies
$H(d,r)$ for every $d\geq1$ and $r\geq0$.
\end{theorem}

The induction in Theorem~\ref{thm:intro-Q-induction} uses the family of pencil maps constructed in \cite[Lemma 3.7]{XY22}. In this way, we can produce varieties $V_i$ of dimension $d-1$ over fields
$F_i$ of transcendence degree $1$ over $F$ for $i=1,2$. We can identify $F_i$ as subfields of the function field $F(V)$ of $V/F$. We can appropriately choose $F_1,~F_2$ so that
\[
 \mathrm{tr.deg}_{F}F_1F_2=2,\quad\mathrm{and}\quad F_1\cap F_2=F.
\]
A key argument on the relative Deligne pairing of adelic line bundles shows that the numerical triviality
of vector-valued heights associated with $\overline{M}$ along all algebraic points of $V$ is preserved on the generic fibers through this pencil map. The hypothesis of
$H(d-1,r+1)$ therefore implies that the restriction of $\overline M$
to $\Spec F(V)$ is a pullback from both $F_1$ and $F_2$. A descent argument in Lemma \ref{Descent} then
forces $\overline M$ to be a pullback from $F$. Finally, injectivity under the restriction to the
generic point gives $H(d,r)$. We will apply $H(d,r)$ to prove the arithmetic Hodge index type theorem \ref{chap1-Rdr} that 
claims that $R(d,r)$ holds. We will treat the curve case $R(1,r)$ first and reduce the higher dimensional case to the $d=1$ case of Theorem \ref{chap1-Rdr}. The assertion $R(1,r)$ used in this argument is the following.

\begin{theorem}[$R(1,r)$~holds]
\label{thm:intro-R1r}
Let $C$ be a normal geometrically integral quasi-projective curve over
$F$, and let $\pi_C\colon C\rightarrow\operatorname{Spec}F$
be the structure morphism. Let $\overline M\in
 \widehat{\mathrm{Pic}}(C/\mathbb{Z})
 _{\mathrm{int},\mathbb Q}$ be an integrable adelic line bundle satisfying
\[
 \deg(\widetilde M)=0
 \qquad\text{and}\qquad
 \pi_{C,*}\langle\overline M,\overline M\rangle
 \equiv_{\mathrm{num}}0.
\]
Then we have
\[
 \overline M\in
 \pi_C^*\widehat{\mathrm{Pic}}(F/\mathbb{Z})_{\mathbb Q}.
\]
\end{theorem}

For $r=0$, the assertion $R(1,0)$ follows from the
case of curves over a number field. For $r\geq1$, the curve case
$R(1,r)$ relies on an explicit calculation of the local Green functions. After passing to a projective model and modifying the model adelic line bundles to have
degree zero while still approximating $\overline{M}$, we can then apply results on admissible metrics and the N\'eron--Tate height
pairing in \cite{YZ26}. We will show that the vanishing of the global self-intersection first forces $\widetilde{M}$ to be trivial. The remaining vertical class is then
controlled by arguments on finite and infinite places of $F$. All local contributions are non-positive, thereby tending to zero. Poincar\'e's inequality
in the archimedean case and a pathwise Cauchy--Schwarz argument in the
non-archimedean case show that the limiting Green functions are
constant. This shows that $\overline{M}$ is a pullback from $\operatorname{Spec}F$.

Finally, $R(d,r)$ is deduced from $H(d,r)$ and $R(1,r)$. Given two
regular closed points of $V$, we can choose a normal, geometrically 
integral complete intersection curve $C$ passing through them by \cite{Sei50}. 
Then we introduce the notion of b-divisors associated with an adelic divisor as the limit of push-forwards of the 
geometric parts of model adelic divisors on dominating models which 
converge to it, which has been discussed in \cite{BGK24} and \cite{BBGHdJ21}. It leads us to show that  
our $\overline{M}$ is indeed derived from a model adelic line bundle endowed with 
flat metrics. Then we show that the hypothesis of
$R(d,r)$ is restricted to the hypothesis of $R(1,r)$ on this curve; hence
$\overline M|_C$ comes from the base, and any two
points on the same curve have the same vector-valued height. It follows that
$\mathfrak h_{\overline M}$ is constant on $V(\overline F)$. After
subtracting the pullback of $\mathfrak h_{\overline M}(x_0)$ for some $x_0\in V(\overline{F})$, we obtain an adelic
line bundle with numerically trivial height at every algebraic point.
The assertion $H(d,r)$ shows that
$\overline M$ itself is a pullback from $\Spec F$.

We conclude by mentioning 
two recent local developments related to the
quasi-projective setting. In \cite{Guo25}, Guo proves an integration formula for
integrable adelic line bundles on a quasi-projective variety over an
archimedean or non-archimedean complete valuation field, identifying
the integral of the Chambert-Loir measure \cite{Gub98,CL06} with the
corresponding geometric intersection number. In \cite{Abboud25}, Abboud defines a local 
intersection pairing
for compactly supported vertical metrized line bundles on
quasi-projective varieties and proves a local Hodge index theorem,
together with an equality criterion and an application to uniqueness
for the Monge--Amp\`ere equation in his theorems $A$ and $B$. 

\noindent \textbf{AI disclosure statement.} I used ChatGPT to assist in adapting the analytic methods of \cite{FRRO22} 
on elliptic regularity in \S\ref{chap5-vertical} to obtain the smoothness of 
Green functions, and adapting the algebraic geometry methods to refine the blowing 
up arguments in Lemma \ref{lem:reduce-to-surface} 
by applying Haution's theorem \cite{Hau15}.
I also used ChatGPT to check for typos, grammar, and mathematical correctness of 
proofs.  I take full responsibility for the mathematical arguments and all text in this paper.
\begin{remark} One may mistake an adelic divisor which is an irrational 
multiple of a principal divisor supported on the boundary as a counterexample 
for Theorem \ref{num-triv}. However, this adelic divisor corresponds to a trivial 
adelic line bundle by definition, therefore, causing no contradiction.
\end{remark}

\subsection{Adelic line bundles}\label{preliminiary}
We recall the constructions of adelic line bundles over finitely generated fields that will be used in this  paper, mainly following 
\cite{YZ26}.   
We first work on a normal quasi-projective arithmetic variety
$\mathcal{U}$, that is, a normal integral scheme which is flat and
quasi-projective over $\Spec \mathbb{Z}$. A projective compactification of $\mathcal{U}$ is an open immersion
$\mathcal{U}\hookrightarrow\mathcal{X}$ with $\mathcal{X}$ projective
and flat over $\Spec \mathbb{Z}$ and we say that $\mathcal{X}$ is a projective model of $\mathcal{U}$. Let $U:= \mathcal{U}_{\mathbb{Q}}$ be the generic fiber of $\mathcal{U}$.



\begin{definition}
A \textit{model adelic line bundle} on $\mathcal{U}$ is a metrized line bundle on $\mathcal{U}^{\mathrm{an}}$ induced by a
hermitian $\mathbb{Q}$-line bundle $\overline{\mathcal{M}}$ on a
projective model $\mathcal{X}$ of $\mathcal{U}$, together with
an identification
\[
    \ell:M\xrightarrow{\sim}\mathcal{M}|_{\mathcal{U}}
\]
in $\mathrm{Pic}(\mathcal{U})_{\mathbb{Q}}$.
Here the hermitian metric is continuous and invariant under complex
conjugation. 
\end{definition}

\begin{definition}
An \textit{adelic line bundle} $\overline{M}$ on a quasi-projective arithmetic variety 
$\mathcal{U}$ is given by a line bundle $\mathcal{M}$ on $\mathcal{U}$ and a sequence of model adelic line bundles $\{(\mathcal{X}_n, \overline{\mathcal{M}}_n, \ell_n)\}_{n \ge 1}$ where:
\begin{itemize}
    \item $\mathcal{X}_n$ is a projective arithmetic compactification
    of $\mathcal{U}$;
    \item $\overline{\mathcal{M}}_n$ is a hermitian
    $\mathbb{Q}$-line bundle on $\mathcal{X}_n$;
    \item $\ell_n:\mathcal{M}\xrightarrow{\sim}\mathcal{M}_n|_{\mathcal{U}}$
    is an identification in $\mathrm{Pic}(\mathcal{U})_{\mathbb{Q}}$;
    \item there are positive rational numbers $\varepsilon_n\to0$
    such that
    \begin{equation}\label{eq:prelim-Cauchy}
        -\varepsilon_n\overline{\mathcal{E}}_0
        \leq \hdiv(\ell_m\ell_n^{-1})
        \leq \varepsilon_n\overline{\mathcal{E}}_0,
        \qquad m\geq n.
    \end{equation}
\end{itemize}
Here $\ell_m\ell_n^{-1}$ is viewed as a rational section of
$\mathcal{M}_m\otimes \mathcal{M}_n^{\vee}$ on a common model and $\overline{\mathcal{E}}_0$ is a fixed strictly effective boundary divisor on a model $\mathcal{X}_0$ with support $\supp(\mathcal{E}_0) = |\mathcal{X}_0 \setminus \mathcal U|$. Two representatives are
equivalent if, under one fixed identification of their underlying line
bundles, their difference converges to zero in the same boundary
sense.
\end{definition}

We denote by $\widehat{\mathrm{Pic}}(\mathcal{U}/\mathbb{Z})$ the group of isomorphism 
classes of adelic line bundles on $\mathcal{U}$, and add subscripts such as \enquote{int, nef, 
snef} to represent the subsets of integrable, nef, strongly nef ones, which will be 
introduced in the end of this subsection. And we write $\Pic(K/\mathbb{Z}):= \Pic(\Spec O_K/\mathbb{Z})$ for simplicity when $K$ is a number field. 
The variety $V$ over $F$ in this paper is defined as an integral scheme $V$ that is flat, and of finite type over 
$\Spec F$. We say that $V$ is a quasi-projective variety over $F$ if $V$ is an open subscheme of a projective 
variety over $\Spec F$.

Over a finitely generated field $F/\mathbb{Q}$, let $V\rightarrow \Spec F$ be a 
quasi-projective variety. In this paper, we use the symbol $V$ to distinguish it from a
quasi-projective variety $U$ or quasi-projective arithmetic variety $\mathcal{U}$ over $O_K$. An arithmetic model of 
$V\to\Spec F$ is a flat morphism
\[
    \mathcal{V}\longrightarrow\mathcal{B}
\]
of integral quasi-projective arithmetic varieties such that $\mathcal{B}$ has function field $F$, and the generic fiber of $\mathcal{V}\rightarrow\mathcal{B}$
is $V\to\Spec F$. Following
\cite[\S1.1.4 $\&$ \S2.5.3]{YZ26}, we set
\begin{equation}\label{eq:prelim-direct-limit}
    \widehat{\mathrm{Pic}}(V/\mathbb{Z})
    :=\varinjlim_{\mathcal{V}\to\mathcal{B}}
        \widehat{\mathrm{Pic}}(\mathcal{V}/\mathbb{Z}).
\end{equation}
The transition maps are pullbacks on common dominating models and
restrictions to smaller open models containing the generic fiber.
Normal models form a cofinal system in the present setting. In particular,
\[
    \widehat{\mathrm{Pic}}(F/\mathbb{Z})
    :=\widehat{\mathrm{Pic}}(\operatorname{Spec}F/\mathbb{Z})
    =\varinjlim_{\mathcal{B}}
        \widehat{\mathrm{Pic}}(\mathcal{B}/\mathbb{Z}).
\]

A model adelic line bundle on $V/\mathbb{Z}$ means the image of a
model adelic line bundle on one of these arithmetic models
$\mathcal{V}$. Consequently, any finite collection of adelic line bundles on $V/F$
can be represented by adelic line bundles on a common arithmetic model $\mathcal{V}$ of $V$ and
then by a Cauchy sequence of hermitian $\mathbb{Q}$-line bundles on projective arithmetic models of $\mathcal{V}$. All projective 
compactifications dominating a
fixed one form a cofinal family. Consequently, the models in a representing sequence can be replaced by
a sequence of common dominating normal models; see also
\cite[Lemmas~2.3.3 $\&$ 2.3.4]{YZ26}.

Recall that the \emph{geometric part}
\[
    \widetilde{M}\in
       \widehat{\mathrm{Pic}}(V/F)_{\mathbb{Q}},
\]
is defined to be the image of $\overline{M}$ under the varying-base map of
\cite[\S2.5.5]{YZ26}. When $V$ is a curve, we have a well-defined degree:
\[
    \deg(\widetilde{M})
       :=\lim_{n\to\infty}\deg(\mathcal{M}_{n,F}).
\]

Finally, we introduce the notion of nefness in the quasi-projective setting:

\begin{definition}
    We say that an adelic line bundle $\overline{M}$ on $U$ is:
    \begin{enumerate}[$\bullet$]
        \item \textit{strongly nef}, if $\overline{M}$ is isomorphic to $(M, \{(\mathcal{X}_n, \overline{\mathcal{M}}_n, \ell_n)_{n \ge 1}\})$ with $\overline{\mathcal{M}}_n$ nef on $\mathcal{X}_n$;
        \item \textit{nef}, if there exists a strongly nef $\overline{L}$ such that $a\overline{M} + \overline{L}$ is strongly nef for all positive integers $a$;
        \item \textit{integrable}, if it is isomorphic to a difference of two strongly nef adelic line bundles on $U$.
    \end{enumerate}
\end{definition}

\subsection{Vector-valued heights}  \label{vector_heights}
In this subsection, we will introduce the vector-valued heights of 
$\overline{M}\in \widehat{\mathrm{Pic}}(V/\mathbb{Z})_{\mathrm{int},\mathbb{Q}}$ following \cite[$\S$5.3]{YZ26}.
Let $F/\mathbb{Q}$ be finitely generated of transcendence degree $r$,
let $V/F$ be a normal geometrically integral quasi-projective variety
of dimension $d$, and let
\[
 \pi_V\colon V\longrightarrow\operatorname{Spec}F
\]
be the structure morphism. For $x\in V(\overline F)$ and its associated closed point $x^\prime$ with residue
field $L$, an adelic line bundle $\overline M$ induces a vector-valued height
\[
 \mathfrak h_{\overline M}(x)
 :=
 \frac{1}{[L:F]}
 N_{L/F}(\overline M|_{x^\prime})
 \in
 \widehat{\mathrm{Pic}}(F/\mathbb{Z})
 _{\mathrm{int},\mathbb{Q}}.
\]
Here $N_{L/F}(\cdot)$ means the Deligne pairing of relative dimension $0$, 
which is introduced in \cite[$\S$4.5]{YZ26}. We say that
$\overline{N}\in\widehat{\mathrm{Pic}}(F/\mathbb{Z})
_{\mathrm{int},\mathbb{Q}}$ is \textit{numerically trivial}, and write
$\overline{N}\equiv_{\mathrm{num}}0$, if
\[
    \overline{N}\cdot\overline{H}_1\cdots\overline{H}_r=0
\]
for every such integrable tuple. When $r=0$, this means
$\widehat{\deg}(\overline{N})=0$, and we write
$h_{\overline{M}}(x)=\widehat{\deg}(\mathfrak{h}_{\overline{M}}(x))$. 

For the reader's convenience, we introduce and prove two elementary lemmas that will be repeatedly used later and are not explicitly stated in \cite{YZ26}.


\begin{lemma}\label{height1}
    Let $F$ be a finitely generated field over $\mathbb{Q}$, let $f: Y \rightarrow X$ be a morphism
of quasi-projective varieties over $F$, and let $\overline{M} \in \widehat{\mathrm{Pic}}(X/\mathbb{Z})_{\mathrm{int},\mathbb{Q}}$ be an adelic $\mathbb{Q}$-line bundle on $X$. Then for every $y \in Y(\overline{F})$, we have
$$
\mathfrak{h}_{f^*\overline{M}}(y)=\mathfrak{h}_{\overline{M}}(f(y)).
$$
\end{lemma}
\begin{proof}
    Let $y^\prime$ be the closed point of $Y$ corresponding to $y$ and 
    $x^\prime=f(y^\prime)$. Denote the residue fields by $E:= \kappa(y^\prime), ~L:=\kappa(x^\prime)$. Then 
    $f^*\overline{M}|_{y^\prime}$ is the pullback of 
    $\overline{M}|_{x^\prime}$ via $\Spec E \rightarrow\Spec L$. Therefore, we have
\begin{equation*}
\begin{aligned}
\mathfrak{h}_{f^*\overline{M}}(y)&= \frac{1}{[E:F]}N_{E/F}( f^*\overline{M}|_{y^\prime})  \\
&=\frac{1}{[E:F]}[E:L]  N_{L/F}(\overline{M}|_{x^\prime})\\
&=\frac{1}{[L:F]}N_{L/F}(\overline{M}|_{x^\prime})=\mathfrak{h}_{\overline{M}}(f(y)),
\end{aligned}
\end{equation*}
by applying \cite[Lemma 4.6.1(4)]{YZ26} to the second equality.
\end{proof}

\begin{lemma}\label{height2}
Suppose $\pi: V\rightarrow \Spec F$ is the structure 
morphism of a quasi-projective variety $V$. 
Let $\overline{N}$ be an integrable adelic line bundle on $\Spec F$. Then for all $x\in V(\overline{F})$, we have
$$
\mathfrak{h}_{\pi^*\overline{N}}(x)=\overline{N}.
$$
\end{lemma}
\begin{proof}
    Let $x^\prime$ be the closed point corresponding to $x$ with residue field $L:=\kappa(x^\prime)$. Suppose that 
    $\mathcal{B}$ is a quasi-projective model 
    of $F$ over $\mathbb{Z}$ such that there is a model $\mathcal{V}\rightarrow \mathcal{B}$ whose
    generic fiber isomorphic to $V\rightarrow \Spec F$. Let 
    $\mathcal{Z}$ be the Zariski closure of $x^\prime$ in $\mathcal{V}$. Then we have a morphism 
    $\alpha:\mathcal{Z}\to \mathcal{B}$. We can replace $\mathcal{B}$ by a smaller $\mathcal{B}^\circ$ and 
    $\mathcal{Z}$ by $\mathcal{Z}^\circ$ so that   
    $\alpha:\mathcal{Z}^\circ \rightarrow \mathcal{B}^\circ$ is finite and flat of degree 
    $[L:F]$, and $\overline{N}$ is represented by an 
    adelic line bundle in $\widehat{\mathrm{Pic}}(\mathcal{B}^\circ /\mathbb{Z})$, which we still denote by 
    $\overline{N}$. Then we have 
    \begin{equation*}
        \begin{aligned}
    \mathfrak{h}_{\pi^*\overline{N}}(x)&=\frac{1}{[L:F]}{N_{L/F}}(\pi^*\overline{N}|_{x^\prime})\\
    &=\frac{1}{[L:F]}\langle\alpha^*\overline{N}\rangle_{\mathcal{Z}^\circ/\mathcal{B}^\circ}\\
    &=\frac{1}{[L:F]}([L:F]\cdot\overline{N})=\overline{N},
        \end{aligned}
    \end{equation*}
by applying \cite[Lemma 4.6.1(4)]{YZ26} again.
\end{proof}


\section{Curves over number fields}
Let $C/K$ be a smooth geometrically integral and geometrically connected quasi-projective curve and let $K$ be a number field.
In this section, we will prove the following theorem, which states that $R(1,0)$ holds:
\begin{theorem}\label{curve-vert-K}
    Let $C$ be a smooth 
    quasi-projective curve over a number field $K$ and $\pi_C:C \rightarrow \Spec K$ be the structure morphism. Suppose $\overline{M}$ is an integrable adelic line bundle on $C$ such that $\deg(\widetilde{M})=0 $ and
    $$\overline{M}^2 = 0.$$
     Then there is an $\overline{N}\in \Pic(K/\mathbb{Z})_\mathbb{Q}$ such that $\overline{M}=\pi_C^*\overline{N}$.
\end{theorem}
 This theorem indicates the implication that numerical triviality $(a)$ implies 
$(c)$ for $1$-dimensional quasi-projective $C$ in Theorem \ref{num-triv}. Assuming $(a)$, we have
$$
\overline{M}\cdot \pi_C^*\overline{H}=\deg({\widetilde{M}})\cdot \widehat{\deg}(\overline H)= 0,
$$
for arbitrary adelic line bundle $\overline{H}$ on $\Spec O_K$ from the 
projection formula. Then we have $\deg({\widetilde{M}})=0$ by taking an $\overline{H}$ with $\widehat{\deg}(\ol H)> 0$. Moreover, by 
letting $\overline{L}_1=\overline{M}$ in $(a)$, we see that $\overline{M}$ satisfies the two conditions in Theorem \ref{curve-vert-K}. Therefore
$\overline{M}=\pi^*\overline{N}$ for some $\overline{N}\in \Pic(K/\mathbb{Z})_\mathbb{Q}$. Then we conclude that $\widehat\deg(\ol N)=0$ since we have
$$
\pi^*\overline{N}\cdot \overline{L}= \widehat\deg(\ol N)\cdot \deg(\widetilde{L})=0,
$$
for arbitrary adelic line bundle $L$ on $C$ by projection formula. This gives $(c)$ in Theorem  \ref{num-triv}.

\subsection{Vertical case} \label{Section-vertcurve}
In this section, we will prove the vertical case of Theorem \ref{curve-vert-K}. First, we introduce the notion of 
vertical integrable adelic line bundles on quasi-projective varieties.
\begin{definition}
     Let $K$ be a number field and $O_K$ be its ring of integers. A hermitian line bundle $\overline{\mathcal{M}}$ on an arithmetic variety $\mathcal{X}/O_K$ of generic fiber $X:=\mathcal{X}_K$ is called \textit{vertical} if the generic fiber $\mathcal{M}_K$ is isomorphic to the trivial line bundle $\mathcal{O}_{X}$ on $X$.
     
    An adelic line bundle $\overline{M}$ on a quasi-projective variety $V$ is called \textit{vertical} if $\overline{M}$ can be represented by a Cauchy sequence of model adelic line bundles $\{\overline{\mathcal{M}}_n\}$ whose generic fiber $\mathcal{M}_{n,K}$ is trivial for all $n$.

\end{definition}

We now prove the following proposition:

\begin{proposition}[Numerically trivial vertical adelic line bundles on curves] \label{curve-vertical}
Let $\overline{M}$ be a vertical integrable adelic line bundle on a smooth 
quasi-projective curve $C$. Suppose that 
$$
\overline{M}^2 = 0.
$$
Then 
$\overline{M}= \pi^* \overline{N}$, for some $ \overline{N}\in \widehat{\mathrm{Pic}}(K/\mathbb{Z})$. 
\end{proposition}

\begin{proof} For convenience, we represent $\overline{M}$ by an adelic $\mathbb{Q}$-Cartier divisor $\overline{D}$. Let
$X$ be a smooth projective model of $C$, and write 
$$ X\setminus C=\{p_1,\dots,p_r\}, $$
which is a finite set. Choose a Cauchy sequence of model adelic divisors 
$$ \{\overline{\mathcal{D}}_n \}_{n\geq 1} $$
representing $\overline{D}$ with $\mathcal{D}_{n,K}=0$. We may assume that 
$\overline{\mathcal{D}}_n$ is defined on an inverse system of projective 
arithmetic models such that all their generic fibers are identified with $X$.

For a place $v$ of $K$, let
$g_{n,v}:=-\log\|s\|_{n,v}$
be the Green function associated with $\overline{\mathcal D}_n$ and $g_v=-\log\|s\|_{\overline M,v}$ be the Green function 
associated with $\overline{D}$ on $C_v^{\an}$. Since 
$\mathcal{D}_{n,K}=0$, the horizontal part $\mathcal{D}_{n}^{\mathrm {hor}}=0$. We see that $g_{n,v}$ is a model 
function on the compact analytic curve $X_v^{\an}$. The  convergence of model adelic line bundles $\overline{\mathcal{D}}_n$ implies that
\[
g_{n,v}\longrightarrow g_v
\]
locally uniformly on $C_v^{\an}$. We aim to show that $g_v$ is a constant function for all $v\in M_K$.

For a vertical model adelic divisor $\overline{\mathcal{D}}_n$  
on a projective arithmetic surface, the self-intersection
$\overline{\mathcal{D}}_n^2$ decomposes into a sum of 
non-positive numbers at each place of $K$ as follows.
\begin{equation}\label{nonpositive-decomposition}
\begin{aligned}
    \overline{\mathcal D}_n^{2}
=&~\overline{\mathcal{D}}_n\cdot\mathcal{D}_{n}^{\mathrm {hor}}
+\sum_{v\in M_K}{\epsilon_v\int_{X_v^{\mathrm{an}}}{g_{n,v}\,c_1(\overline{\mathcal{D}}_n)}_v}\\
=&\sum_{v\in M_K}{\epsilon_v\int_{X_v^{\mathrm{an}}}{g_{n,v}\,c_1(\overline{\mathcal{D}}_n)}_v},
\end{aligned}
\end{equation}
where $\epsilon_v:={[K_v:\mathbb{Q}_v]}/[K:\mathbb{Q}]$ for every place $v$ of $K$.
Each term in the sum is non-positive, which is explained as follows: at an archimedean place $v$, the integration
\begin{equation} \label{arch}
    \int_{X_v^{\mathrm{an}}}{g_{n,v}\,c_1(\overline{\mathcal{D}}_n)}_v
=\int_{X_v(\mathbb{C})}{g_{n,v}\,{dd^c} g_{n,v}}=\frac{1}{\pi i}\int_{X_v(\mathbb{C})} 
\partial g_{n,v}\wedge \overline{\partial} g_{n,v} \le 0,
\end{equation}
since $\frac{1}{\pi i}\partial h\wedge \overline{\partial}h=-\frac{1}{2\pi}(h_x^2+h_y^2)dx\wedge dy \le 0$ for any smooth real function $h$ on $X_v(\mathbb{C})$.

When $v$ is a fixed non-archimedean place, for every $n$ we replace the model $\mathcal{X}_n$ over which $\mathcal{D}_n$ is defined by a 
regular semistable model dominating it after passing to a finite extension of $K$. The reduced special fiber of $\mathcal{X}_n$ can be written as
$$\mathcal{X}_{n,s}=\sum_{i=1}^r{E_i}.$$
And write
$\mathcal{D}_n=\sum_{i=1}^r{b_iE_i}$, which is a vertical $\mathbb{Q}$-Cartier divisor.

Now let $\Gamma_n$ be the reduction graph of $\mathcal{X}_{n}$, which can be embedded as a skeleton of $X_{v}^{\mathrm{an}}$ 
with respect to the semistable vertex set consisting of preimages $\xi_i$ of 
generic points of $E_i$ under the reduction map $X_{v}^{\mathrm{an}}\rightarrow \mathcal{X}_{n,s}$; see \cite{BPR13,Nic14}. 
We assign length $\lambda_v:=-\log|\varpi_v|_v $ to each edge, 
where $\varpi_v$ is a uniformizer of $K_v$. The model function $g_{n,v}$ satisfies 
$$g_{n,v}(\xi_i)=b_i \lambda_v,$$ 
and $g_{n,v}$ is affine on every edge of $\Gamma_n$. Thus $|g_{n,v}^\prime|=|b_i-b_j|$ on the edge connecting $\xi_i$ and $\xi_j$. We have
\begin{equation} \label{vert-intersect}
    \int_{\Gamma_n}{|g_{n,v}^\prime|^2dt}=\lambda_v\sum_{i<j}m_{ij}(b_i-b_j)^2.
\end{equation}
Here $m_{ij}=E_i\cdot E_j$. We claim that the right-hand side of \eqref{vert-intersect} is 
$-\lambda_v\mathcal D_n^2$, which is explained as follows. Since $\mathcal{X}_{n,s}$
is a principal divisor, we have
$$0=\mathcal{X}_{n,s} \cdot E_i=E_i^2+\sum_{i\neq j}{E_i \cdot E_j}=E_i^2+\sum_{j\ne i}m_{ij}. $$
Consequently, the self-intersection number of $\mathcal{D}_n$ is
\begin{align}
\mathcal{D}_n^2
&=\sum_i b_i^2E_i^2
  +2\sum_{i<j}b_ib_jm_{ij} \notag\\
&=-\sum_{i<j}m_{ij}(b_i-b_j)^2
.\label{self-intersection}
\end{align}

The model function $g_{n,v}$ on $X_v$ induced by a vertical 
divisor is affine on each fiber of the contraction map  
$\tau_n: X_v^{\mathrm{an}} \rightarrow \Gamma_n$ by \cite[Theorem 5.15]{BPR13}; see also \cite{Thu05}. 
Hence ${dd^c}{g_{n,v}}$ vanishes on $X_v^{\mathrm{an}}\setminus \Gamma_n$. We see that
\begin{equation}
    \int_{X_v^{\mathrm{an}}}{g_{n,v}\,{dd^c} g_{n,v}}=-\int_{X_v^{\mathrm{an}}}{|g_{n,v}^\prime|^2}dt=-\int_{\Gamma_n}{|g_{n,v}^\prime|^2}dt.
\end{equation}
We now conclude that every term of the right-hand side of \eqref{nonpositive-decomposition} 
$${\int_{X_v^{\mathrm{an}}}{g_{n,v}\,c_1(\overline{\mathcal{D}}_n)}_v}  
$$
is nonpositive. Since the left-hand side of \eqref{nonpositive-decomposition} satisfies  $\overline{\mathcal{D}}_n^2\rightarrow\overline{\mathcal{D}}^2=0$ by assumption, we get 
$$ {\int_{X_v^{\mathrm{an}}}{g_{n,v}\,c_1(\overline{\mathcal{D}}_n)}_v} \rightarrow 0,~\text{as $n\to \infty$ for every place }v\in M_K.$$
Noticing that $g_{n,v}\rightarrow g_v$ locally uniformly on $C_v^{\mathrm{an}}$, we will show that $g_v$ is a constant function, which is 
treated in the following two lemmas for archimedean and non-archimedean cases separately. 
Once we deduce that the local Green functions associated with $\overline{M}$ on $C_v^{\mathrm{an}}$ are constant for any non-archimedean and 
archimedean places of $K$, then the vertical $\overline{M}$ is a pullback from $\Pic(K/\mathbb{Z})_\mathbb{Q}$, which proves Proposition \ref{curve-vertical}. 
\end{proof}

\begin{lemma}[Archimedean case]\label{lem:zero-energy}
Let $U$ be a connected open subset of a compact Riemann surface $S$.  Let $(g_n)_{n\ge 1}$ be a sequence of smooth
real functions on $S$ that converges uniformly on every compact
subset of $U$ to a continuous function $g$ on $U$.  If
\[
\int_S {\partial g_{n}\wedge \overline{\partial} g_{n}}\longrightarrow 0,
\]
then $g$ is constant on $U$.
\end{lemma}
\begin{proof} We mainly apply Poincar\'e's inequality, which states that on a bounded connected smooth domain $B\subset \mathbb{R}^2$, 
\begin{equation*}
  \norm{u-u_B}_{L^2(B)}
  \le C_B\norm{\nabla u}_{L^2(B)}
\end{equation*}
for any smooth function $u$ on $B$, where $u_B=\frac1{|B|}\int_Bu$ and $C_B$ is a positive constant.
Now we take a coordinate chart $B\subset U$ which is biholomorphic to an open disk in $\mathbb{C}$.  The $L^2$-norm of the gradient of $u$ used here is 
\begin{equation*}\begin{aligned}
  \norm{\nabla u}_{L^2(B)}
  =&\left(\int_B |du|_g^2\,dV_g\right)^{1/2}\\
  =&\left(\int_B {(u_x^2+u_y^2})dx\wedge dy\right)^{\frac{1}{2}} \\
  =& \left(2i\int_B {\partial u\wedge \overline{\partial}}u\right)^{\frac{1}{2}}.
\end{aligned}
\end{equation*}
Since $g_n \rightarrow g$ locally uniformly, we have the $L^2(B)$-convergence and average convergence:
\[
  g_n\to g\text{ in }L^2(B),~\,\text{and}
  \quad (g_n)_B\to g_B.
\]
Poincar\'e's inequality implies that
\[
  \norm{g_n-(g_n)_B}_{L^2(B)}
  \le C_B\norm{\nabla g_n}_{L^2(B)}\longrightarrow 0.
\]
Letting $n\rightarrow \infty$, we obtain
\[
  \norm{g-g_B}_{L^2(B)}=0.
\]
Thus, $g=g_B$ almost everywhere on $B$. Since $g$ is continuous and $U$ is connected, we deduce that $g$ is a constant on $U$.
\end{proof}

\begin{lemma}[Non-archimedean case]\label{nonarch}
Let $U$ be a connected open subset of a projective Berkovich curve $X$ which is the analytification of a smooth quasi-projective curve.  Let $\{g_n\}$ be a sequence of
model functions on $X$ that converges uniformly on every compact
subset of $U$ to a continuous function $g$.  If
\[
\int_{X}{|g_n'|^2\,dt}\longrightarrow 0,~\text{as}~ n\to \infty,
\]
then $g$ is constant on $U$.
\end{lemma}

\begin{proof}
We assume that $X$ is defined over a complete algebraically closed field over $K_v$. Let $x,y\in U$ be points of type $2$ or $3$.  Since a connected open Berkovich curve
is path connected, there is a piecewise geodesic segment
\[
\gamma\subset U
\]
connecting $x$ to $y$ with length $L(\gamma)< \infty$. For each $n$, we can choose a skeleton $\Gamma_n$ on which $g_n$ is supported and $\Gamma_n$ is large enough to 
contain $\gamma$. The restriction of
$g_n$ to $\gamma$ is piecewise affine. 
Parametrizing $\gamma$ by arc length and applying the Cauchy--Schwarz inequality, we obtain
\begin{equation}
\begin{aligned}
|g_n(x)-g_n(y)|^2
&=
\left|\int_{\gamma}g_n'(t)\,dt\right|^2  \\
&\le
\left(\int_{\gamma}1\,dt\right)\left(\int_{\gamma}|g_n'(t)|^2\,dt\right)\\
&\le
L(\gamma) \int_{X}{|g_n'(t)|^2\,dt} 
\longrightarrow 0.
\end{aligned}
\end{equation}
Since we have $g_n(x)\rightarrow g(x)$ and $g_n(y)\rightarrow g(y)$, we get $g(x)=g(y)$ by passing 
to the limit. Moreover, we know that points of type $2$ and $3$ are 
dense in a Berkovich analytic curve, and $g$ is continuous. Consequently, the function $g$ is a constant on $U$.
\end{proof}

\subsection{General case} \label{curve-K}
In this subsection, we will reduce the proof of Theorem \ref{curve-vert-K} to the vertical case.
Given an adelic line bundle $\overline{M}$ on $C$ such that 
$\deg(\widetilde{M})=0$ and $\overline{M}^2=0$, we will show that we can modify the Cauchy representatives of $\overline{M}$ so that indeed $\overline{M}$ is vertical.
\begin{proof}[Proof of Theorem \ref{curve-vert-K}]
For convenience, we represent $\overline{M}$ by an adelic $\mathbb{Q}$-Cartier divisor $\overline{D}$. Let
$X$ be a smooth projective model of $C$, and write 
$$ X\setminus C=\{p_1,\dots,p_r\}, $$
which is a finite set. We assume that $X\setminus C$ is non-empty since the projective case has already been solved by Yuan--Zhang. Choose a Cauchy sequence of model adelic divisors 
$$ \{\overline{\mathcal{D}}_n \}_{n\geq 1} $$
representing $\overline{D}$. We may assume that 
$\overline{\mathcal{D}}_n$ is defined on an inverse system of projective 
arithmetic models $\mathcal{X}_n$ such that all their generic fibers are identified with $X$. Let $D_n:=\mathcal{D}_{n,K}$. We have 
 $$D_n|_C =D|_C, \mathrm{~for~every~}n. $$ 
Now we fix a strictly effective boundary divisor $\overline{\mathcal{E}}$ whose generic part is 
$$E:=\mathcal{E}_K=\sum_{1\leq i\leq r}{e_i}p_i,$$
with $e_i>0$
for all $i$. By definition, we have a sequence of positive rational numbers $\varepsilon_n\rightarrow 0$ as $n\to \infty$ such that
$$ -\varepsilon_nE \leq D_m-D_n \leq \varepsilon_nE $$
for all $m\geq n$.

By the assumption that $ \deg(\widetilde{D})=0 $ on $C$, we can modify the representative $\{\overline{\mathcal{D}}_n\}$ so that $\deg({D}_n)=0$, for all $n\ge 0$. Since
$\deg(D_n)\to\deg(\widetilde{D})=0$, we can set 
\begin{equation}\label{modify}
 c_n := -\frac{\deg_X(D_n)}{\deg_X(E)},\quad D_n^\prime:=D_n+c_nE.
\end{equation} 
Then $c_n\to 0$ as $n\to \infty$, and $D_n^\prime$ has degree zero. 
Let $\overline{\mathcal{E}}_n$ be the pullback of $\overline{\mathcal{E}}$ to the $n$-th model, and 
replace $\overline{\mathcal D}_n$ by 
\[ \overline{\mathcal{D}}_n^\prime := \overline{\mathcal{D}}_n+c_n\overline{\mathcal{E}}_n. \] 
Since $\overline{\mathcal{E}}_n$ is a boundary divisor, the sequence $\{\overline{\mathcal{D}}_n^\prime\}$ 
converges to the same limit $\overline{\mathcal{D}}$ as $\{\overline{\mathcal{D}}_n\}$ does. The generic divisors have 
degree zero and still converge coefficient-wise to $\widetilde{D}$. The generic divisor of $\overline D_n^\prime$ has degree zero for all $n$. Consequently, we may assume that \[ \deg_X(D_n)=0 \] for every $n$.

 Next we will show that $ \widetilde{D}\equiv_{\mathrm{lin}} 0$ in $\operatorname{Div}^0(X)_{\mathbb{R}} $.
 We can take $\overline{\mathcal{P}}_n \in \Pic(X)_{\mathbb{Q},\mathrm{int}}$
to be an admissible extension of $\mathbb{Q}$-line bundle $\mathcal{O}_X(D_n)$ such that $\overline{\mathcal{V}}_n\colon=\overline{\mathcal{D}}_n-\overline{\mathcal{P}}_n$ is vertical and
$$ \overline{\mathcal{P}}_n \cdot \overline{\mathcal{V}}_n=0, $$
by \cite[Lemma 2.6]{YZ17}. Then we have $\overline{\mathcal{V}}_n^2 \leq 0$ by the arithmetic 
Hodge index theorem over projective varieties. Hence
$$ \overline{\mathcal{D}}_n^2=(\overline{\mathcal{P}}_n+\overline{\mathcal{V}}_n)^2=\overline{\mathcal{P}}_n^2+\overline{\mathcal{V}}_n^2=-2\hat{h}(D_n)+\overline{\mathcal{V}}_n^2 \leq -2\hat{h}(D_n) \leq 0. $$
Here $\hat{h}:\mathrm{Pic}^0(X_{\overline{K}})\rightarrow \mathbb{R}$ is the N\'eron--Tate 
height function introduced in \cite[$\S$6.5.1]{YZ26}. Since $\{\overline{\mathcal{D}}_n\}$ represents 
$\overline{D}$, the intersection pairing gives:
$$\lim_{n\rightarrow \infty}{\overline{\mathcal{D}}_n^2}=\overline{D}^2=0.$$
It follows that
$$ \lim_{n\rightarrow \infty}\hat{h}(D_n)=0 $$
Since $\hat{h}$ is continuous on $J(K)_{\mathbb{R}}$, where $J=\mathrm{Pic}^0_{X/K}$ is the Jacobian of $X$, we have
$$ \hat{h}(\widetilde{D})=0. $$
By the positive definiteness of $\hat{h}$, we see that
$\widetilde{D}\equiv_{\mathrm{lin}} 0 \mathrm{~in~} {\mathrm{Pic}^0(X)}_\mathbb{R}$, 
and $D|_C \equiv_{\mathrm{lin}} 0$ in ${\mathrm{Pic}^0(C)}_\mathbb{Q}$.

Now we aim to show that $\overline{\mathcal{D}}$ can 
be represented by a Cauchy sequence $\{\overline{\mathcal{D}}_n\}$ such that the generic 
part $D_n$ is exactly zero for every $n$. Since $\widetilde{D}\equiv_{\mathrm{lin}}0$ in ${\mathrm{Pic}^0(X)}_\mathbb{R}$ and $D|_C\equiv 0$ in ${\mathrm{Pic}^0(C)}_\mathbb{Q}$,
there exist rational functions $u_1,\ldots,u_t\in \Gamma(C,\mathcal{O}_C^{\times})$
and real numbers $\lambda_1,\ldots,\lambda_t$ such that
\[
    \widetilde{D}
    =
    \sum_{j=1}^t
    \lambda_j\operatorname{div}_X(u_j).
\]
Choose sequences of rational numbers $\lambda_{j,n}\rightarrow\lambda_j$ as $n\to \infty$ for all $j$,
and replace $\{\overline{\mathcal{D}}_n\}$ by
\[\{\overline{\mathcal{D}}_n-\sum_{j=1}^t
    \lambda_{j,n}\widehat{\operatorname{div}}(u_j)\},
\]
which represents the same adelic line
bundle as $\{\overline{\mathcal{D}}_n\}$, and its generic divisor is supported on $X \setminus C=\{p_1,\dots,p_r\}$. Therefore, there exists a sequence of rational numbers $\{r_{i,n}\}_{1\leq i\leq r,\,n\geq1}$ such that
\[
D_n=\mathcal{D}_{n,K}=\sum_{1\leq i\leq r} r_{i,n}p_i,
\quad \mathrm{with~} \max_i|r_{i,n}|\longrightarrow 0,~\text{as}~n \to \infty.
\]
Thus we can replace $\{\mathcal{O}( \overline{\mathcal{D}}_n)\}$ by 
$\{\mathcal{O}( \overline{\mathcal{D}}_n)\otimes \mathcal{O}(-\sum_{1\leq i\leq r} r_{i,n}\overline{p}_i)\}$, where $\overline{p}_i$ is the Zariski closure of 
$p_i$ in the arithmetic model $\mathcal{X}_n$. The geometric part $\mathcal D_{n,K}$ is exactly $0$ 
for all $n$ and we can assume now $\overline{\mathcal{D}}$ is a 
vertical adelic divisor, which can be approximated by vertical 
model adelic divisors. Hence the theorem can be reduced to the vertical case, which is proved in Proposition \ref{curve-vertical}. 
\end{proof}

\section{Curves over finitely generated fields}\label{curve over F}
 \subsection{Adelic line bundles with numerically trivial heights}\label{Q1r}
Suppose that we are given an integrable adelic line bundle $\overline{M}$ on a normal, 
geometrically integral quasi-projective curve $C$ over a finitely generated field $F/\mathbb{Q}$. In this section, we show that if $\mathfrak{h}_{\overline{M}}(x)\equiv_{\mathrm{num}} 0$ for all $x\in C(\overline{F})$, then $\deg(\widetilde{M})=0$ and
$$ \overline{M}^2 \cdot \pi_C^*\overline{H}_1 \cdots \pi_C^*\overline{H}_r =0,
$$
for $r$ arbitrary integrable adelic line bundles $\overline{H}_1,\dots,\overline{H}_r$. 
This exactly says that the hypothesis of assertion $H(1,r)$ implies the hypothesis of assertion $R(1,r)$. Hence we can reduce proving $H(1,r)$ to proving $R(1,r)$, while $R(1,r)$ will be shown later in $\S$\ref{explicit}.
\begin{theorem}[$H(1,r)$ holds]\label{curve:fin-gen}
  Let $F$ be a finitely generated 
field of transcendence degree $r\geq 1$ over $\mathbb{Q}$, and let $C$ be a normal geometrically integral quasi-projective curve over 
$F$. If $ \overline{M} \in \widehat{\mathrm{Pic}}(C/\mathbb{Z})_{\mathrm{int},\mathbb{Q}} $ satisfies that
$$
\mathfrak{h}_{\overline{M}}(x)\equiv_{\mathrm{num}} 0 \in \widehat{\mathrm{Pic}}(F/\mathbb{Z})_{\mathbb{Q}},
$$
for every $x\in C(\overline{F})$,
then we have 
$ \overline{M}\in\pi_{C}^{*}\widehat{\mathrm{Pic}}(F/\mathbb{Z})_{\mathbb{Q}}. $
\end{theorem}
For such $\overline{M}$ with numerically trivial height at every algebraic point, we first show that the geometric part $\widetilde{M}$ is numerically trivial.
To see this, choose a finitely generated field $F_0\subseteq F$ such that
$$
\mathrm{tr.deg}_{\mathbb{Q}}F_0=r-1,\quad\mathrm{tr.deg}_{F_0}F=1,
$$
and $F_0$ is algebraically closed in $F$. Hence we can take a quasi-projective model $\mathcal{C}\rightarrow 
S$ of $C\rightarrow \Spec F$, where $S$ is an integral quasi-projective
curve over $F_0$ with function field identified with $F$.
Let $\rho:\Spec F \rightarrow \Spec F_0$ be the structure morphism. Then we can utilize the 
results in \cite[$\S$5.3.7]{YZ26} dealing with the function field analogue 
of essential minimum and pseudo-effectivity; see also \cite{QY24}.

There is a canonical varying base map  
$$ v_C:\widehat{\mathrm{Pic}}(C/\mathbb{Z})_{\mathrm{int},\mathbb{Q}} \longrightarrow \widehat{\mathrm{Pic}}(C/F_0)_{\mathrm{int},\mathbb{Q}} $$
introduced in \cite[$\S$2.5.5]{YZ26}. We denote the image of 
$\overline{M}$ by $\overline{M}^{\mathrm{rel}}$, which only carries the 
adelic information over $F_0$. For $z\in C(\overline{F})$, let 
$z^\prime$ be the corresponding closed point and $E:= \kappa(z^\prime)$ be its residue field. Then the vector-valued 
height of $z$ defined by $\overline{M}^{\mathrm{rel}}$  is
$$ 
\mathfrak{h}_{\overline{M}^{\mathrm{rel}}}(z):=\frac{1}{\left[E:F\right]}{\langle \overline{M}^{\mathrm{rel}}|_{z^\prime} \rangle} \in \widehat{\mathrm{Pic}}(F/F_0)_{\mathrm{int},\mathbb{Q}}.
$$
We can choose a representative of $\mathfrak{h}_{\overline{M}^{\mathrm{rel}}}(z)$ on a quasi-projective model
of $F/F_0$. Applying the degree map to $\mathfrak{h}_{\overline{M}^{\mathrm{rel}}}(z)$, 
we obtain a real-valued height  
$$
h_{\overline{M}^{\mathrm{rel}}}(z):=\widehat{\deg}_{F/F_0}( \mathfrak{h}_{\overline{M}^{\mathrm{rel}}}(z)).
$$
Then we will show that $\mathfrak{h}_{\overline{M}}(x)\equiv_{\mathrm{num}} 0$ for 
every algebraic point $x$ implies the two hypotheses of $R(1,r)$ in the following two propositions.

\begin{proposition} \label{geom-trivial}
    Let $\overline{M}$ be an adelic line bundle on $C$ as above, which satisfies
    $$\mathfrak{h}_{\overline{M}}(z)\equiv_{\mathrm{num}} 0 \in \widehat{\mathrm{Pic}}(F/\mathbb{Z})_{\mathbb{Q}}$$
    for every $z\in C(\overline{F})$. Suppose $F_0\subseteq F$ is a subfield with $\mathrm{tr.deg}_{\mathbb{Q}}F_0=r-1$, 
    $\mathrm{tr.deg}_{F_0}F=1$, and $F_0$ is algebraically closed in $F$. Let $\overline{M}^{\mathrm{rel}}$ be the image of 
    $\overline{M}$ in $\widehat{\mathrm{Pic}}(C/F_0)_{\mathrm{int},\mathbb{Q}}$. Then the geometric height
    $$ h_{\overline{M}^{\mathrm{rel}}}: C(\overline{F})\longrightarrow \mathbb{R} $$
     induced by $\overline{M}^{\mathrm{rel}}$ is identically zero. In particular, we have 
     $
     \deg(\widetilde{M})=0.
     $
\end{proposition}
\begin{proof} Let $z^\prime$ be the associated closed point of $z$ in $C$ and $E:=F(z^\prime)$ be its residue field.
    We denote the varying base map on $\Spec E$ and $\Spec F$ by $v_E$ and $v_F$ respectively. The horizontal arrows in the 
    following diagram are the norm functors induced by the finite 
    extension $E/F$ as the relative $0$-dimensional Deligne pairing, which is compatible with base change:
    \[
\begin{tikzcd}
\widehat{\mathrm{Pic}}(E/\mathbb{Z})_{\mathrm{int},\mathbb{Q}} \arrow{r}{N_{E/F}} \arrow[swap]{d}{v_E} & \widehat{\mathrm{Pic}}(F/\mathbb{Z})_{\mathrm{int},\mathbb{Q}} \arrow{d}{v_F} & \overline{M}|_{z^\prime}\arrow[|->]{r}\arrow[|->]{d}  & \langle\overline{M}|_{z^\prime}\rangle \arrow[|->]{d}\\
\widehat{\mathrm{Pic}}(E/F_0)_{\mathrm{int},\mathbb{Q}} \arrow{r}{N_{E/F}^{\mathrm{rel}}}& \widehat{\mathrm{Pic}}(F/F_0)_{\mathrm{int},\mathbb{Q}}, & \overline{M}^{\mathrm{rel}}|_{z^\prime}\arrow[|->]{r}  & \langle\overline{M}^{\mathrm{rel}}|_{z^\prime} \rangle
\end{tikzcd}
\]
Therefore, by definition we have 
\begin{equation}\label{varying-base}
\mathfrak{h}_{\overline{M}^{\mathrm{rel}}}(z)=v_F(\mathfrak{h}_{\overline{M}}(z)).
\end{equation}
Now we observe that for arbitrary integrable adelic $\mathbb{Q}$-line bundles 
$\overline{H}_1,\dots,\overline{H}_r\in \widehat{\mathrm{Pic}}(F_0/\mathbb{Z})_{\mathrm{int},\mathbb{Q}}$ and every $z\in C(\overline{F})$, we have
\begin{equation}\label{h=h}
\mathfrak{h}_{\overline{M}}(z)\cdot \rho^*\overline{H}_1\cdots\rho^*\overline{H}_r =h_{\overline{M}^{\mathrm{rel}}}(z) \cdot(\overline{H}_1\cdots\overline{H}_r).
\end{equation} 
This is obtained by representing $\mathfrak{h}_{\overline{M}}(z)$ by an adelic $\mathbb{Q}$-line bundle $\overline{N}_z$ on a 
quasi-projective model $q: \mathcal{B} \rightarrow \mathcal{T}$ over $\mathbb{Z}$ of $\Spec F\rightarrow \Spec F_0$ such 
that $C\rightarrow \Spec F$ is the generic fiber of 
$\mathcal{C}\rightarrow \mathcal{B}$; such models exist by \cite[Lemma 2.3.3]{YZ26}. Suppose $\mathcal{T}\rightarrow \mathcal{S}$ is flat with base $\mathcal{S}=\Spec O_K$ for some number field $K$. Then we apply \cite[Lemma 4.6.1(3)]{YZ26} to
$\mathcal{B}\rightarrow \mathcal{T} \rightarrow \mathcal{S}$ to get:
$$ \pi_{\mathcal{B}/\mathcal{S},*}\langle\overline{N}_z, \rho^*\overline{H}_1,\dots,\rho^*\overline{H}_r\rangle = (\widehat{\deg}_{F/F_0}(\overline{N}_z^{\mathrm{rel}}))\cdot\pi_{\mathcal{T}/\mathcal{S},*}\langle\overline{H}_1,\cdots,\overline{H}_r\rangle \in \widehat{\mathrm{Pic}}(\mathcal{S}/\mathbb{Z})_{\mathrm{int},\mathbb{Q}}. $$
Here $\overline{N}_z^{\mathrm{rel}}:= v_F(\overline{N}_z)$ is the 
image of $\overline{N}_z$ under the varying base map, which represents 
$\mathfrak{h}_{\overline{M}^{\mathrm{rel}}}(z)$ by \eqref{varying-base}. Therefore, taking the arithmetic degree of both sides on 
$\mathcal{S}$ gives exactly \eqref{h=h}. Hence the assumption 
$\mathfrak{h}_{\overline{M}}(z)\equiv_{\mathrm{num}} 0 $ and \eqref{h=h} together imply that 
$h_{\overline{M}^{\mathrm{rel}}}(z)=0$ for every $z\in C(\overline{F})$  
after choosing nef $\overline{H}_i$ on $\Spec F_0$ with $\overline{H}_1\cdots\overline{H}_r>0$.

 Recall that the essential minimum of $C\rightarrow \Spec F$, where $F$ is 
 considered as a function field of one variable over $F_0$, is defined by:
$$
e_1(C,\overline{M}^{\mathrm{rel}})= \sup_{U \subset C}\inf_{x\in U(\overline{F})}{h_{\overline{M}^{\mathrm{rel}}}(x)},
$$
where the supremum is taken over all non-empty Zariski open subschemes $U$ of $C$. Therefore, we have 
$$ e_1(C,\overline{M}^{\mathrm{rel}})=0=e_1(C,-\overline{M}^{\mathrm{rel}}). $$
Therefore, both $\overline{M}^{\mathrm{rel}}$ and 
$-\overline{M}^{\mathrm{rel}}$ are pseudo-effective in $\widehat{\mathrm{Pic}}(C/F_0)_{\mathrm{int},\mathbb{Q}}$ 
by the result for function fields in \cite[Theorem 1.8]{QY24};
see also \cite[Theorem 5.3.6]{YZ26}.
Pseudo-effectivity of $\pm \overline{M}^{\mathrm{rel}}$ implies
that for any nef adelic line bundle $ \overline{H}\in\widehat{\mathrm{Pic}}(F/F_0)_{\mathbb{Q},\mathrm{nef}}$, 
$$ 
\overline{M}^{\mathrm{rel}} \cdot \pi_C^*\overline{H}\ge 0 \text{ and } -\overline{M}^{\mathrm{rel}} \cdot \pi_C^*\overline{H}\ge 0.
$$
Therefore, we have 
$\overline{M}^{\mathrm{rel}} \cdot \pi_C^*\overline{H}= 0$. Moreover, we have 
$$
\overline{M}^{\mathrm{rel}} \cdot \pi_C^*\overline{H}= \deg(\widetilde{M})\cdot \widehat{\deg}_{F/F_0}(\overline{H}), 
$$
from the projection formula. Since $\overline{H}$ is arbitrary, after choosing one $\overline{H}$ of positive degree, we have
$ 
\deg(\widetilde{M})=0.
$
\end{proof}

We have the following proposition, which reduces proving $H(1,r)$ to proving $R(1,r)$.
\begin{proposition}\label{Q-implies-H}
    Let $F$ be a finitely generated 
field of transcendence degree $r\geq 1$ over $\mathbb{Q}$, and $C$ be 
a normal geometrically integral quasi-projective curve over 
$F$, with structure morphism $\pi_C:C\to\Spec F$. 
Suppose $ \overline{M} \in \widehat{\mathrm{Pic}}(C/\mathbb{Z})_{\mathrm{int},\mathbb{Q}} $ is an integrable adelic 
line bundle on $C$ satisfying that the vector-valued height 
$$
\mathfrak{h}_{\overline{M}}(x)\equiv_{\mathrm{num}} 0 \in \widehat{\mathrm{Pic}}(F/\mathbb{Z})_{\mathbb{Q}}.
$$
for every $x\in C(\overline{F})$.
Then for any integrable adelic line bundle $\overline{A}\in\widehat{\mathrm{Pic}}(C/\mathbb{Z})_{\mathrm{int},\mathbb{Q}}$ and $r$ 
integrable $\overline{H}_1,\dots,\overline{H}_r \in \widehat{\mathrm{Pic}}(F/\mathbb{Z})_{\mathrm{int},\mathbb{Q}}$, 
we have 
$$ \overline{M}\cdot\overline{A} \cdot \pi_C^*\overline{H}_1 \cdots \pi_C^*\overline{H}_r =0.$$
In particular, we have 
$$\overline{M}^2 \cdot \pi_C^*\overline{H}_1 \cdots \pi_C^*\overline{H}_r =0.$$
\end{proposition}
\begin{proof}
First, 
we choose a normal, flat, integral projective arithmetic model $\mathcal{B}$ of $F$ over $\mathbb{Z}$ such that there is a projective model
$$
\pi_0: \mathcal{X}_0 \longrightarrow \mathcal{B},
$$
whose generic fiber is $ X\rightarrow \Spec F$ and $X$ is a projective model 
of $C$ over $F$. We can assume $X$ is smooth by resolution of singularities \cite{Lip78}. 
We can then take a cofinal family of normal, generically smooth projective models
$$
p_n: \mathcal{X}_n \longrightarrow \mathcal{X}_0,
$$
with $\pi_n=\pi_0 \circ p_n: \mathcal{X}_n\rightarrow \mathcal{B}$ being a projective arithmetic model of $X\rightarrow \Spec F$.

 Since  
$\overline{A},\overline{H}_1,\dots \overline{H}_r$ can be approximated by model adelic line bundles and integrable adelic line bundles are differences of strongly nef ones, it suffices to prove for any nef 
hermitian $\mathbb{Q}$-line bundle $\overline{\mathcal{A}}\in \widehat{\mathrm{Pic}}(\mathcal{X}_0)_{\mathbb{Q}}$ and nef 
$\overline{\mathcal{H}}_1,\dots \overline{\mathcal{H}}_r \in \widehat{\mathrm{Pic}}(\mathcal{B})_{\mathbb{Q}}$. Let

$$ 
\{\overline{\mathcal{M}}_n\}_{n \ge 0}, \qquad \overline{\mathcal{M}}_n\in \widehat{\mathrm{Pic}}(\mathcal{X}_n)_{\mathbb{Q}}
$$
be a Cauchy sequence of hermitian $\mathbb{Q}$-line bundles representing 
$\overline{M}$. Here the limit is taken with respect to  
$\overline{\mathcal{E}}_0$ which is a strictly effective boundary divisor 
on $\mathcal{X}_0$, whose support contains $\mathcal{X}_0 \setminus \mathcal{C}$, where $\mathcal C \subset \mathcal{X}_0$ is a quasi-projective arithmetic model of $C$.
Its smooth metric is chosen so that the norm of its canonical 
section is at most $1$. Write $\overline{\mathcal{E}}_n:=p_n^*\overline{\mathcal{E}}_0$. The Cauchy condition reads as
\begin{equation}\label{Cauchy}
-\varepsilon_n \overline{\mathcal{E}}_n \le \overline{\mathcal{M}}_n- \overline{\mathcal{M}}_m \le \varepsilon_n \overline{\mathcal{E}}_n, \quad \text{where }\varepsilon_n \rightarrow 0 \text{ as }n\to\infty, 
\end{equation}
for all $m> n\ge 0$. From Proposition \ref{geom-trivial}, we have $\deg(\widetilde{M})=0$. Then we can normalize each $\overline{\mathcal{M}}_n$ similarly as in \eqref{modify} and assume that 
\begin{equation}\label{degree-zero}
    \deg(\mathcal{M}_{n,F})=0.
\end{equation}

Now fix $n\geq 0$. By openness of the flat map $\mathcal{X}_n\to \mathcal{B}$ and geometric integrality, we can take a dense open subset
$$
B_n^{\circ} \subset \mathcal{B}_{\mathbb{Q}},
$$
such that $\pi_{n,\mathbb{Q}}$ is flat with geometrically integral 
fibers. Since $\deg_X(\mathcal M_{n,F}) =0$, we can shrink it further so 
that the degree of $\mathcal M_{n,F}$ on every fiber over $B_n^{\circ}$ is zero.
To apply the arithmetic Demailly approximation theorem of Qu--Yin in \cite{QY24}, we need to take a closed subset $Y_n$ of 
$\mathcal{X}_{n,\mathbb{Q}}$ and one ample line bundle for every 
$\mathcal{X}_{n}$. We 
define $Y_n \subset \mathcal{X}_{n,\mathbb{Q}}$ to be the finite union 
of the following three kinds of prime divisors by the dimensions of its 
images under map $\mathcal{X}_{n,\mathbb{Q}} \rightarrow \mathcal{B}_{\mathbb{Q}}$: 
\begin{enumerate}[1.]
    \item  every prime divisor contained in $\supp(\mathcal{E}_{n,\mathbb{Q}})$;
    \item every divisorial component whose image in $\mathcal{B}_{\mathbb{Q}}$ is a codimension one component of $\mathcal{B}_{\mathbb{Q}}\setminus B_n^{\circ}$;
    \item every prime divisor whose image in $\mathcal{B}_{\mathbb{Q}}$ has codimension at least two.
\end{enumerate}
Note that if a prime divisor $Z$ of $\mathcal{X}_{n,\mathbb{Q}}$ has image of codimension at least two, then the generic point of $Z$ over its 
image has relative dimension of at least two. The fiber dimension function is known to be upper-semicontinuous. Thus, $Z$ is contained in the inverse image of a 
proper closed subset of $\mathcal{B}_{\mathbb{Q}}$, and    $\mathcal{X}_{n,\mathbb{Q}}$ contains only finitely many such components by noetherian property.

As the pullback of $\mathcal{A}$ under birational morphism is nef, not ample in general, we choose for every $n\ge 0$ an ample 
hermitian line bundle $\overline{\mathcal{B}}_n$ on $\mathcal{X}_n$. We define
$$
\overline{\mathcal{A}}_{n,\delta} :=p_n^*\overline{\mathcal{A}}+\delta\cdot \overline{\mathcal{B}}_{n},\quad \overline{\mathcal{A}}_{n}= p_n^*\overline{\mathcal{A}},
$$
for every rational number $\delta > 0$. Then $\overline{\mathcal{A}}_{n,\delta}$ is ample.

Fix one rational $\delta > 0$. For every $n\ge 0$, we apply the arithmetic Demailly approximation in \cite[Theorem 4.10]{QY24} to
$$
(\mathcal{X}_n, \mathcal{A}_{n,\delta},Y_n ),
$$ 
and then we get integers $N_{n,k} \rightarrow \infty$ and a sequence of sections $s_{n,k} \in H^0 (\mathcal{X}_n, N_{n,k}\mathcal{A}_{n,\delta} )$ such that $N_{n,k}\mathcal{A}_{n,\delta}$ is integral and

\begin{enumerate}[$\bullet$]
    \item  $\frac{1}{N_{n,k}}\log\|s_{n,k}\|\rightarrow 0$ in $L^1(\mathcal{X}_n(\mathbb{C}))$;
    \item  $\div(s_{n,k})_\mathbb{Q}$ is smooth and does not contain any irreducible components of $Y_n$;
    \item  $\div(s_{n,k})$ has no vertical components over $\Spec \mathbb Z$.
\end{enumerate}
Therefore, we can write
\begin{equation} \label{div(s)}
\div(s_{n,k})= 
 \sum_{\alpha}c_{\alpha,n,k}\mathcal{Z}_{\alpha,n,k}
 +\sum_{\beta}d_{\beta,n,k}\mathcal{W}_{\beta,n,k},
\end{equation}
where $\mathcal{Z}_{\alpha,n,k}$ is the prime divisor that dominates $\mathcal{B}$ and $\mathcal{W}_{\beta,n,k}$ maps to a horizontal prime divisor of $\mathcal{B}$. Now we consider the intersection product
\begin{equation*}
I_{n,\delta} := \overline{\mathcal{M}}_n \cdot \overline{\mathcal{A}}_{n,\delta}\cdot \pi_n^*\overline{\mathcal{H}}_1\cdots \pi_n^*\overline{\mathcal{H}}_r.
\end{equation*}
Using the inductive formula by taking above section $s_{n,k}$ of $N_{n,k} \mathcal{A}_{n,\delta}$, we have
\begin{equation} 
\begin{aligned}
   I_{n,\delta}=&\frac{1}{N_{n,k}} (\overline{\mathcal{M}}_n \cdot  \pi_n^*\overline{\mathcal{H}}_1\cdots \pi_n^*\overline{\mathcal{H}}_r|\div(s_{n,k})) 
    -\frac{1}{N_{n,k}}\int_{\mathcal{X}_n(\mathbb{C})}{\log\|s_{n,k}\|}\,\mu_n\\
   =&\frac{1}{N_{n,k}}\sum_{\alpha}c_{\alpha,n,k}(\overline{\mathcal{M}}_n \cdot  \pi_n^*\overline{\mathcal{H}}_1\cdots \pi_n^*\overline{\mathcal{H}}_r|\mathcal{Z}_{\alpha,n,k})\\
   &+ 
   \frac{1}{N_{n,k}}\sum_{\beta}d_{\beta,n,k}(\overline{\mathcal{M}}_n \cdot  \pi_n^*\overline{\mathcal{H}}_1\cdots \pi_n^*\overline{\mathcal{H}}_r|\mathcal{W}_{\beta,n,k}) 
   -\frac{1}{N_{n,k}}\int_{\mathcal{X}_n(\mathbb{C})}{\log\|s_{n,k}\|} \,\mu_n, \label{In} \\
  \text{where}&~\mu_n:= c_1(\overline{\mathcal M}_n)\wedge
   \pi_n^*c_1(\overline{\mathcal H}_1)\wedge\cdots\wedge
   \pi_n^*c_1(\overline{\mathcal H}_r).
\end{aligned}
\end{equation}
    Since $\mathcal{Z}_{\alpha,n,k}$ dominates $\mathcal{B}$, 
its generic point defines an algebraic point
$z_{\alpha,n,k}\in C(\overline{F})$ . We denote the corresponding closed point by $z_{\alpha,n,k}^\prime$, and let $E_{\alpha,n,k} := F(z_{\alpha,n,k}^\prime)$ be its residue field.
The choice of $Y_n$ guarantees that $z_{\alpha,n,k}^\prime$ lies in $C\subset \mathcal{X}_{n,\mathbb{Q}}$.

We define the Moriwaki height of $z\in C(\overline{F})$ with respect to $\overline{\mathcal{M}}_n$ by
\begin{equation*}
h_{\overline{\mathcal{M}}_n}^{\overline{\mathscr{H}}}(z):=\mathfrak{h}_{\overline{\mathcal{M}}_n}(z)\cdot \overline{\mathcal{H}}_1 \cdots \overline{\mathcal{H}}_r,
\end{equation*}
for $z\in C(\overline{F})$ with respect to a nef tuple 
$\overline{\mathscr{H}}=(\overline{\mathcal{H}}_1, \cdots, \overline{\mathcal{H}}_r)$ and similarly for 
$\overline{M}$. Let $m \rightarrow \infty$ in \eqref{Cauchy}, we have
$$-\varepsilon_n \overline{\mathcal{E}}_n \le \overline{\mathcal{M}}_n- \overline{M} \le \varepsilon_n \overline{\mathcal{E}}_n, \quad \varepsilon_n \rightarrow 0. $$
Taking the height of $z_{\alpha,n,k}$ with respect to each term, we get
 $$ 
 -\varepsilon_n h_{\overline{\mathcal{E}}_n}^{\overline{\mathscr{H}}}(z_{\alpha,n,k})\leq h_{\overline{\mathcal{M}}_n}^{\overline{\mathscr{H}}}(z_{\alpha,n,k})-h_{\overline{M}}^{\overline{\mathscr{H}}}(z_{\alpha,n,k}) \leq \varepsilon_n h_{\overline{\mathcal{E}}_n}^{\overline{\mathscr{H}}}(z_{\alpha,n,k}),
 $$
 by the effectivity of the boundary divisors and the nefness of all $\overline{\mathcal{H}}_i$. Since $\mathfrak{h}_{\overline{M}}(z_{\alpha,n,k})$ is numerically trivial, we have $h_{\overline{M}}^{\overline{\mathscr{H}}}(z_{\alpha,n,k})=0$ and hence
 $$
   \left|h_{\overline{\mathcal{M}}_n}^{\overline{\mathscr{H}}}(z_{\alpha,n,k})\right| \leq \varepsilon_n h_{\overline{\mathcal{E}}_n}^{\overline{\mathscr{H}}}(z_{\alpha,n,k}).
 $$

Now we discuss the intersection products on $\mathcal{Z}_{\alpha,n,k}$ and $\mathcal{W}_{\beta,n,k}$ in the following lemmas:

\begin{lemma}\label{Z}
    For a component $\mathcal{Z}=\mathcal Z_{\alpha,n,k}$, with residue field $E=E_{\alpha,n,k}$ and an algebraic point $z=z_{\alpha,n,k}\in C(\overline{F})$, we have
\begin{equation} \label{eqht}
 (\overline{\mathcal{M}}_n\cdot
 \pi_n^*\overline{\mathcal{H}}_1\cdots \pi_n^*\overline{\mathcal{H}}_r
 |\mathcal{Z})
 =[E:F]\cdot (\mathfrak{h}_{\overline{\mathcal{M}}_n}(z)\cdot  \overline{\mathcal{H}}_1 \cdots \overline{\mathcal{H}}_r)
 = [E:F]\cdot  h_{\overline{\mathcal{M}}_n}^{\overline{\mathscr{H}}}(z).
\end{equation}
The same formula holds with $\ol{\mathcal M}_n$ replaced by
$\ol{\mathcal{E}}_n$.
\end{lemma}
\begin{proof}
 Let $\alpha: \Spec E \rightarrow \Spec F$ be the morphism induced by $z$. We have
  \begin{equation*}
  \begin{aligned}
    (\overline{\mathcal{M}}_n\cdot
 \pi_n^*\overline{\mathcal{H}}_1\cdots \pi_n^*\overline{\mathcal{H}}_r
 |\mathcal{Z})= &
    \overline{\mathcal{M}}_n|_{z^\prime} \cdot \alpha^*\overline{\mathcal{H}}_1 \cdots \alpha^*\overline{\mathcal{H}}_r \\
    = & N_{E/F}(\overline{\mathcal{M}}_n|_{z^\prime}) \cdot \overline{\mathcal{H}}_1 \cdots \overline{\mathcal{H}}_r,    
  \end{aligned}
  \end{equation*}
  from the projection formula of relative dimension $0$.  Recall that 
  $\mathfrak{h}_{\overline{\mathcal{M}}_n}(z)=\frac{1}{\left[E:F\right]}{N_{E/F}(\overline{\mathcal{M}}_n|_{z^\prime})}$, so this gives the equality in \eqref{eqht}.
\end{proof}

\begin{lemma}\label{W}
    For any prime divisor $\mathcal{W}=\mathcal{W}_{\beta,n,k}$ in  \eqref{div(s)}, we have
    \begin{equation*}
        (\overline{\mathcal{M}}_n \cdot \pi_n^*\overline{\mathcal{H}}_1\cdots \pi_n^*\overline{\mathcal{H}}_r
 |\mathcal{W})=0.
    \end{equation*}
    Moreover, we have
    \begin{equation*}
        (\overline{\mathcal{E}}_n \cdot \pi_n^*\overline{\mathcal{H}}_1\cdots \pi_n^*\overline{\mathcal{H}}_r
 |\mathcal{W})\ge 0.
    \end{equation*}
\end{lemma}
\begin{proof}
    Let $\mathcal{D}:= \pi_n(\mathcal{W})$ be the scheme-theoretic image of $\mathcal{W}$ in $\mathcal{B}$. We can take a projective flat morphism of relative dimension $1$ between arithmetic varieties:
    $$
    q:\mathcal{W}^\prime \longrightarrow \mathcal{D}^\prime,
    $$
 where $\mathcal{W}^\prime$ and $\mathcal{D}^\prime$ are projective models of $\mathcal{W}$ and $\mathcal{B}$.
    The intersection product is invariant when pulling back to $\mathcal{W}^\prime$, and then we get
\begin{equation*}
    \begin{aligned}
        (\overline{\mathcal{M}}_n \cdot \pi_n^*\overline{\mathcal{H}}_1\cdots \pi_n^*\overline{\mathcal{H}}_r
 |\mathcal{W})=\deg(\mathcal{M}_{n,F}) \cdot \overline{\mathcal{H}}_1|_{\mathcal{D}^\prime}\cdots \overline{\mathcal{H}}_r|_{\mathcal{D}^\prime},
    \end{aligned}
\end{equation*}
by applying the projection formula in \cite[Lemma 4.6.1(3)]{YZ26}. Since 
$\deg(\mathcal{M}_{n,F})=0$ by \eqref{degree-zero}, this proves the equality in this lemma.

As for the boundary divisor $\overline{\mathcal{E}}_n$, we have the same 
equality with $\deg(\mathcal{E}_{n,F})>0$. Since 
$\overline{\mathcal{H}}_i$ is nef and so is the restriction 
$\overline{\mathcal{H}}_i|_{\mathcal{D}^\prime}$. Thus, the 
intersection of $\overline{\mathcal{E}}_n$ and 
$\pi_n^*\overline{\mathcal{H}}_1,\dots, \pi_n^*\overline{\mathcal{H}}_r$ 
over $\mathcal{W}$ is nonnegative.
\end{proof}

It remains to bound the height term and integral term in $I_{n,\delta}$. By Lemma \ref{Z} and \ref{W}, we can rewrite the inductive formula in \eqref{In} as
\begin{equation}\label{In-final}
I_{n,\delta}=\frac{1}{N_{n,k}}\sum_{\alpha}c_{\alpha,n,k} [E_{\alpha,n,k}:F]h_{\overline{\mathcal{M}}_n}^{\overline{\mathscr{H}}}(z_{\alpha,n,k})
   -\frac{1}{N_{n,k}}\int_{\mathcal{X}_n(\mathbb{C})}{\log\|s_{n,k}\|} \,\mu_n.
\end{equation}
For fixed $n$ and $\delta >0$, the signed measure $\mu_n$ has smooth and bounded density with respect to a smooth volume form 
on $\mathcal{X}_n(\mathbb{C})$. Therefore, we have
$$
\lim_{k\rightarrow \infty}{\frac{1}{N_{n,k}}\int_{\mathcal{X}_n(\mathbb{C})}{\log\|s_{n,k}\|} \,\mu_n}=0,
$$
since $\frac{1}{N_{n,k}}{\log\|s_{n,k}\|}\rightarrow 0$ in $L^1$-sense as $k\to \infty$. We define 
$$
J_{n,\delta}:= \overline{\mathcal{E}}_n \cdot \overline{\mathcal{A}}_{n,\delta}\cdot \pi_n^*\overline{\mathcal{H}}_1\cdots \pi_n^*\overline{\mathcal{H}}_r,
$$
which is a nonnegative number by ampleness of $\overline{\mathcal{A}}_{n,\delta}$ and nefness of $\pi_n^*\overline{\mathcal{H}}_i$.
We apply the same induction formula to $\overline{\mathcal{E}}_n$ as in \eqref{In}, as the archimedean integration term tends to zero and the intersection over $\mathcal{W}$ is nonnegative for $\overline{\mathcal{E}}_n$, we get
\begin{equation}\label{Jn-bound}
\limsup_{k\rightarrow \infty} \left|{\frac{1}{N_{n,k}}\sum_{\alpha}c_{\alpha,n,k} [E_{\alpha,n,k}:F]h_{\overline{\mathcal{E}}_n}^{\overline{\mathscr{H}}}(z_{\alpha,n,k})}\right| \leq J_{n,\delta}. 
\end{equation}
Note that $I_{n,\delta}$ and $J_{n,\delta}$ are both independent of $k$. Letting $k\rightarrow \infty$ in \eqref{In-final}, we have
\begin{equation*}
    |I_{n,\delta}|\le \varepsilon_nJ_{n,\delta}.
\end{equation*}
  As $\overline{\mathcal{A}}_{n,\delta}=\overline{\mathcal{A}}_{n}+\delta\cdot \overline{\mathcal{B}}_{n} $, letting $\delta  \rightarrow 0$, we have
  \begin{equation}\label{ineq-final}
 \left|\overline{\mathcal{M}}_n \cdot \overline{\mathcal{A}}_{n}\cdot \pi_n^*\overline{\mathcal{H}}_1\cdots \pi_n^*\overline{\mathcal{H}}_r\right| 
 \leq 
 \varepsilon_n\,\overline{\mathcal{E}}_n \cdot \overline{\mathcal{A}}_{n}\cdot \pi_n^*\overline{\mathcal{H}}_1\cdots \pi_n^*\overline{\mathcal{H}}_r,
  \end{equation}
by multilinearity and continuity of intersection products. Since 
$$
\overline{\mathcal{E}}_{n}=p_n^* \overline{\mathcal{E}}_0, \quad \overline{\mathcal{A}}_{n}=p_n^* \overline{\mathcal{A}}, 
$$
every line bundle in the right-hand side of \eqref{ineq-final} is a 
pullback from $\mathcal{X}_0$. Applying the projection formula of 
birational morphism, we have
$$
J_n:=J_{n,0}= \overline{\mathcal{E}}_n \cdot \overline{\mathcal{A}}_{n}\cdot \pi_n^*\overline{\mathcal{H}}_1\cdots \pi_n^*\overline{\mathcal{H}}_r = \overline{\mathcal{E}}_0 \cdot \overline{\mathcal{A}}\cdot \pi_0^*\overline{\mathcal{H}}_1\cdots \pi_0^*\overline{\mathcal{H}}_r = J_0.
$$
We write $I_n:=I_{n,0}$. Since $\varepsilon_n \rightarrow 0$ as $n\rightarrow \infty$, we obtain 
$$
I_n \longrightarrow 0.
$$
By the definition of arithmetic intersection of adelic line bundles, this exactly says
$$
\overline{M}\cdot\overline{A} \cdot \pi_C^*\overline{H}_1 \cdots \pi_C^*\overline{H}_r =0.
$$
Hence, we finish the proof of Proposition \ref{Q-implies-H}.
\let\qed\relax
\end{proof}

\subsection{Inequality case}
Now we can prove the inequality case of the arithmetic Hodge index over 
quasi-projective curves.
\begin{proposition} \label{ineqq}
Let $F$ be a finitely generated 
field of transcendence degree $r\geq 1$ over $\mathbb{Q}$, and $C$ be a normal, geometrically integral quasi-projective curve over 
$F$. Suppose $ \overline{M} \in \widehat{\mathrm{Pic}}(C/\mathbb{Z})_{\mathrm{int},\mathbb{Q}} $ is an integrable adelic line bundle on $C$ satisfying $\deg(\widetilde{M})=0$. Then
$$
\overline{M}^2\cdot \pi_C^*\overline{H}_1\cdots \pi_C^*\overline{H}_r \leq 0,
$$
for any $r$ nef adelic line bundles $\overline{H}_1,\dots, \overline{H}_r\in \widehat{\mathrm{Pic}}(F/\mathbb{Z})_{\mathbb{Q},\mathrm{nef}}$.
\end{proposition}
\begin{proof}
We can choose a normal projective arithmetic model $\mathcal{B}$ of $F$ over $\mathbb{Z}$. For each $n$, let 
$$
\pi_n:\mathcal{X}_n \longrightarrow \mathcal{B}_n 
$$
be a normal arithmetic projective model of $\pi_C: C\rightarrow \Spec F$ with birational morphism $\mathcal{B}_n\rightarrow \mathcal{B}$, with
$\mathcal{B}_0=\mathcal{B}$, and 
$\overline{\mathcal{M}}_n$ is defined over $\mathcal{X}_n $ by 
passing to a common dominating model. Since $C/F$ is a curve, we 
can identify every $\mathcal{X}_{n,F}$ with a unique projective 
compactification $X$ of $C$. Such $X$ is unique up to isomorphism.

 We let $E:=X\setminus C$ be 
the subscheme of $X$ with reduced structure, regarded as an effective divisor on $X$. Since the projective case has already been 
solved by Yuan--Zhang in \cite{YZ13}, we always assume that $E$ is non-empty. 
   
Suppose that $\overline{M}=(M,(\overline{\mathcal{M}}_n,\mathcal{X}_n,\ell_n)_n)$ is an integrable adelic line bundle with respect to a boundary divisor $\overline{\mathcal{E}}_0$ on $\mathcal{X}_0$ with generic fiber ${\mathcal{E}}_{0,F}=E$. Replacing $\overline{\mathcal{M}}_n$ by
   $$
   \overline{\mathcal{M}}_n-\frac{\deg_{X_n}(\mathcal{M}_{n,F})}{\deg_X(E)}\cdot \overline{\mathcal{E}}_n,
   $$
where $\overline{\mathcal{E}}_n$ is the pullback of $\overline{\mathcal{E}}_0$ to $\mathcal{X}_n$, we can assume that $\deg(\mathcal{M}_{n,F})=0$ for all $n$ similarly as in \eqref{modify}.

Since we can choose $\{\overline{\mathcal{M}}_{n}\}$ representing $\overline{M}$ with $\deg(\mathcal{M}_{n,F})=0$, we can apply the inequality case of arithmetic Hodge index theorem of Yuan–Zhang to each $\overline{\mathcal{M}}_n$. Hence, 
    $$
     \overline{\mathcal{M}}_n^2\cdot \pi_n^*\overline{H}_1\cdots \pi_n^*\overline{H}_r
     \le 0.
   $$
    Letting $n\rightarrow \infty$, we get $\overline{M}^2\cdot \pi_C^*\overline{H}_1\cdots \pi_C^*\overline{H}_r \leq 0$.

\end{proof}

\subsection{Equality case}\label{Sec-Modify}
In this subsection, we will show that for any adelic line bundle $\overline{M}$ that satisfies the 
hypothesis of Theorem \ref{thm:intro-Q1r}, we can find a Cauchy 
sequence of model adelic $\mathbb{Q}$-line bundles $\{\overline{\mathcal{M}}_n\}$ 
representing $\overline{M}$, such that the geometric part 
$$ \mathcal{M}_{n,F}= 0$$
in analogy with the process in $\S$\ref{Section-vertcurve}.

Let $F$ be a finitely generated 
field of transcendence degree $r\geq 1$ over $\mathbb{Q}$, and let $C$ 
be a normal geometrically integral quasi-projective curve over 
$F$. Suppose $ \overline{M} \in \widehat{\mathrm{Pic}}(C/\mathbb{Z})_{\mathrm{int},\mathbb{Q}} $ is an integrable adelic line bundle on $C$ satisfying $\deg(\widetilde{M})=0$ and
$$
\overline{M}^2\cdot \pi_C^*\overline{H}_1\cdots \pi_C^*\overline{H}_r =0,
$$
for any $r$ integrable adelic line bundles $\overline{H}_1,\dots, \overline{H}_r\in \widehat{\mathrm{Pic}}(F/\mathbb{Z})_{\mathrm{int},\mathbb{Q}}$. 





We write $M_n:= \mathcal{M}_{n,F}$.
Identifying every $\mathcal{X}_{n,F}$ with $X$, we can regard
 $M_n$ as a line bundle in $\mathrm{Pic}^0(X)$, and interpret $\widetilde{M}$ as the limit of $M_n$ in $\mathrm{Pic}^0(X)_{\mathbb{R}}$, which we still denote by $\widetilde{M}$ by abuse of notation. This is possible because 
 $$
 -\varepsilon_nE \leq M_n-M_m \leq \varepsilon_n E, ~~~\text{as}~n\rightarrow \infty
 $$
for $m>n\geq 0$ and $\varepsilon_n \rightarrow 0$, from the definition of the adelic line bundle $\overline{M}$.

 By \cite[Theorem 6.5.1]{YZ26}, for each $n$,
there is an adelic $\mathbb{Q}$-line bundle $\overline{M}_{n,0} \in \widehat{\mathrm{Pic}}(X/\mathbb{Z})_{\mathrm{int},\mathbb{Q}}$ with underlying line bundle $M_n$ such that
$$
\pi_{X,*}\langle\overline{M}_{n,0},\overline{V}\rangle = 0,~~~  \forall ~\overline{V}\in  \widehat{\mathrm{Pic}}(X/\mathbb{Z})_{\mathrm{vert},\mathbb{Q}},
$$
and,
$$
\pi_{X,*}\langle\overline{M}_{n,0},\overline{M}_{n,0}\rangle = -2\,\hat{\mathfrak{h}}(M_n) \in \widehat{\mathrm{Pic}}(F/\mathbb{Z})_{\mathrm{int},\mathbb{Q}}.
$$
Now we fix a polarization $\overline{\mathscr{H}}:=(\overline{\mathcal{H}},\dots,\overline{\mathcal{H}})$ where 
 $\overline{\mathcal{H}}$ is an arithmetically ample hermitian line 
 bundle on a projective arithmetic model of $\Spec F$. Let $J:=\mathrm{Pic}^0_{X/F}$ be the Jacobian variety 
of $X$.
We define the Moriwaki height by
\begin{equation*}
    \begin{aligned}
\hat{h}^{\overline{\mathscr{H}}}: J(F)_{\mathbb{Q}}~  &\longrightarrow~ \mathbb{R}_{\geq 0}, \\
M_n ~~~ &\longmapsto~ \hat{\mathfrak{h}}(M_n)
\cdot \overline{\mathcal{H}}^r.
    \end{aligned}
\end{equation*}

 Then $\hat{h}^{\overline{\mathscr{H}}}$ extends to a semipositive 
quadratic form on $J(F)_{\mathbb{R}}=J(F)\otimes_{\mathbb{Z}}{\mathbb{R}}$. 
By the theorem of Lang--N\'eron, $J(F)_{\mathbb{R}}$ is a 
finite-dimensional $\mathbb{R}$-linear space
for abelian varieties over a finitely generated field $F$.  
By the Northcott property of the Moriwaki height introduced in \cite[Theorem 5.3.1]{YZ26}, we see that $\hat{h}^{\overline{\mathscr{H}}}$ 
is positive definite on $J(F)_{\mathbb{R}}$, since tensoring with $\mathbb{R}$ kills the torsion part of $J(F)$. 
    
     We are going to show that $\hat{h}^{\overline{\mathscr{H}}}(M_n)$ converges to $0$. Since
    $
    \overline{V}_n := \overline{M}_n - \overline{M}_{n,0}
    $
is vertical, we have
\begin{equation*}
    \begin{aligned}
 \pi_{X,*}\langle\overline{M}_{n},\overline{M}_{n}\rangle \cdot \overline{\mathcal{H}}^r
 &=
 \pi_{X,*}\langle\overline{M}_{n,0},\overline{M}_{n,0} \rangle\cdot \overline{\mathcal{H}}^r + \pi_{X,*}\langle\overline{V}_n,\overline{V}_n \rangle\cdot \overline{\mathcal{H}}^r\\
 &=-2\hat{h}^{\overline{\mathscr{H}}}(M_n) + \pi_{X,*}\langle\overline{V}_n,\overline{V}_n \rangle\cdot \overline{\mathcal{H}}^r \longrightarrow 0.
    \end{aligned}
\end{equation*}
By Proposition \ref{ineqq}, the second term on the right-hand side is non-positive. Since $-2\hat{h}^{\overline{\mathscr{H}}}(M_n)\leq 0$ is also non-positive, we have
$$ \hat{h}^{\overline{\mathscr{H}}}(M_n)\rightarrow 0, $$
as $n\to\infty$. By the continuity of the quadratic form $\hat{h}^{\overline{\mathscr{H}}}$ on the finite-dimensional $\mathbb{R}$-vector space, we obtain
$$
\hat{h}^{\overline{\mathscr{H}}}(\widetilde{M})=0.
$$
Since the restriction of $\widetilde{M}$ to $C$ is $M$, and the natural map $J(F)_\mathbb{Q} \rightarrow J(F)_\mathbb{R}$ is injective, we have
$$
\widetilde{M}=0\in\mathrm{Pic}^0(X)_{\mathbb{R}},~\text{and}~ M=0\in\mathrm{Pic}^0(C)_{\mathbb{Q}}.
$$
Therefore, for any Cauchy sequence of model adelic $\mathbb{Q}$-divisors $\{\overline{\mathcal{D}}_n\}$ representing $\overline{M}$, we can modify it so that $\mathcal{D}_{n,F}=0$ by repeating the process in $\S$\ref{curve-K}.

\subsection{An explicit calculation}\label{explicit}
Let $\overline{M}$, $\overline{\mathcal{M}}_n$, $X$ and $\pi_n: \mathcal{X}_n \longrightarrow \mathcal{B}_n$ be as in Proposition \ref{ineqq}.
We can assume that $\mathcal{M}_{n,F}=0$ as discussed above. We 
denote $\mathcal{B}:=\mathcal{B}_0$. Let $\overline{\mathcal{H}}$ be an arithmetically ample hermitian 
line bundle on $\mathcal{B}$, and use the same symbol 
$\overline{\mathcal{H}}$ simply for its pullback to 
$\mathcal{B}_n$. Note that the measure 
$c_1(\overline{\mathcal{H}})^r$ has zero mass over proper Zariski closed
subsets outside the open subset where $\mathcal{B}_n\rightarrow \mathcal{B}$ is an isomorphism. We set
$$
\overline{\mathcal{Q}}_n= \langle \overline{\mathcal{M}}_n,\overline{\mathcal{M}}_n \rangle_{\mathcal{X}_n/\mathcal{B}_n}.
$$
Suppose $\{\overline{\mathcal{D}}_n\}$ is a Cauchy sequence of model $\mathbb{Q}$-divisors converging to $\overline{M}$ with 
$\overline{\mathcal{M}}_n=\mathcal{O}(\overline{\mathcal{D}}_n)$. The 
canonical section $s_n$ of $\mathcal{D}_n$ induces a section $\langle s_n,s_n \rangle$ of 
$\mathcal{Q}_n$ as follows: since $\mathcal{M}_{n,F}=0$, we can fix an 
isomorphism $\alpha_n:\mathcal{O}_X \rightarrow \mathcal{M}_{n,F}$ which maps $1$ to $s_{n,F}$. Then we get an isomorphism of one-dimensional $F$-linear space:
$$
F \xrightarrow{\sim} \langle\mathcal{O}_X,\mathcal{O}_X \rangle_{X/F} \xrightarrow{\langle\alpha_n,\alpha_n \rangle}  \langle \mathcal{M}_{n,F},\mathcal{M}_{n,F} \rangle_{X/F} \xrightarrow{\sim} \mathcal{Q}_{n,F} 
$$
Therefore, the trivializing section of $\mathcal{Q}_{n,F}$ extends to a unique rational section of $\mathcal{Q}_n$ on $\mathcal{B}_{n}$, which we 
denote by $\langle s_n,s_n \rangle$. Note that we actually take a rational 
section $s_n^{\otimes a}$ for some integer $a$ such that $\mathcal{M}_n^{\otimes a}$ is integral and write formally $\langle s_n,s_n \rangle=\langle s_n^{\otimes a},s_n^{\otimes a} \rangle^{\otimes 1/a^2}$. 
We write its corresponding arithmetic $\mathbb{Q}$-divisor on $\mathcal{B}_n$ as
$$
\widehat{\div}(\langle s_n,s_n \rangle)= \left(\sum_{Z\text{ prime}}{q_{n,Z}Z},~ \varphi_n  \right).
$$
Here the sum is taken over prime divisors $Z$ of $\mathcal{B}_n$ 
and $\varphi_n:=-\log\|\langle s_n,s_n \rangle \|_{\overline{\mathcal{Q}}_n}$ is the model 
Green function of $\overline{\mathcal{Q}}_n$ on 
$\mathcal{B}_n(\mathbb{C})$. Hence the following arithmetic intersection product can be written as
\begin{equation} \label{QH^r}
I_n :=\overline{\mathcal{Q}}_n\cdot \overline{\mathcal{H}}^r=  \sum_{Z\text{ prime}}{q_{n,Z}(\overline{\mathcal{H}}^r|Z)}+ \int_{\mathcal{B}_n(\mathbb{C})}\varphi_n c_1(\overline{\mathcal{H}})^r.
\end{equation}
And $I_n \rightarrow 0$ as $n\to \infty$ by the assumption on $\overline{M}$. Since $\overline{\mathcal{H}}$ is arithmetically ample, we see that 
$(\overline{\mathcal{H}}^r|Z)\geq 0$ and 
$c_1(\overline{\mathcal{H}})^r$ is semipositive. We will first prove that for all $n$, we have
\begin{equation}\label{nonpositive}
q_{n,Z} \leq 0, \quad \mathrm{and}\quad \varphi_n \leq 0 \text{ almost everywhere}.
\end{equation}

Now fix one $\pi_n: \mathcal{X}_n \longrightarrow \mathcal{B}_n$, we will omit some $n$ in discussion of the following local Deligne pairing for simplicity. 
For a prime divisor $Z$ of $\mathcal{B}$, let $\eta_Z$ be its generic point. Since $\mathcal{B}$ is normal, the 
    localization $\mathcal{O}_{\mathcal{B},\eta_Z}$ is a 
    discrete valuation ring. The base change of $\pi_n$ under $\Spec \mathcal{O}_{\mathcal{B},\eta_Z}\rightarrow \mathcal{B} $  can be written as
    $$
    \pi_Z: \mathcal{X}_{Z}=\mathcal{X}\times_{\mathcal{B}} \Spec \mathcal{O}_{\mathcal{B},\eta_Z} \longrightarrow \Spec \mathcal{O}_{\mathcal{B},\eta_Z}.
    $$
    Then $\mathcal{X}_{Z}$ is a relative arithmetic surface over 
$\mathcal{O}_{\mathcal{B},\eta_Z}$, whose generic fiber is $X \rightarrow \Spec F$. Since Deligne pairing is compatible with base change, we have
    $$
    \overline{\mathcal{Q}}_n|_{\Spec \mathcal{O}_{\mathcal{B}_n,\eta_Z}}= \langle \overline{\mathcal{M}}_n|_{\mathcal{X}_{Z}},\overline{\mathcal{M}}_n|_{\mathcal{X}_{Z}} \rangle.
    $$
We define $$\mathcal{D}_Z:=\mathcal{D}_n|_{\mathcal{X}_{Z}}.$$
Since $ \mathcal{D}_{n,F}=0$, the divisor $\mathcal{D}_Z$ is vertical.

We claim that $q_{Z} = \mathcal{D}_{Z}^2$. In fact, we 
let $s_Z:=s_n|_{\mathcal{X}_Z}$ be the canonical section of $\mathcal{D}_Z$. Then the restriction of  $\langle s_n,s_n \rangle$ to $\mathcal{Q}_n|_{\Spec \mathcal{O}_{\mathcal{B},\eta_Z}}$ is 
$$
\langle s_Z,s_Z \rangle=\langle s_{n}|_{\mathcal{X}_Z},s_n|_{\mathcal{X}_Z} \rangle.
$$
We denote the closed point of $\Spec \mathcal{O}_{\mathcal{B},\eta_Z}$ by $v=v_Z$, then by 
definition,
$$
q_Z=\ord_{v}(\langle s_Z,s_Z \rangle),
$$
which is the vanishing order of section $\langle s_Z,s_Z \rangle$
at $v$. Thus, we have $\div(\langle s_Z,s_Z \rangle)=q_Z\cdot [v]. $

    Moreover, the geometric Deligne pairing for $\pi_Z$ gives
    $
    \div(\langle s_Z,s_Z \rangle)=\pi_{Z,*}\langle\mathcal{D}_Z,\mathcal{D}_Z\rangle.
    $
Since $\mathcal{D}_Z$ is vertical, and $\pi_{Z,*}\langle\mathcal{D}_Z,\mathcal{D}_Z\rangle$ is supported on the closed point $v$, we have
$$
\pi_{Z,*}\langle\mathcal{D}_Z,\mathcal{D}_Z\rangle=(\mathcal{D}_Z^2)\cdot [v].
$$
As a consequence, we have 
\begin{equation}\label{q=D^2}
q_Z=\mathcal{D}_Z^2.
\end{equation}

Now by abuse of notation, we let $v=v_Z$ be the divisorial valuation on $F$ corresponding to $Z$ and $F_v$ be the completion of $F$. Then the 
completion $\widehat{\mathcal{O}}_{\mathcal{B},\eta_Z}$ of 
$\mathcal{O}_{\mathcal{B},\eta_Z}$ has fraction field $F_v$. We get a model
$$
\widehat{\mathcal{X}}_Z:=\mathcal{X}_Z \times_{\mathcal{O}_{\mathcal{B},\eta_Z}}\widehat{\mathcal{O}}_{{\mathcal{B},\eta_Z}},
$$
and the analytification of its generic fiber is a Berkovich curve
$$ X_v^{\mathrm{an}}:=(X \times_F F_v)^{\mathrm{an}}. $$
Therefore the vertical divisor $\mathcal{D}_{Z}$ defines a model function on $ X_v^{\mathrm{an}}$, which we denote by
$$
g_{n,v}(x)=-\log\|1 \|_{\mathcal{D}_{Z}} (x).
$$
Denote by $\langle1,1\rangle$ the canonical section of $\langle \overline{\mathcal{O}(g_{n,v})},\overline{\mathcal{O}(g_{n,v})} \rangle$. 
Therefore, the metric of the local pairing $\langle \overline{\mathcal{O}(g_{n,v})},\overline{\mathcal{O}(g_{n,v})} \rangle$ at $v$ is
\begin{equation}\label{11}
-\log\|\langle 1,1\rangle\|_{n,v}=\int_{X_v^{\mathrm{an}}}g_{n,v}c_1(\overline{\mathcal{O}(g_{n,v})})=\int_{X_v^{\mathrm{an}}}g_{n,v}dd^cg_{n,v},
\end{equation}
by the fibral metric formula in \cite[$\S$4.6.2]{YZ26}.
Also we have
\begin{equation} \label{11'}
-\log\|\langle 1,1\rangle\|_{n,v}=\lambda_vq_{n,Z}, 
\end{equation}
since $q_{n,Z}$ is the vanishing order of $\langle 1,1\rangle$; here 
$\lambda_v=-\log|\varpi_v|_v$ and $\varpi_v$ is a uniformizer of 
$\mathcal{O}_{\mathcal{B},\eta_Z}$; see \cite[$\S$2.5]{CL06}. Therefore, we are about to show the first part of \eqref{nonpositive} in the following lemma:

\begin{lemma} \label{qnZ}
    We have $q_{n,Z} = \mathcal{D}_{n,Z}^2 = \frac{1}{\lambda_{v}} \int_{X_{v}^\mathrm{an}} {g_{n,v} dd^c g_{n,v}}$, and $q_{n,Z} \leq 0$ for 
    every prime divisor $Z$ of $\mathcal{B}_n$.
\end{lemma}
\begin{proof}
    The equality is obtained from \eqref{q=D^2}, \eqref{11} and \eqref{11'}. Write the special fiber of $\mathcal{X}_Z$ as $\mathcal{X}_{Z,s}=\sum_{1\leq i \leq r}{\nu_iV_i}$. Let 
    $\mathcal{D}_Z=\sum_{1\leq i \leq r}{d_iV_i}$.  By defining 
    $E_i:= \nu_iV_i$, $b_i:=d_i/\nu_i$ and passing to a regular semistable model, we can get the same equality as in \eqref{self-intersection}:
    $$
\mathcal{D}_{n,Z}^2 =-\sum_{i<j}m_{ij}(b_i-b_j)^2,~\text{where}~ m_{ij}:=\nu_i\nu_j V_i\cdot V_j\geq 0.
    $$
Hence, we have $q_{n,Z} = \mathcal{D}_{n,Z}^2 \leq 0$.
    
\end{proof}

The second part of the claim \eqref{nonpositive} is included in the following lemma:
\begin{lemma}
    For every $w\in \mathcal{B}_n(\mathbb{C})$, let 
    $X_w(\mathbb{C}):=\pi_n^{-1}(w)$ be the one-dimensional complex fiber. 
    Let $U_n$ be the open subset of $\mathcal{B}_n(\mathbb{C})$ on which 
    $\pi_n$ is smooth. Shrink $U_n$ so that the canonical section of $\mathcal{D}_n$ has no zero or pole on the fibers over $w\in U_n$. Then we have 
    \begin{equation}\label{phi}
    \varphi_n(w)= \int_{X_w(\mathbb{C})}g_{n,w}(z)dd^cg_{n,w}(z),~\text{for}~ w \in U_n.
    \end{equation}
 Here $g_{n,w}(z)=-\log\| s_n(z)\|_{\overline{\mathcal{M}}_n|_{\mathcal{X}_w}}$ is the 
 model Green function of $\overline{\mathcal{D}}_n|_{\mathcal{X}_w}$. Moreover, we have $\varphi_n(w) \leq 0$ on $U_n$.
\end{lemma}
\begin{proof}
    Recall that $\langle s_n,s_n\rangle$ is a section of 
    $\overline{\mathcal{Q}}_n=\langle \overline{\mathcal{M}}_n,\overline{\mathcal{M}}_n \rangle$ and its metric at $w$ is given by
    $$
    -\log\|\langle s_n,s_n\rangle \|(w)= \int_{X_w(\mathbb{C})}{g_{n,w}dd^c g_{n,w}},
    $$
by the archimedean case of fibral metric in \cite[$\S$4.2.2]{YZ26}. Since $\varphi_n=-\log\|\langle s_n,s_n\rangle \|$ by 
definition, this proves \eqref{phi}. Since $g_{n,w}$ is a smooth real-valued function, we have
$$
\varphi_n(w)= \int_{X_w(\mathbb{C})}{g_{n,w}\,{dd^c} g_{n,w}}=\frac{1}{\pi i}\int_{X_w(\mathbb{C})} 
\partial g_{n,w}\wedge \overline{\partial} g_{n,w} \le 0,
$$
as in \eqref{arch}.
\end{proof}

Now we return to the expression of $I_n$:
$$
I_n =  \sum_{Z\text{ prime}}{q_{n,Z}(\overline{\mathcal{H}}^r|Z)}+ \int_{\mathcal{B}_n(\mathbb{C})}\varphi_n(w)c_1(\overline{\mathcal{H}})^r.
$$
Since $q_{n,Z}\leq  0$ and $(\overline{\mathcal{H}}^r|Z)\geq  0$, 
each summand is non-positive. As $\varphi_n\leq 0$ on $U_n$ and $c_1(\overline{\mathcal{H}})^r \geq 0$ and has 
mass $0$ on any proper Zariski closed subset, the integral term is also non-positive.  As $n\rightarrow\infty$, since $I_n$ tends to zero, we have
$$
q_{n,Z} (\overline{\mathcal{H}}^r|Z)\longrightarrow 0, ~\text{and}~\int_{\mathcal{B}_n(\mathbb{C})}\varphi_n(w)c_1(\overline{\mathcal{H}})^r\longrightarrow 0.  
$$
Note that $(p_n^*\overline{\mathcal{H}}^r|Z)=0$ when $Z$ is exceptional over $\mathcal{B}_0$ on which 
$\overline{\mathcal{H}}$ is defined.
However, since $Z$ represents a divisorial valuation $v_Z$ and $q_{n,Z}$ is 
independent of the choice of $\overline{\mathcal{H}}$, we can always 
pass to a model of $F$ on which $v_Z$ is represented by a prime 
divisor and choose another polarization $\overline{\mathcal{H}}^\prime$ with
$$ {(\overline{\mathcal{H}}^\prime}^r|Z^\prime)>0, $$
where $Z^\prime$ is the strict transform of $Z$. From Lemma \ref{qnZ}, the number $q_{n,Z}$ is independent of the model realizing the divisorial valuation 
$v=v_Z$ of $Z$ over $F$. Therefore, we conclude that as $n\rightarrow \infty$, we have 
$$q_{n,Z} \longrightarrow 0,  $$
for a fixed divisorial valuation $v_Z$.
Since $c_1(\overline{\mathcal{H}})^r$ is semipositive and $\varphi_n\leq 0$ on a full measure subset of $\mathcal{B}(\mathbb{C})$, we see that 
$$ \varphi_n \longrightarrow 0,$$
almost everywhere after passing to a subsequence $\{\varphi_{n_j}\}$ of $\{\varphi_n\}$. For simplicity, we replace $\{\varphi_{n}\}$ by this subsequence.

Let $g_v$ be the Green function of $\overline{M}$ on $C_v^{\mathrm{an}}$.
The sequence $\{g_{n,v}|_{C_v^{\mathrm{an}}}\}$ locally uniformly converges to 
$g_v$ on $C_v^{\mathrm{an}}$, which is explained as follows: since $g_{n,v}$ 
is the Green function induced by a common section $1$ on the model $(\mathcal{X}_Z,\mathcal{D}_{n}|_{\mathcal{X}_Z})$ for all $n$, we have
$$
|g_{n,v}-g_v|\leq\varepsilon_n\tilde{g}_0,
$$
where $\tilde{g}_0$ is the Green function induced by the boundary divisor $\mathcal{E}_0$ on $\mathcal{X}_Z$. As $\tilde{g}_0$ is bounded on every compact subset of
$C_v^{\mathrm{an}}$, we are done. We refer to \cite[Proposition 3.4.1]{YZ26} for details in analytifications of Green functions. Now we can prove Theorem \ref{thm:intro-R1r} which states that $R(1,r)$ holds.

\begin{proof}[Proof of Theorem \ref{thm:intro-R1r}]
Now we can apply Lemma \ref{nonarch} to see that $g_v$ is a constant 
on $C_v^{\an}$, for every divisorial valuation $v=v_Z$ of $F$. As for 
the archimedean case when $w\in \mathcal{B}(\mathbb{C})$, similarly we have
$$
g_{n,w} \longrightarrow g_{w}, 
$$
locally uniformly on $C_w(\mathbb{C})$ as $n\to \infty$ by the Cauchy condition. Applying Lemma \ref{lem:zero-energy} to $\{g_{n,w}\}$ we see that $g_w$ is 
constant for almost every $w\in \mathcal{B}(\mathbb{C})$. It remains to show that this suffices to imply that $\overline{M}$ is a pullback from $\Spec F$.

First, we choose one $c\in C(\overline{F})$. Let $c'$ be its corresponding closed point, and write $L:=\kappa(c^\prime)$ for its residue field.
Define
\begin{equation}\label{eq:curve-descent-N}
    \overline{N}
    :=\mathfrak{h}_{\overline{M}}(c)=\frac{1}{[L:F]}N_{L/F}\bigl(\overline{M}|_{c'}\bigr)
    \in\widehat{\mathrm{Pic}}(F/\mathbb{Z})_{\mathrm{int},\mathbb{Q}}.
\end{equation}
We aim to show that $\overline{M}-\pi_C^*\overline{N}=0\in \widehat{\mathrm{Pic}}(C/\mathbb{Z})_{\mathrm{int},\mathbb{Q}}.$ Suppose that $\{\overline{\mathcal{P}}_n\}$ is a sequence of model adelic divisors  
converging to an adelic divisor $\overline{\mathcal{P}}$, which represents $\overline{M}-\pi_C^*\overline{N}$. Then we can assume that $\mathcal{P}_{n,F}=0$ as in $\S$\ref{Sec-Modify}. Hence it suffices to show that $\overline{\mathcal{P}}=0$ is a zero adelic divisor, i.e.
\begin{equation}\label{tends-zero}
-\varepsilon_n\overline{\mathcal{E}}_n \leq \overline{\mathcal{P}}_n \leq \varepsilon_n \overline{\mathcal{E}}_n, ~~~\text{as}~n\rightarrow \infty,
\end{equation}
for a sequence $\{\varepsilon_n\}$ of positive rational numbers decreasing to $0$.

\noindent\emph{Geometric part:}
We claim that for any divisorial valuation $u$ of $F(C)$, 
$$ \lim_{n\rightarrow \infty} u(\mathcal{P}_n)=0. $$
We normalize $u$ with value group $\mathbb{Z}$.
Then for any prime divisor $Y$ of $\mathcal{X}_n$, we have $|\ord_Y(\mathcal{P}_n)|\leq \varepsilon_n \ord_{Y}(\mathcal{E}_n)$. Since $\mathcal{X}_n$ is normal, we have
$$
-\varepsilon_n\mathcal{E}_n
\leq \mathcal{P}_n\leq\varepsilon_n\mathcal{E}_n,
$$
which is the geometric part of \eqref{tends-zero}.

Now we proceed by proving the above claim.
If $u|_F$ is trivial, then $u$ corresponds to a closed point of the smooth projective curve $X/F$. Since $\mathcal{P}_{n,F}=0$, we have
$$
u(\mathcal{P}_n)=0.
$$

If $u|_F$ is nontrivial, then we see that $u|_F$ is still a divisorial valuation on $F$ by \cite[Lemma 4.5]{HKP21}. Write $u|_F=e \cdot v$, where 
$e\ge 1$ is an integer and $v$ is a divisorial valuation on $F$ of value group $\mathbb{Z}$. 

For any $a\in F$, we have $|a|_v=\exp(-\lambda_v v(a))$. Then the multiplicative seminorm
$$
|f|_{\xi_u}:= \exp(-\lambda_v u(f)/e), \mathrm{~for~} f\in F(C),
$$
restricts to $|\cdot|_v$ on $F$, hence defining a point $\xi_u$ in $C_v^{\mathrm{an}}$. Let 
 $g_{\overline{\mathcal{P}}_{n},v} $ be the local Green function of $\overline{\mathcal{P}}_n$ at $v$, then we get
$$
g_{\overline{\mathcal{P}}_{n},v}(\xi_u )=\lambda_v u(\mathcal{P}_n)/e;
$$
We refer to the analytification of adelic divisors in \cite[$\S$3.3.1]{YZ26}. Since $g_{\overline{\mathcal{P}}_{n},v} \rightarrow 0$ locally uniformly, we see that 
$$
u(\mathcal{P}_n) \longrightarrow 0,
$$
which proves the claim.

\noindent\emph{Arithmetic part:} By definition, we have 
$$
-\varepsilon_n\overline{\mathcal{E}}_n \leq \overline{\mathcal{P}}_n-\overline{\mathcal{P}}_m \leq \varepsilon_n \overline{\mathcal{E}}_n, ~\mathrm{for}~ m>n,
$$
Let $\mathcal C$ be a quasi-projective arithmetic model of $C$, identified with an open set of $\mathcal X_n$ for every $n$. Then we have an inequality of Green functions:
$$
-\varepsilon_n g_{\overline{\mathcal{E}}_n} \leq g_{\overline{\mathcal{P}}_n}-g_{\overline{\mathcal{P}}_m} \leq \varepsilon_n g_{\overline{\mathcal{E}}_n}~\text{on}~ \mathcal C(\mathbb{C}).
$$
We will show that $g_{\overline{\mathcal{P}}}$ vanishes everywhere in Lemma \ref{zero-g} 
below. Once we know that $g_{\overline{\mathcal{P}}}\equiv 0$ on $C(\mathbb{C})$, then letting $m\to \infty$ gives 
$$
-\varepsilon_n g_{\overline{\mathcal{E}}_n} \leq g_{\overline{\mathcal{P}}_n} \leq \varepsilon_n g_{\overline{\mathcal{E}}_n}~\text{on}~ \mathcal C(\mathbb{C}).
$$
This proves the arithmetic part of \eqref{tends-zero}.
Thus, we conclude that $\overline{\mathcal{P}}=0$ in the group of adelic divisors on $C$, which implies that $\overline{M}-\pi_C^*\overline{N}=0$. Hence we are done.
\end{proof}
\begin{lemma} \label{zero-g}
    The Green function $g_{\overline{\mathcal{P}}}$ of $\overline{\mathcal{P}}$ on $\mathcal C(\mathbb{C})$ satisfies 
    $$
    g_{\overline{\mathcal{P}}}\equiv 0.
    $$
    And we have $g_{\overline{\mathcal{P}}_n} \rightarrow 0$ on $\mathcal C(\mathbb{C})$ as $n\rightarrow \infty$.
\end{lemma}
\begin{proof}
    
Choose a common quasi-projective arithmetic model of $C$ on which $M$ and $\pi_C^*\overline{N}$ are represented, and shrink the base
to a non-empty open subscheme $\mathcal{B}^{\circ}$ of $\mathcal{B}$ so that
\[
    \pi:\mathcal{C}\longrightarrow\mathcal{B}^{\circ}
\]
is a smooth model of $C\to\Spec F$, the base $\mathcal{B}^\circ(\mathbb{C})$ is smooth, and
the closure $\mathcal{T}$ of $c^\prime$ in $\mathcal{C}$ is a multisection which is finite \'etale of
degree $\ell$ over $\mathcal{B}^{\circ}$. Let $p:\mathcal{T} \rightarrow \mathcal{B}^\circ$ be the surjective morphism. On
$\mathcal{B}^{\circ}$, we can represent $\overline{N}\in \widehat{\mathrm{Pic}}(F/\mathbb{Z})_{\mathrm{int},\mathbb{Q}}$ by
\[
    \frac{1}{\ell}N_{\mathcal{T}/\mathcal{B}^{\circ}}
        \bigl(\overline{M}|_{\mathcal{T}}\bigr) \in \widehat{\mathrm{Pic}}(\mathcal{B}^{\circ}/\mathbb{Z})_{\mathrm{int},\mathbb{Q}}.
\]
Identifying $\mathcal{B}(\mathbb{C})$ with an open subset of 
$\mathcal{B}_n(\mathbb{C})$ for every $n$, we can set $\Sigma$ as the 
intersection of  $(\cap_{n}{U_n}) \cap\mathcal{B}^{\circ}(\mathbb{C})$ and the full-measure subset where $\varphi_n \to 0$, 
which is still a full-measure
subset so that $g_w$ is a constant on $\Sigma$. Suppose that
$$g_w(z)\equiv c_w, \text{ for all }w\in\Sigma.$$  
Then we have
\[
    g_{\overline{N}}(w)
    =\frac{1}{\ell}\sum_{t\in p^{-1}(w)}g_w(t)
    =c_w.
\]
Since $g_{\overline{M}}$ and $g_{\pi^*\overline{N}}$ are both equal to 
constant $c_w$ on each analytification of fiber $\pi^{-1}(w)$ for every 
$w\in \Sigma$,
we see that the Green function $g_{\overline{\mathcal{P}}}$ vanishes on
$\pi^{-1}(\Sigma)$.  Since
$\Sigma$ is dense with respect to the strictly positive measure $c_1(\overline{\mathcal{H}})^r$ and $\pi$ is an open mapping, we see that $\pi^{-1}(\Sigma)$ is dense in $\mathcal{C}(\mathbb{C})$.
Since $g_{\overline{\mathcal{P}}}$ is continuous, we obtain that
\begin{equation*}
    g_{\overline{\mathcal{P}}}=0
    \quad\text{on }\mathcal{C}(\mathbb{C}).
\end{equation*}
Letting $m\to \infty$, we have $|g_{\overline{\mathcal{P}}_n,w}| \leq \varepsilon_n g_{\overline{\mathcal{E}}_n,w}$.
\end{proof}



\section{Proof of assertion $H(d,r)$} \label{Hdr}
 In this section, we will prove the following theorem, which asserts that an integrable adelic 
 line bundle $\overline{M}$ on quasi-projective $V/F$ is a pullback from $\Spec F$ if its 
 associated vector-valued height at every algebraic point of $V(\overline{F})$ is numerically trivial.
 \begin{theorem}[$H(d,r)$ holds] \label{H(d,r)}
  Let $F$ be a finitely generated 
field of transcendence degree $r\ge 0$ over $\mathbb{Q}$, and $V$ be a normal, geometrically integral quasi-projective variety over 
$F$ of dimension $d$. Suppose that $ \overline{M} \in \widehat{\mathrm{Pic}}(V/\mathbb{Z})_{\mathrm{int},\mathbb{Q}} $ is an 
integrable adelic $\mathbb Q$-line bundle on $V$ such that 
$$
\mathfrak{h}_{\overline{M}}(x)\equiv_{\mathrm{num}} 0 \in \widehat{\mathrm{Pic}}(F/\mathbb{Z})_{\mathbb{Q}},$$
for every $x\in V(\overline{F})$.
Then $$ \overline{M}\in\pi_{V}^{*}\widehat{\mathrm{Pic}}(F/\mathbb{Z})_{\mathbb{Q}}. $$
\end{theorem}

\subsection{Construction of generic pencils}\label{pencil}
We recall the following construction of Xie--Yuan in \cite{XY22}, which is essential in our future discussions:
\begin{lemma}[{\cite[Lemma 3.7]{XY22}}]\label{lem3.7}
    Let $k$ be an algebraically closed field. Let $K/k$ be a
finitely generated extension of transcendence degree $\mathrm{tr.deg}(K/k) > 1$.
Let $(S,\mathcal{M})$ be a polarization of $K/k$, consisting of a projective
model $S$ of $K/k$ and an ample line bundle $\mathcal{M}$ over $S$. Assume that $\mathcal{M}$ is very ample
over $S$.
Then there are infinitely many intermediate fields $k^\prime$
of $K/k$, algebraically closed in $K$ and of transcendence degree $1$ over $k$, together with
a geometrically integral subvariety $H$ in $S_{k^\prime}$ of codimension $1$ satisfying the following properties:
\begin{enumerate}
    \item [(1)] the divisor $H$ of $S_{k^\prime}$ is linearly equivalent to $\mathcal{M}_{k^\prime}$;
    \item [(2)] the composition $H$ → $S_{k
^\prime}\rightarrow S$ induces an isomorphism between
the generic points of $H$ and $S$.
\end{enumerate}
\end{lemma}
The proof of the above lemma works also for the finitely generated field 
$F$. Briefly, we start with a very ample line bundle on a normal 
projective model $X$ of $V$, and two independent sections $s_0, s_1$. Then we can get a flat morphism from 
$\widetilde{X}$ to $\mathbb{P}^1$ by blowing up the base locus. We denote the generic fiber of $\widetilde X$ by $\widetilde{X}_\eta$ under this map. Then 
there is a geometrically integral subvariety $H$, whose generic point is isomorphic to that of $X$.
Passing to the open immersions $V\subseteq X$ and $V_\eta \subseteq \widetilde{X}_\eta$ and 
choosing two such algebraically independent pencils, we obtain two morphisms of schemes:
$$\rho_i:V_i \longrightarrow V,~~i=1,\,2$$
where $V_i$ is a quasi-projective subvariety over $F_i=F(t_i)$ with $\dim_{F_i}{V_i}=d-1$ and function field $F_i(V_i)=F(V)$.

\subsection{Key lemmas on the generic fibers}
Let $T\subseteq \mathbb{P}^1_F$ be a non-empty open subset of the projective line over finitely generated field $F$. Suppose that there is a normal quasi-projective curve $Z$, endowed with a finite flat morphism $p \colon Z \rightarrow T$ of degree $m$ and a morphism $s : Z \rightarrow V$.

Let $F^\prime=F(T)$ be the function field of $\mathbb{P}^1_F$, which is identified with a subfield of $F(V)$ via $\widetilde{X}\rightarrow \mathbb{P}^1_F$.
Note that a point of $V_{i}(\overline{F^\prime})$ corresponds to an $\overline{F^\prime}$-point of $V_{\overline{F^\prime}}$,
not to a point of $V(\overline{F})$ covered directly by the original height hypothesis. However, we have obtained the following result:

\begin{lemma}[Numerical triviality of generic-fiber point heights]\label{height pullback}
    Let $f:V\to T$ be one of
the above pencil maps, where $T$ is a non-empty open subset of $\mathbb{P}^1_F$. Let $\rho\colon W \rightarrow V$ be one of the  
morphisms $\rho_i$. Let $\overline{M}$ be an integrable adelic 
$\mathbb{Q}$-line bundle on $V$ such that for every $v\in V(\overline{F})$,
$$\mathfrak{h}_{\overline{M}}(v) \equiv_{\mathrm{num}} 0 ~~\text{in}~ \widehat{\mathrm{Pic}}(F/\mathbb{Z})_{\mathrm{int},\mathbb{Q}} .$$
Then we have 
$$ \mathfrak{h}_{\rho^*\overline{M}}(x)\equiv_{\mathrm{num}} 0~~\text{in}~ \widehat{\mathrm{Pic}}(F^{\prime}/\mathbb{Z})_{\mathrm{int},\mathbb{Q}} ,$$
for every $x\in W(\overline{F^\prime})$.
\end{lemma}
\begin{proof}
    We fix one $x\in W(\overline{F^\prime})$, let $x^\prime$ be the closed point of $W$ containing $x$, and put
$E_x = \kappa(x^\prime)$, and $m = \left[E_x : F^\prime\right]$.
Let $Z$ be the normalization of $T$ in $E_x$, then $Z$ is a quasi-projective model of $\Spec E_x$. After shrinking $T$, the morphism $p:Z\rightarrow T$
is finite and flat of degree $m$, and the point $x^\prime$
extends to a morphism
$s \colon Z \rightarrow V$.
We define 
$$\overline{H}_x := \frac{1}{m}p_*\langle s^*\overline{M} \rangle= \frac{1}{m}N_{Z/T}(s^*\overline{M}) \in  \widehat{\mathrm{Pic}}(T/\mathbb{Z})_{\mathrm{int},\mathbb{Q}}.$$
Since the Deligne pairing of adelic line bundles is compatible with the base change $\alpha:\Spec F^\prime \rightarrow T$ by \cite[Theorem 4.1.3]{YZ26}, the restriction of $\overline{H}_x$ to the generic point of $T$ is exactly
\begin{equation*}
    \begin{aligned}
        \left.\overline{H}_x\right|_{\Spec F^\prime}
        =&\frac{1}{m}\alpha^*\left(N_{Z/T}(s^*\overline{M})\right)\\
        =&\frac{1}{m}N_{E_x/F^\prime}(x^*\rho^*\overline{M}) \\
        =&\frac{1}{m}\left\langle\left.\rho^*\overline{M}\right|_{x^\prime}\right\rangle=\mathfrak{h}_{\rho^*\overline{M}}(x),
    \end{aligned}
\end{equation*}
by the definition of vector-valued height. For a point $t\in T(\overline{F})$ , let $t^\prime$ be its corresponding closed 
point and $L:=\kappa(t^\prime)$ be its residue field. Let $z_1,\dots,z_\ell$ be the preimage of $t^\prime$ under 
$p:Z\rightarrow T$. Let $e_i$ be the local multiplicity of $z_i$. Denote the $0$-dimensional subscheme $p^{-1}(t^\prime)$ by $Z_{t^\prime}$, then there is a commutative diagram:
\[ 
\begin{tikzcd}
Z_{t^\prime} \arrow{r}{\iota} \arrow[swap]{d}{q} & Z \arrow{d}{p} \\
\Spec L \arrow{r}{t}& T
\end{tikzcd}
\]
Therefore we observe that
\begin{equation*}
\begin{aligned}
    \overline{H}_x|_{t^\prime}=&{\frac{t^*(p_*\langle s^*\overline{M} \rangle )}{m}}\\
    =& \frac{1}{m}q_*\left\langle s^*\overline{M}|_{Z_{t^\prime}} \right\rangle\\
    =&
    \sum_{i=1}^{\ell}{\frac{e_i}{m}{N_{\kappa(z_i)/L}\left(s^* \overline{M} |_{z_i}\right)}}.
\end{aligned}
\end{equation*}
This expresses  
\begin{equation*}
\begin{aligned}
    \mathfrak{h}_{\overline{H}_x}(t)=&\frac{1}{\left[ L:F\right]}N_{L/F}(\overline{H}_x|_{t^\prime})\\
    =&\frac{1}{\left[ L:F\right]}{\sum_{i=1}^{\ell}{\frac{e_i}{m}{N_{\kappa(z_i)/F}\left(s^* \overline{M} |_{z_i}\right)}}} \\
    =& {\sum_{i=1}^{\ell}{\frac{e_i\left[\kappa(z_i):L\right]}{m}{\mathfrak{h}_{s^* \overline{M}}}(z_i)}},
\end{aligned}
\end{equation*}
as a linear combination of $\mathfrak{h}_{\overline{M}}(s(z_i))=\mathfrak{h}_{s^* \overline{M}}(z_i)$. Consequently, 

$$ \mathfrak{h}_{\overline{H}_x}(t) \equiv_{\mathrm{num}} 0,\quad\text{for all }t\in T(\overline{F}).$$
Applying $H(1,r)$ to $\overline{H}_x$, we know that $\overline{H}_x$ is a numerically trivial adelic $\mathbb{Q}$-line bundle on $T$.
It remains to prove that $\mathfrak{h}_{\rho^*\overline{M}}(x)=\left.\overline{H}_x\right|_{\Spec F^\prime}$ is also 
numerically trivial.

Let $\overline{A}_0,\dots,\overline{A}_r$ be $r+1$ arbitrarily chosen adelic line bundles in $\widehat{\mathrm{Pic}}(F^\prime/\mathbb{Z})_{\mathrm{int},\mathbb{Q}}.$ We 
can take a non-empty open subvariety $T^\circ \subseteq T$ over $F$ 
which is a quasi-projective model for all of the 
$\overline{A}_j$. Since $\mathfrak{h}_{\overline{H}_x}(t) \equiv_{\mathrm{num}} 0,$ for all $t\in T^{\circ}(\overline{F})$, 
we see that $\left.\overline{H}_x\right|_{T^{\circ}}$ is numerically trivial by $H(1,r)$. Therefore
$$ \left.\overline{H}_x\right|_{T^{\circ}} \cdot \overline{A}_0 \cdots \overline{A}_r =0. $$
By the definition of arithmetic intersection product on $\Spec F^{\prime}$, this equality is exactly
$$ \mathfrak{h}_{\rho^*\overline{M}}(x) \cdot \overline{A}_0 \cdots \overline{A}_r =0. $$
We conclude that $\mathfrak{h}_{\rho^*\overline{M}}(x)$ is numerically trivial for every $x\in W(\overline{F^\prime})$.
\end{proof}

We will need this lemma on injectivity of the pullback maps of adelic line bundles over finitely 
generated fields in the following sections.
\begin{lemma}\label{extension injectivity}
    Let $F \subseteq E$ be an extension of finitely generated fields over $\mathbb{Q}$. Then the pullback map 
    $$ \alpha^*:~ \widehat{\mathrm{Pic}}(F/\mathbb{Z})_{\mathbb{Q}} \longrightarrow \widehat{\mathrm{Pic}}(E/\mathbb{Z})_{\mathbb{Q}}$$ 
is injective via $\alpha:\Spec E \rightarrow \Spec F$.
\end{lemma}
\begin{proof}
    Suppose $\overline{A} \in \widehat{\mathrm{Pic}}(F/\mathbb{Z})_{\mathbb{Q}}$ is an adelic line bundle with 
    $\alpha^*\overline{A}=0 \in \widehat{\mathrm{Pic}}(E/\mathbb{Z})_{\mathbb{Q}} $. Choose a smooth integral 
quasi-projective variety $Y/F$ with function field $E$ and structure morphism $\pi_Y: Y\rightarrow \Spec F$. Then $\left.\pi_Y^*\overline{A}\right|_{\Spec E}=0 $. Applying the injectivity of map
$$ \widehat{\mathrm{Pic}}(Y/\mathbb{Z})_{\mathbb{Q}} \longrightarrow \widehat{\mathrm{Pic}}(E/\mathbb{Z})_{\mathbb{Q}}. $$
we see that $\pi_Y^*\overline{A}=0$ by \cite[Corollary 3.4.2]{YZ26}.
Now we choose a closed point $y \in Y$ which corresponds to $\beta: \Spec L\rightarrow Y$. Hence $\beta^*\pi_Y^*\overline{A}=0$. Denote by $\pi_y:\Spec L\rightarrow \Spec F$ the composition $\pi_Y\circ\beta$, then we have
\begin{equation*}
    \begin{aligned}
        0=&N_{L/F}(\beta^*\pi_Y^*\overline{A})\\
        =&N_{L/F}(\pi_y^*\overline{A})=\left[L:F\right]\overline{A}.
    \end{aligned}
\end{equation*}
Hence $\overline{A}=0 \in \widehat{\mathrm{Pic}}(F/\mathbb{Z})_{\mathbb{Q}}$.
\end{proof}

\subsection{Proof of Theorem \ref{thm:intro-Q-induction}}
In this subsection, we will prove the induction argument that claims the assertion $H(d-1,r+1)$ implies $H(d,r)$. Therefore, we can deduce that $H(d,r)$ holds from the curve case by induction.

   Let $V/F$ be a normal, geometrically integral quasi-projective variety and $E:= F(V)$ be the 
   function field of $V$. We can assume that $F$ is algebraically closed in $E$. Suppose 
   $\overline{M}$ is an integrable adelic 
   line bundle on $V$.  Following Lemma \ref{lem3.7}, we can choose two varieties $V_1$ and 
   $V_2$, which are defined over fields $F_1,~F_2$, identified with the intermediate fields $F\subseteq F_i \subseteq E$, such that
   $$ F_i=F(t_i),\quad\dim_{F_i}(V_i)=d-1, \text{ and }  F_i(V_i)=E,  $$
   where $t_i\in F_i$. And the two pencil maps are chosen so that $F_1\neq F_2$ and $F_1$,$F_2$ 
   are algebraically independent. Then we have the following descent lemma, which asserts that if a given integrable adelic line bundle $\overline{P}$ on $\Spec E$ is a pullback from both $\Spec F_1$ and $\Spec F_2$, then $\overline{P}$ indeed comes from $\Spec F$.

\begin{lemma}\label{Descent} 
    Let $F_1,\,F_2\subseteq E$ be two finitely generated fields over $F$ such 
    that $F_i/F$ has transcendence degree $1$, $F_1F_2/F$ has transcendence 
    degree $2$ and $F_1 \cap F_2 = F$. Let $\alpha_i:\Spec E \rightarrow \Spec F_i$ denote the corresponding morphism. Suppose that there is an adelic 
    $\mathbb{Q}$-line bundle 
    $\overline{P}\in \widehat{\mathrm{Pic}}(E/\mathbb{Z})_{\mathrm{int},\mathbb{Q}}$ contained in 
    
    $$ \alpha_1^*\,\widehat{\mathrm{Pic}}(F_1/\mathbb{Z})_{\mathbb{Q},\mathrm{int}}\, \cap\alpha_2^*\,\widehat{\mathrm{Pic}}(F_2/\mathbb{Z})_{\mathbb{Q},\mathrm{int}}. $$
Then $\overline{P}$ is contained in the image 
$$ \mathrm{im}(\widehat{\mathrm{Pic}}(F/\mathbb{Z})_{\mathrm{int},\mathbb{Q}}\longrightarrow \widehat{\mathrm{Pic}}(E/\mathbb{Z})_{\mathrm{int},\mathbb{Q}}).$$
\end{lemma}
\begin{proof}
    By assumption, there are two integrable adelic $\mathbb{Q}$-line bundles 
    $\overline{N}_1 \in \widehat{\mathrm{Pic}}(F_1/\mathbb{Z})_{\mathrm{int},\mathbb{Q}}$ 
    and $\overline{N}_2\in \widehat{\mathrm{Pic}}(F_2/\mathbb{Z})_{\mathrm{int},\mathbb{Q}}$
    such that
    $$
    \overline{P}= \alpha_1^*\overline{N}_1=\alpha_2^*\overline{N}_2. 
    $$
    
    Let $F_{12}:=F_1F_2$ and denote by $\beta_i:\Spec F_{12}\rightarrow \Spec F_i$ the natural morphism. Lemma \ref{extension injectivity} indicates the injectivity of
    $\widehat{\mathrm{Pic}}(F_{12}/\mathbb{Z})_{\mathrm{int},\mathbb{Q}}\longrightarrow \widehat{\mathrm{Pic}}(E/\mathbb{Z})_{\mathbb{Q}}$. Hence, we have 
    $$
    \beta_1^*\overline{N}_1=\beta_2^*\overline{N}_2.
    $$

    This reduces to the case where $E=F_{12}$.
Now we choose smooth projective geometrically integral curves $X_i/F$ whose
function field is $F(X_i) = F_i$ for $i=1,2$. There are non-empty open subschemes $C_i \subseteq X_i$ such that $\overline{N}_i$ can be defined over $C_i$. 

 Since both $C_1$ and $C_2$ are geometrically integral, the function field of $C_1\times_F C_2$ is $F_{12}$. 
This fiber product is also normal and geometrically integral.
Let $p_i:C_1\times_F C_2 \rightarrow C_i$ be the projection onto the $i$-th factor for $1\leq i\leq 2$. Applying the injectivity argument in \cite[Corollary 3.4.2]{YZ26}, we see that

$$  p_1^*\overline{N}_1=p_2^*\overline{N}_2\quad \text{in }  \widehat{\mathrm{Pic}}(C_1 \times_F C_2/\mathbb{Z})_{\mathrm{int},\mathbb{Q}}.  $$

Next we choose a closed point $c_2 \in C_2$ with residue field $\kappa(c_2)=L$. We consider the restriction of both sides in $ p_1^*\overline{N}_1=p_2^*\overline{N}_2$ to $C_1\times\left\{ c_2\right\}\simeq C_{1,L}$, then we obtain 
\begin{equation}\label{C1}
\overline{N}_{1,L}=\pi_{1,L}^*(c_2^*\overline{N}_2) \quad\text{in }  \widehat{\mathrm{Pic}}(C_{1,L}/\mathbb{Z})_{\mathrm{int},\mathbb{Q}},
\end{equation}
through the following commutative diagram:
\[ 
\begin{tikzcd}
 C_1\times \{c_2\} \arrow{r}{\pi_{1,L}} \arrow[swap]{d}{} & \Spec L \arrow{d}{c_2} \\
C_1\times C_2 \arrow{r}{p_2} \arrow[swap]{d}{p_1} & C_2 \arrow{d}{} \\
C_1 \arrow{r}{\pi_1}& \Spec F
\end{tikzcd}
\]
Applying the Deligne pairing of relative dimension $0$ or equivalently the norm map to both sides of \eqref{C1}, we get 
$$
\left[L:F\right]\overline{N}_1=N_{C_{1,L}/C_1}(\pi_{1,L}^*(c_2^*\overline{N}_2))=\pi_1^*(N_{L/F}(c_2^*\overline{N}_2)),
$$
where the second equality is due to the compatibility of the Deligne pairing and base changes. Therefore, we have 
$$
\overline{N}_1= \pi_1^*\overline{N}, \quad\text{with } \overline{N}:=  \frac{1}{\left[L:F\right]}(N_{L/F}(c_2^*\overline{N}_2))\in  \widehat{\mathrm{Pic}}(F/\mathbb{Z})_{\mathrm{int},\mathbb{Q}}.
$$
This proves the lemma since $\overline{P}=\alpha_1^*\overline{N}_1=\alpha_1^*\pi_1^*\overline{N}$.
    
 \end{proof}

  Now we assume that $H(d-1,r+1)$ holds and indicate that it implies $H(d,r)$.
Let $V$ be a normal, geometrically integral and geometrically connected variety over a finitely generated field $F/\mathbb{Q}$ of 
$\mathrm{tr.deg}_{\mathbb{Q}}(F)=r \geq 0$. Suppose $\overline{M} \in \widehat{\mathrm{Pic}}(V/\mathbb{Z})_{\mathbb{Q}}$ is an integrable adelic line bundle whose vector-valued height satisfies that for every $x\in V(\overline{F})$, we have
$$ \mathfrak{h}_{\overline{M}}(x)\equiv_{\mathrm{num}} 0 \in \widehat{\mathrm{Pic}}(F/\mathbb{Z})_{\mathrm{int},\mathbb{Q}}, 
$$
if $r>0$, or
$ h_{\overline{M}}(x)=0$ for all $x\in V(\overline{F})$ if $r=0$. 

Let $\rho_i:V_i\rightarrow V,~(i=1,2)$ 
be the two morphisms which are obtained by the generic pencil construction in $\S$\ref{pencil}. By Lemma \ref{height pullback}, we have 
$$ \mathfrak{h}_{\rho_i^*\overline{M}}(x)\equiv_{\mathrm{num}} 0~~\text{in}~ \widehat{\mathrm{Pic}}(F_i/\mathbb{Z})_{\mathrm{int},\mathbb{Q}}, 
$$
for every $x\in V_i(\overline{F}_i)$.
Since $\rho_1^*\overline{M}$ and $\rho_2^*\overline{M}$ satisfy the hypothesis of $H(d-1,r+1)$, we see that $\rho_i^*\overline{M}$ is a pullback of an adelic $\mathbb Q$-line bundle in 
$\widehat{\mathrm{Pic}}(F_i/\mathbb{Z})_{\mathbb{Q}}$. Hence $\left.\overline{M}\right|_{\Spec E} $ is a pullback from both 
$\widehat{\mathrm{Pic}}(F_1/\mathbb{Z})_{\mathbb{Q}}$ and $\widehat{\mathrm{Pic}}(F_2/\mathbb{Z})_{\mathbb{Q}}$. Applying Lemma \ref{Descent} to $\left.\overline{M}\right|_{\Spec E} $, we see that 
$$
\left.\overline{M}\right|_{\Spec E}=\left.\pi_V^*\overline{N}\right|_{\Spec E}
$$
for some $\overline{N}\in \widehat{\mathrm{Pic}}(F/\mathbb{Z})_{\mathbb{Q}}$, 
where $\pi_V:V\rightarrow \Spec F$ is the structure morphism. Since the inclusion 
map $j:\Spec E\rightarrow V$ of generic point of $V$ induces an injection of 
arithmetic Picard 
groups of essentially quasi-projective varieties by \cite[Corollary 3.4.2]{YZ26}, we have 
$$
\overline{M}=\pi_V^*\overline{N}.
$$
Since we prove the curve case $H(1,d+r-1)$ in $\S$\ref{curve over F}, we finish our 
proof of $H(d,r)$ for all $d \ge 1, r\geq 0$ by starting at $H(1,d+r-1)$ and applying the induction argument in Theorem \ref{thm:intro-Q-induction}.

\section{Hodge index type theorem and numerical triviality}\label{chap5-Rdr}
In this section, we will prove Theorem \ref{chap1-Rdr} by reducing to the case of complete intersection curves and applying the assertions $R(1,r)$ and $H(d,r)$ which are proved in $\S$\ref{curve over F} and $\S$\ref{Hdr} respectively.

\begin{theorem}[Theorem \ref{chap1-Rdr}] \label{Rdr}
 Let $V$ be a normal geometrically integral quasi-projective variety of dimension $d\ge 1$ over a finitely generated field $F$ of 
     $\mathrm{tr.deg}_{\mathbb{Q}}(F)=r\ge 0$. Suppose
$\overline M$ is an integrable adelic line bundle on $V$.
Assume that for every collection
 $\overline A_1,\ldots,\overline A_{d-1} \in \widehat{\mathrm{Pic}}(V/\mathbb{Z}) _{\mathrm{int},\mathbb Q}$, one has
\[
 \widetilde M\cdot
 \widetilde A_1\cdots\widetilde A_{d-1}=0
\]
and
\[
 \pi_{V,*}\bigl\langle
 \overline M,\overline M,
 \overline A_1,\ldots,\overline A_{d-1}
 \bigr\rangle
 \equiv_{\mathrm{num}}0.
\]
When $d=1$, the first equality is interpreted as
$\deg(\widetilde M)=0$. Then
\[
 \overline M\in
 \pi_V^*\widehat{\mathrm{Pic}}(F/\mathbb{Z})_{\mathbb Q}.
\]
\end{theorem}

 The hypothesis of Theorem \ref{Rdr} is equivalent to that $ \overline{M} \in \widehat{\mathrm{Pic}}(V/\mathbb{Z})_{\mathrm{int},\mathbb{Q}} $
 is represented by a Cauchy sequence $\{\overline{\mathcal{M}}_n\}$ of hermitian 
 $\mathbb{Q}$-line bundle on a cofinal 
family of projective arithmetic models $\{\mathcal{X}_n\}$ of $V$ over $\mathbb{Z}$ such that 
for arbitrarily chosen model adelic line bundles 
$\overline{A}_{i} \in  \widehat{\mathrm{Pic}}(V/\mathbb{Z})_{\mathrm{mod},\mathbb{Q}}$ on $V$, we have 
$$
\mathcal{M}_{n,F}\cdot \mathcal{A}_{1,F}\cdots\mathcal{A}_{d-1,F}\rightarrow 0,~\text{as }n\rightarrow\infty,
$$ 
and 
\begin{equation}\label{numtrivial}
\pi_{V,*}\left\langle\overline{M},\overline{M},\overline{A}_{1},\dots,\overline{A}_{d-1}\right\rangle \equiv_{\mathrm{num}}0 \in \widehat{\mathrm{Pic}}(F/\mathbb{Z})_{\mathrm{int},\mathbb{Q}},
\end{equation} 
for hermitian $\mathbb{Q}$-line bundles $\overline{\mathcal{A}}_{i}$ representing 
$\overline{A}_i$ on a common arithmetic model of $V$ and with $\mathcal{A}_{i,F}$ 
being its generic fiber of the underlying line bundle $\mathcal{A}_{i}$. 

Throughout this section, products indexed by an empty set are interpreted as $1$,
and empty sums as $0$.
In particular, all statements involving $r$ base factors include the
case $r=0$. We use additive notation for tensor products. 

\subsection{Inequality and orthogonality}
\label{subsec-Rdr-ineq}

Let $\overline M\in
    \widehat{\mathrm{Pic}}(V/\mathbb{Z})_{\mathrm{int},\mathbb{Q}}$ be an integrable adelic line 
    bundle satisfying the hypotheses of Theorem \ref{Rdr}. Thus, for all integrable adelic
$\mathbb{Q}$-line bundles $\overline A_1,\ldots,\overline A_{d-1}$
on $V/\mathbb{Z}$ and
$\overline H_1,\ldots,\overline H_r$ on $\Spec F$, one has
\begin{equation}\label{hypotheses}
\begin{aligned}
    \widetilde M \cdot\widetilde A_1 &\cdots
        \widetilde A_{d-1}=0,\\
    \overline M^2\cdot\overline A_1\cdots\overline A_{d-1}
        &\cdot\pi_V^*\overline H_1\cdots\pi_V^*\overline H_r=0.
\end{aligned}
\end{equation}

As in the curve case in Proposition \ref{ineqq}, we first prove the 
inequality case of the arithmetic Hodge index theorem for higher dimensional varieties as follows. 
\begin{proposition}[Inequality case of arithmetic Hodge index over quasi-projective varieties]
\label{lem:smR-inequality}
Let $\overline A$ be a model adelic line bundle on $V$ induced by an arithmetically ample hermitian
$\mathbb{Q}$-line bundle $\overline{\mathcal A}$ on a normal projective
arithmetic model $\pi_{\mathcal X}:\mathcal X\to\mathcal B$,
where $\mathcal B$ is a normal projective arithmetic model of $F$
and the generic fiber $X$ is a compactification of $V$. Let
$\overline L$ be an integrable adelic line bundle on $V$. If
\[
    \widetilde L \cdot \widetilde{A}^{d-1}=0,
\]
then, for all nef adelic $\mathbb{Q}$-line bundles
$\overline H_1,\ldots,\overline H_r$ on $F/\mathbb{Z}$,
\[
    \overline L^2\cdot\overline A^{d-1}
    \cdot\pi_V^*\overline H_1\cdots\pi_V^*\overline H_r\leq 0.
\]
\end{proposition}

\begin{proof} 
Choose a sequence of model adelic line bundles $\{\overline{\mathcal{L}}_n\}_{n\ge 0}$ which approximates $\overline{L}$
on a sequence of dominating normal models
$p_n:\mathcal X_n\to\mathcal X_0$, where we identify $\mathcal{X}$ with $\mathcal{X}_0$. Write $X_n:=\mathcal X_{n,F}$, and $\overline{\mathcal A}_n:=p_n^*\overline{\mathcal A}$. Then we have
$$\qquad a:=\mathcal{A}_F^d= \mathcal{A}_{n,F}^d>0,$$
by the projection formula.
 Let
\[
    \delta_n:=
    \frac{\mathcal L_{n,F}\cdot\mathcal A_{n,F}^{d-1}}{a}.
\]
The hypothesis on $\widetilde{L}$ implies that $\delta_n\to 0$ as $n\to \infty$. And the
$\mathbb{Q}$-line bundle
$\mathcal L_{n,F}-\delta_n\mathcal A_{n,F}$ satisfies 
$$
(\mathcal L_{n,F}-\delta_n\mathcal A_{n,F})\cdot \mathcal{A}_{n,F}^{d-1}=0.
$$
We write $\pi_n:=\pi_{\mathcal X}\circ p_n$ for the morphism
$\mathcal X_n\to\mathcal B$. Since $\mathcal A_{n,F}$ is nef and big,
the inequality part of the projective arithmetic Hodge index theorem
\cite[Theorem~1.2]{YZ13} gives
\[
    (\overline{\mathcal L}_n-\delta_n\overline{\mathcal A}_n)^2
    \cdot\overline{\mathcal A}_n^{d-1}
    \cdot\pi_n^*\overline H_1\cdots\pi_n^*\overline H_r\leq 0.
\]
Consequently, we obtain the wanted inequality by taking $n\to \infty$.
\end{proof}

We have the following useful corollary that plays an essential role in this section. 
\begin{corollary}[Orthogonality]
\label{Cor-Orthogonality}
Let $\overline A$ be a model adelic line bundle on $V$ induced by an arithmetically ample hermitian line bundle, and let
$\overline L$ be an integrable adelic line bundle satisfying
\[
    \widetilde L\cdot\widetilde A^{d-1}=0.
\]
Then for any $\overline{M}$ satisfying \eqref{hypotheses} and for every nef $r$-tuple $\overline H_1,\ldots,\overline H_r$ of adelic line bundles on $\Spec F$, one has
\begin{equation}\label{orthogonality}
    \overline M\cdot\overline L\cdot\overline A^{d-1}
    \cdot\pi_V^*\overline H_1\cdots\pi_V^*\overline H_r=0.
\end{equation}
In particular, this holds for every vertical model line bundle
$\overline T$ on an arithmetic model of $V$ since $\widetilde{T}=0$.
\end{corollary}

\begin{proof}
Applying Proposition~\ref{lem:smR-inequality} to
$\overline M+t\overline L$ for every $t\in\mathbb{Q}$, we have
$$
2t\cdot\overline M\cdot\overline L\cdot\overline A^{d-1}
    \cdot\pi_V^*\overline H_1\cdots\pi_V^*\overline H_r+t^2\cdot \overline{L}^2\cdot\overline A^{d-1}
    \cdot\pi_V^*\overline H_1\cdots\pi_V^*\overline H_r \leq 0
$$
Therefore, we see that the coefficient of $t$ must be zero, which yields
$$
\overline M\cdot\overline L\cdot\overline A^{d-1}
    \cdot\pi_V^*\overline H_1\cdots\pi_V^*\overline H_r=0.
$$
Hence we are done. 
\end{proof}

We say an integrable adelic line bundle $\overline{M}$ on $V$ 
satisfies \emph{orthogonality} if \eqref{orthogonality} holds for $\overline{M}$, 
every arithmetically ample $\overline{A}$, every integrable adelic 
line bundle $\overline{L}$ satisfying $ \widetilde L\cdot\widetilde A^{d-1}=0 $, and every nef $r$-tuple 
$\overline H_1,\ldots,\overline H_r$ of adelic line bundles on $\Spec F$.

\subsection{Push-forward construction} \label{subsec-pushforward}
Let $\overline M$ and $V$ be as in $\S$\ref{subsec-Rdr-ineq}.
Fix an $x_0\in V(\overline F)$. Let $x_0'$ be the corresponding closed
point, let $L_0:=\kappa(x_0')$, and put $\ell_0:=[L_0:F]$. Define
\begin{equation}\label{eq:smR-N0}
    \overline N_0:=\mathfrak h_{\overline M}(x_0)
    =\frac{1}{[L_0:F]}N_{L_0/F}(x_0'^*\overline M),
    \qquad
    \overline M':=\overline M-\pi_V^*\overline N_0.
\end{equation}
Then we have $\mathfrak h_{\overline{M}^\prime}(x_0)=0$. Choose a rational section $s$ of $M$ which is regular and nonvanishing at $x_0'$. Let
\[
    \overline D:=\widehat{\operatorname{div}}(s)
\]
be the adelic $\mathbb{Q}$-divisor associated with $\overline{M}$.
Then the restriction of $\overline{D}$
to $x_0'$ is defined. We set $\overline\Delta_0$ to be the adelic divisor
on $\Spec F$ representing $\overline N_0$, as
\begin{equation}\label{eq:smR-P}
    \overline\Delta_0:=\frac{1}{[L_0:F]}
        N_{L_0/F}(x_0'^*\overline D).
    \end{equation}
The norm in this formula is taken with the rational section induced
by $s$. The adelic $\mathbb{Q}$-divisor
$$
\overline P:=\overline D-\pi_V^*\overline\Delta_0
$$
represents the adelic $\mathbb{Q}$-line bundle $\overline M-\pi_V^*\overline N_0=\overline{M}^\prime$.

The geometric part of $\pi_V^* \overline{N}_0$ is zero. Applying the
projection formula therefore shows that $\overline M'$ also satisfies
\eqref{hypotheses}. 

We work on a resolution of $V$ so that we can choose a projective arithmetic model $\mathcal{X}_0$ of $V$ with smooth generic fiber 
$X_0$ and $X_0\setminus V$ being a simple normal crossing divisor. Its generic fiber 
$X_0:= \mathcal{X}_{0,F}$ is a projective compactification of
$V$. Let $\overline{M}$ be an integrable adelic line bundle on $V$ satisfying \eqref{hypotheses} and 
$\{\overline{\mathcal{M}}_n\}_{n\ge 0}$ be a Cauchy sequence of model adelic line bundles approximating $\overline{M}$ with respect to a strictly effective
boundary divisor $\overline{\mathcal E}_0$ on $\mathcal{X}_0$. The arithmetic 
model $\mathcal{X}_n$ on which $\overline{\mathcal{M}}_n$ is defined is cofinal with respect to the following compatible morphisms:
\[
    \qquad p_{m,n}:\mathcal X_m\longrightarrow\mathcal X_n
    \quad(m\geq n),\quad p_{m}:=p_{m,0}.
\]
We use the same notation for the induced morphisms on generic fibers. By the 
definition of adelic line bundles, we can choose a sequence of model adelic 
divisors $\{\overline{\mathcal{P}}_n\}$  representing $\overline P$ such that
\begin{equation}\label{eq:smR-boundary-Cauchy}
    -\varepsilon_n\overline{\mathcal E}_m
    \leq\overline{\mathcal P}_m-p_{m,n}^*\overline{\mathcal P}_n
    \leq\varepsilon_n\overline{\mathcal E}_m
    ~~ \text{for any}~m\geq n\geq 0,
\end{equation}
where $\{\varepsilon_n\}$ is a sequence of positive rational numbers converging to $0$, and we set
\[
    \quad
  P_n:=\mathcal P_{n,F}  ,
    \quad \overline{\mathcal E}_n:=p_n^*\overline{\mathcal E}_0,~~\text{and}~~ E_n:=\mathcal E_{n,F},
\]
as before.

Let $X_n:=\mathcal X_{n,F}$ be the generic fiber of $\mathcal{X}_n$.
The generic fiber $X_n$ can be taken smooth by resolution in characteristic
zero. We do not assume regularity of the arithmetic models $\mathcal{X}_n$. Let $p_n: X_n \to X_0$ be the 
birational morphism isomorphic over $V$.

Let $q_Y:Y\to X_0$ be a birational morphism from a normal projective model $Y$ onto $X_0$ which is isomorphic over
$V$. We can choose a compatible family of morphisms
$$
p_{n,Y}: X_n \longrightarrow Y 
$$
such that $p_{n,Y} \circ p_{m,n}=p_{m,Y}$ for any $m\ge n>N$, where $N$ is a sufficiently large integer.

\begin{definition}
    The \textit{geometric b-divisor} $b(\overline{P})=(P_Y)_Y$ \textit{associated with} $\overline{P}$ is defined by
\begin{equation}\label{limit-divisor}
    P_Y:=\lim_{n\to\infty}(p_{n,Y})_*P_n,
\end{equation}
for all proper normal modifications $Y$ of $X_0$ which are isomorphic over $V$.
\end{definition}
Each $P_Y$ is a Weil $\mathbb{R}$-divisor as
a coefficient-wise limit of push-forwards of Weil $\mathbb{Q}$-divisors. It is well defined since on a fixed model $Y$, there are only 
finitely many boundary supported prime divisors, and
\eqref{eq:smR-boundary-Cauchy} makes their coefficients Cauchy.
The same definition applies to any $\overline{P}$ being adelic
$\mathbb{Q}$-divisors or adelic $\mathbb{R}$-divisors.

\begin{remark}
In \cite{BBGHdJ21} and \cite{BGK24}, they relate a b-divisor $D_\pi$ to an adelic divisor $\overline{D}$ under a map 
\begin{align*}
b:\widehat{\mathrm{Div}}(V/\mathbb{Z})~ &\longrightarrow ~ \text{W-b-}\Div_{\mathbb{R}}(X_0)    \\
                    \overline{D} ~~~~~             ~ &\longmapsto ~~~~~ (D_\pi)_\pi
\end{align*}
from the group $\widehat{\mathrm{div}}(V/\mathbb{Z})$ of 
adelic divisors on $V$ to the group $\text{W-b-}\Div(V)$ of  
Weil b-divisors on $V$, where $\pi:Y\to X_0$ runs over
all possible proper modifications of $X_0$.
In \cite[\S3.1]{BGK24}, the group of Weil b-divisors is defined as
$$
\text{W-b-}\Div_{\mathbb{R}}(X_0):= \varprojlim_{q_Y:Y\to X_0} \Div_{\mathbb{R}}(Y),
$$
where the limit is taken with respect to all possible proper modifications $Y$ of $X_0$. This map $b$ is 
well-defined
by \cite[Lemma 3.9]{BGK24}. And for model adelic divisors $\overline{P}$, the representative of 
$b(\overline{P})$ over a specific model $Y$ is exactly the underlying divisor of $\overline{P}$ over $Y$.
Therefore, their definitions coincide with ours in \eqref{limit-divisor}.
\end{remark}

Suppose that $\overline{\mathcal P}_n$ on
$\mathcal X_n$ represents $\overline P$ and $\mathcal Y$ is a fixed
normal projective arithmetic model of $\mathcal{X}_0$. We can similarly define 
the \textit{b-divisor} $(P_{\mathcal{Y}})_{\mathcal{Y}}$ associated with 
$\overline{P}$ as a family of Weil $\mathbb{R}$-divisors with
\begin{equation}\label{arithmetic-limit}
    P_{\mathcal Y}
       :=\lim_{n\to\infty}(p_{n,\mathcal Y})_*\mathcal P_n,
\end{equation}
where $p_{n,\mathcal Y}:\mathcal{X}_n \to \mathcal{Y}$ are compatible so that 
$p_{n,\mathcal Y} \circ p_{m,n}=p_{m,\mathcal Y}$ for sufficiently large $m$ and $n$.
Both constructions for projective models $Y$ over $\Spec F$ and $\mathcal{Y}$ 
over $\Spec \mathbb{Z}$ are independent
of the representatives and compatible with push-forwards via proper birational
morphisms $\mathcal{Y}^\prime \to \mathcal{Y}$.

\begin{lemma}
\label{lem-limit divisor}
Let $V$, $X_n$, $p_{n}$, $\overline M$, $\overline N_0$ be as above and $\overline{P}$ be an adelic divisor representing 
$\overline M-\pi_V^*\overline N_0$. Assume that $X_n$ is smooth for all $n
\ge 0$. We set $P_{X_0}:= \lim_{n\to \infty}{p_{n,*}}P_n$. Then $P_{X_0}$ is a
numerically trivial $\mathbb R$-Cartier divisor on $X_0$. For every smooth
projective model $q_Y:Y\to X_0$ which is isomorphic over $V$, one has
\begin{equation}\label{eq:smR-trace-pullback}
    P_Y=q_Y^*P_{X_0}.
\end{equation}
Moreover, on $X_n$, we have
\begin{equation}
    -\varepsilon_nE_n\leq P_n-p_n^*P_{X_0}\leq\varepsilon_nE_n.
\end{equation}
\end{lemma}

\begin{proof}
Let $B_1,\ldots,B_{d-1}$ be $d-1$ Cartier divisors on
$Y$. The projection formula and the first equality in
\eqref{hypotheses} imply
\[
    P_Y\cdot B_1\cdots B_{d-1}
    =\lim_{n\to \infty} P_n\cdot p_{n,Y}^*B_1\cdots p_{n,Y}^*B_{d-1}=0.
\]
Since $Y$ is smooth, $P_Y$ is a Cartier $\mathbb{R}$-divisor. Applying the classical Hodge
index theorem to $P_Y$ by choosing ample $B_i$, we see that $P_Y$ is
numerically trivial for all $Y$.

In particular, $P_{X_0}$ is numerically trivial. Then the divisor
$P_Y-q_Y^*P_{X_0}$ is exceptional and numerically trivial on $Y$. Therefore, both $P_Y-q_Y^*P_{X_0}$ and its negative are $q_Y$-nef. The well-known
negativity lemma in \cite[Lemma~3.39]{smR-KM98} states that if $-B$ is a $q_Y$-nef Cartier 
$\mathbb{R}$-divisor on $Y$, then $B$ is effective iff $q_{Y,*}B$ is effective.  It is easy to 
see that $q_{Y,*}(P_Y-q_Y^*P_{X_0})=0$, hence the negativity lemma implies that $P_Y-q_Y^*P_{X_0}$ and 
its negative are both effective. We conclude that $P_Y-q_Y^*P_{X_0}=0$.

Taking the push-forward of the geometric part of \eqref{eq:smR-boundary-Cauchy}, we get 
\begin{equation} \label{geom-push}
-\varepsilon_nE_n\leq (p_{m,n})_{*}P_m-P_n\leq\varepsilon_nE_n.
\end{equation}
As $m\to \infty$, the Cartier $\mathbb{R}$-divisor $(p_{m,n})_{*}P_m$ converges to $P_{X_n}=p_n^*P_{X_0}$ by \eqref{eq:smR-trace-pullback}. Thus, we are done.
\end{proof}

  By abuse of notation, we also write $p_{m,n}:\mathcal{X}_m \to \mathcal{X}_n$ 
  for the dominant morphism between the arithmetic models of $V$. Similarly as in \eqref{geom-push}, we have
$$
-\varepsilon_n \mathcal{E}_n\leq (p_{m,n})_{*} \mathcal{P}_m-\mathcal{P}_n\leq\varepsilon_n \mathcal{E}_n.
$$
As $m\to \infty$, this implies that 
\begin{equation}\label{arith-limit-difference}
-\varepsilon_n \mathcal E_n\leq \mathcal P_n- P_{\mathcal{X}_n}\leq\varepsilon_n \mathcal E_n.
\end{equation}

\subsection{Vertical case}\label{chap5-vertical}
First we are going to prove 
Theorem \ref{Rdr} for the vertical adelic line bundle $\overline{M}$, following the ideas in $\S$\ref{Sec-Modify} and $\S$\ref{explicit}.
\begin{definition}
    We say that an adelic divisor $\overline{D}$ on a quasi-projective variety $V$ over a 
    finitely generated field $F$ is \emph{vertical} if $\overline{D}$ can be 
    represented by a 
sequence of model adelic divisors $\{\overline{\mathcal{D}}_n\}$ such that $\mathcal{D}_{n,F}=0$.

We say that an adelic line bundle $\overline{M}$ on $V$ is \emph{vertical} if it 
can be represented by a vertical adelic divisor. 
\end{definition}

For a vertical adelic divisor $\overline{D}$ on $V$, we have $D_Y=0$ 
for any normal projective model $Y$ of $V$. We will prove the fiberwise constancy of its Green function in the 
following proposition for such $\overline{D}$ which also satisfies orthogonality. 

\begin{proposition}
\label{vertical-case-Green-const}
Let $F$ be a finitely generated 
field of transcendence degree $r$ over $\mathbb{Q}$, and let $V$ be a normal, 
geometrically integral quasi-projective variety over 
$F$ of dimension $d\ge 2$. Suppose that $\{\mathcal{X}_n\}_{n\ge 0}$ is a cofinal sequence of 
projective arithmetic models of $V$. Let $\overline P$ be
an integrable adelic $\mathbb{R}$-divisor
on $V$ such that $P_Y=0$ for any normal projective model $Y$ of 
$X_0=\mathcal{X}_{0,F}$.
Assume that for every vertical model adelic line bundle $\overline T$, every model adelic line bundle $\overline A$ on $V$,
and every tuple of arithmetically ample hermitian line bundles
$\overline H_1,\ldots,\overline H_r$ on an arithmetic model of $\Spec F$, we have
\begin{equation}\label{eq-orthogonal}
    \mathcal O(\overline P)\cdot\overline T\cdot\overline A^{d-1}
       \cdot\pi_V^*\overline H_1\cdots\pi_V^*\overline H_r=0,
\end{equation}
where $\pi_V:V\to \Spec F$ is the structure morphism.
Then $\overline{P}$ is vertical, and the Green function of $\overline P$ is fiberwise constant on 
$\mathcal{V}(\mathbb{C})\subset \mathcal{X}_0(\mathbb{C})$, where 
$\mathcal{V}\subset \mathcal{X}_0$ and $\mathcal V\to \mathcal B^\circ$ is a quasi-projective 
arithmetic model of $V/F$.
\end{proposition}

\begin{proof}
By definition, we can choose model adelic divisors $\overline{\mathcal P}_n$ on 
normal projective arithmetic models $\mathcal X_n$, with morphisms $p_{m,n}:\mathcal X_m\to\mathcal X_n$ for $m\geq n$, such that
\begin{equation*}
     -\varepsilon_n\overline{\mathcal E}_m
    \leq\overline{\mathcal P}_m-p_{m,n}^*\overline{\mathcal P}_n
    \leq\varepsilon_n\overline{\mathcal E}_m
    ~~ \text{for any}~m\geq n\geq 0,
\end{equation*}
where $\overline{\mathcal E}_0$ is a strictly effective boundary divisor on
$\mathcal X_0$, and we put
\[
    p_n:=p_{n,0},\qquad
    \overline{\mathcal E}_n:=p_n^*\overline{\mathcal E}_0,
    ~~\text{and}~~ E_n:=\mathcal E_{n,F}.
\]
Taking the push-forward of the geometric part to $X_n$, we get
$$
-\varepsilon_nE_n\leq (p_{m,n})_{*}P_m-P_n\leq\varepsilon_nE_n,
$$
as in \eqref{geom-push}. Letting $m\to \infty$, we obtain 
\begin{equation}\label{eq:smR-zero-trace-bound}
    -\varepsilon_nE_n\leq P_{n}\leq\varepsilon_nE_n,
\end{equation}
since $P_{X_n}=0$ for all $n$. Hence, 
the underlying divisor of $\overline{P}$ on $V$ is zero and 
$\overline{P}$ is vertical. By shrinking a normal arithmetic model 
$\mathcal{X}_0 \to \mathcal{B}$ of $V/F$, we can therefore choose a smooth
quasi-projective model
\[
    \pi:\mathcal V\longrightarrow\mathcal B^\circ
\]
on which the underlying divisor of $\overline P$ is zero since smoothness is an open property on the target. We may
assume that $\mathcal B^\circ(\mathbb C)$ is smooth. We can assume that
$\mathcal{X}_{0,\mathbb{Q}}$ has smooth geometrically integral compact fibers over $\mathcal{B}_{\mathbb{Q}}$,
with a relative simple normal crossings boundary. The generic fibers of $\mathcal{X}_n$ over $F$, and
their fibers over $\mathbb Q$, can be made smooth by generic resolution for each $n\geq 0$.

Now, we choose an arithmetically ample hermitian
line bundle $\overline{\mathcal A}$ on $\mathcal{X}_0$ endowed with a smooth strictly positive metric. Denote the
corresponding model adelic line bundle by $\overline{A}$ and 
its Chern form by $c_1(\overline{A})$.
We also choose an arithmetically ample hermitian bundle $\overline H$
on the base, with Chern form $c_1(\overline{H})$. For every conjugation-invariant real
function
\[
    f\in \mathscr C_c^\infty(\mathcal V(\mathbb C),\mathbb R),
\]
we can extend $f$ by zero to $\mathcal{X}_0(\mathbb{C})\supset \mathcal{V}(\mathbb{C})$. By endowing the trivial line bundle on $\mathcal{X}_0$ with
the metric $\|1\|=e^{-f}$, we get a vertical model adelic line bundle 
$\overline{\mathcal O(f)}$. We denote by $g$ the Green function of $\overline P$ 
on $\mathcal{V}(\mathbb{C})$.
Letting $\overline T=\overline{\mathcal O(f)}$ in
\eqref{eq-orthogonal}, we have
\begin{equation*}
    \int_{\mathcal V(\mathbb C)}
       f\,dd^c g\wedge c_1(\overline{A})^{d-1}\wedge\pi^* c_1(\overline{H})^r=0.
\end{equation*}
By the Stokes formula, this is equivalent to
\begin{equation} \label{eq:smR-current-test}
\int_{\mathcal V(\mathbb C)}
       g\,dd^c f\wedge c_1(\overline{A})^{d-1}\wedge\pi^* c_1(\overline{H})^r=0,
\end{equation}
since both $c_1(\overline{A})$ and $\pi^* c_1(\overline{H})$ are $\partial$-closed and 
$\bar \partial$-closed on $\mathcal{V}(\mathbb{C})$.

Let $g_n$ be the Green function of
$\overline{\mathcal P}_n$ on
$\mathcal V(\mathbb C)$ where we identify $\mathcal V(\mathbb C)$ as an open subset of every $\mathcal{X}_n(\mathbb{C})$.
For $w\in\mathcal B^\circ(\mathbb C)$, we write $\mathcal V_w(\mathbb{C})$ for the complex fiber of $\mathcal V(\mathbb C)$ over $w$
and $g_w$, $g_{n,w}$, $c_1(\overline{A})_w$ for restrictions of
$g$, $g_n$, $c_1(\overline{A})$ to $\mathcal{V}_w(\mathbb{C})$ respectively. We claim that on $\mathcal V_w(\mathbb{C})$, 
\[
    dd^cg_w\wedge c_1(\overline{A})_w^{d-1}=0
\]
as a current for almost every $w\in \mathcal{B}^\circ(\mathbb C)$. To see this, we fix a point $p\in \mathcal V(\mathbb C)$, and take a neighborhood of $p$ of the form $U\times U^\prime$, with local holomorphic charts
$$
(z,w)=(z_1,\dots,z_d,w_1,\dots,w_r)
$$
and $U,U^\prime$ biholomorphic to complex disks of dimensions $d$ and $r$ respectively. 
We then choose smooth functions $\phi\in \mathscr C_c(U)$ and $\psi\in \mathscr C_c(U^\prime)$ of compact support and take
$$
f(z,w)=\phi(z)\psi(w)
$$
in \eqref{eq:smR-current-test} to get
\begin{align*}
&\int_{U\times U^\prime}g(z,w)\psi(w) \,dd_z^c \phi(z)\wedge c_1(\overline{A})_w^{d-1}\wedge\pi^* c_1(\overline{H})^r \\
    =&\int_{U^\prime} \left(\int_U g(z,w)\,dd_z^c \phi(z)\wedge c_1(\overline{A})_w^{d-1} \right) \psi(w)\,c_1(\overline{H})^r=0.
\end{align*}
Here, we use the notation $dd_z^c$ to denote the operator $dd^c$ on $U$ with respect to the local coordinate $z$.
As $\psi$ runs over all compactly supported smooth functions on $U^\prime$, we have
$$
\int_U g(z,w)\,dd_z^c \phi(z)\wedge c_1(\overline{A})_w^{d-1} =0, ~\text{for almost every }w\in U^\prime.
$$
Moreover, this integral is continuous in $w$, which implies that
$$
\int_U g_w(z)\,dd_z^c \phi(z)\wedge c_1(\overline{A})_w^{d-1} =0, ~\text{for every }w\in U^\prime.
$$
Applying the Stokes formula again, we see that
$$
\int_U  \phi(z)\,dd_z^c g_w(z) \wedge c_1(\overline{A})_w^{d-1} =0.
$$
Since $\phi$ is arbitrary, this proves the claim that the current 
$dd_z^cg_w\wedge c_1(\overline{A})_w^{d-1}=0$ for every $w\in \mathcal{B}^\circ(\mathbb C)$ by varying $p$ and $U^\prime \subset \mathcal{B}^\circ(\mathbb C)$.

Next, we will show that the continuous function $g_w$ is indeed smooth for all $w$. Choose local unitary frames $dz_1,\dots,dz_d$ such that
$$
c_1(\overline{A})_w=i\sum_{1\leq j,\,k\leq d}h_{jk}\,{dz_j \wedge d \bar{z}_k},
$$
where $H_w=(h_{jk})$ is a positive definite Hermitian matrix with smooth entries.
For any smooth function $u$ on $\mathcal{V}_w(\mathbb{ C})$, we denote its Hessian matrix by
$$
H(u)= \left( \frac{\partial^2u}{\partial z_j\partial\bar z_k}\right)_{1\le j,k\le d}.
$$
Therefore, we have
$$
dd^c u \wedge c_1(\overline{A})_w^{d-1}=\frac{1}{d} \mathrm{tr}(H_w^{-1}H(u)) c_1(\overline{A})_w^{d}.
$$
This defines an elliptic operator $L_w:u \mapsto \mathrm{tr}(H_w^{-1}H(u))$. By the elliptic 
regularity result in \cite[Theorem 2.35]{FRRO22}, we deduce that $g_w$ is a smooth function from the equation
$L_w(g_w)=0$. Therefore, the current $dd^c g_w$ is indeed a smooth $(1,1)$-form on $\mathcal{V}_w(\mathbb{C})$. 

Then we are going to show $dd^c g_w=0$. For any smooth function $h$ on $\mathcal{V}(\mathbb{C})$, we can modify the metric on $\overline{A}=(A,\| \cdot\|)$ to get $\overline{A}_h=(A,\|\cdot\|_h)$ with metric
$$
\| \cdot\|_h=e^{-h}\| \cdot\|.
$$
By the assumption on $\overline{P}$, we have
$$
dd^cg_w\wedge c_1(\overline{A}_h)_w^{d-1}=0,
$$
for all $h$ such that $\overline{A}_h$ is strictly semipositive. We can find a smooth $h$ extending $g_w$, then we deduce that
\begin{equation}\label{ddcg}
dd^cg_w\wedge(c_1(\overline{A})_w+t\,dd^c g_w)^{d-1}=0,    
\end{equation}
for sufficiently small $t\in \mathbb{R}$.

Fix a point 
$x\in \mathcal{V}_w(\mathbb{C})$ and a local frame $dz_1,\dots,dz_d$ of sheaf of 
differentials near $x$. The $dz_i$ are chosen so that for a fixed $\overline{A}$, we have 
$$
c_1(\overline{A})_w=i\sum_{j=1}^d{dz_j \wedge d \bar{z}_j},~~\text{and}~~ dd^cg_w=i\sum_{1\le j,~k \le d} \lambda_{jk}\,{dz_j \wedge d \bar{z}_k},
$$
at $x$ with $\lambda_{jk}\in \mathbb{C}$. Considering the coefficient of $t$ in \eqref{ddcg}, we see that
$$
-\frac{1}{d} \sum_{j,~k}\left|\lambda_{jk}\right|^2=0.
$$
For the detailed calculation, we refer to the proof in \cite[Lemma 1.1.2]{Mor96}.
Thus, we have shown that
\begin{equation}\label{eq:smR-pluriharmonic}
    dd^cg_w=0,
\end{equation}
on every fiber $\mathcal{V}_w(\mathbb{C})$ of $\mathcal{V}(\mathbb{C})\to \mathcal{B}(\mathbb C)$.
The function $g_w$ is 
pluriharmonic on $\mathcal{V}_w(\mathbb{C})$ for each $w$.

To show that $g_w$ is a constant,
we next study the behavior of $g_w$ near $\mathcal{X}_0(\mathbb{C})\setminus \mathcal{V}(\mathbb{C})$. Let $g_{E,n}$ be the Green function of
$\overline{\mathcal E}_n$. 
The Cauchy condition on $\overline{\mathcal P}_n$ shows that $g_n\to g$ locally uniformly on 
$\mathcal{V}(\mathbb{C})$. By \eqref{eq:smR-zero-trace-bound}, the $\mathbb{R}$-divisors $\varepsilon_nE_n\pm P_{n}$ are effective. 
Their Green functions 
 on $\mathcal{V}_w(\mathbb{C})$ are bounded below. Hence for each $n$ and $w\in \mathcal{B}^\circ (\mathbb{C})$, there is a constant $c_{n,w}\in \mathbb{R}$ such that
\[
    |g_{n,w}|\leq\varepsilon_ng_{E,w}+c_{n,w}
\]
on $\mathcal{X}_{n,w}(\mathbb{C})$, where $g_{E,w}$ is the restriction of $g_E$ to $\mathcal{X}_{n,w}(\mathbb{C})$.
We also have
\[
    |g_{n,w}-g_w|\leq\varepsilon_ng_{E,w}
\]
on $\mathcal{V}_w(\mathbb{C})$. Consequently, for every positive number $\lambda>0$, there is a constant $c_{\lambda,w}\in\mathbb R$ such that
\[
    |g_w|\leq\lambda g_{E,w}+c_{\lambda,w}.
\]
 Since the boundary in the fiber $\mathcal{X}_{0,w}$ is a simple normal crossing divisor, the Green
functions $g_{E,w}$ is bounded above by a positive linear combination
of $-\log|z_j|$, where $z_1,\dots,z_{d}$ are local coordinates near 
the boundary $\mathcal{X}_{0,w}(\mathbb{C})\setminus \mathcal{V}_w(\mathbb{C})$. It is a 
standard result to show that the pluriharmonic function $g_w$ bounded by $\lambda$ times a 
logarithmic growth near boundary for every $\lambda >0$ is indeed a constant. And we include 
the proof in Lemma~\ref{lem:smR-removability}.

\end{proof}

\begin{lemma}[Removability]
\label{lem:smR-removability}
Let $X$ be a connected compact complex manifold, and $B\subset X$
be a simple normal crossings divisor. Suppose $u$ is a real-valued
pluriharmonic function on $X\setminus B$. Suppose that on every local
coordinate chart on which $B$ is given by equation
$z_1\cdots z_k=0$, and for every $\lambda>0$, there is a constant
$c_\lambda\in \mathbb{R}$ such that
\[
    |u(z)|\leq -\lambda \sum_{j=1}^k\log|z_j|+c_\lambda.
\]
Then $u$ extends to a pluriharmonic function on $X$ and is thereby a constant.
\end{lemma}

\begin{proof}
We work in a local coordinate polydisc with $|z_j|<1$ and put
$\psi:=\sum_{1\le j \le k}\log|z_j|$. For each $\lambda>0$, the functions
$u+\lambda\psi$ and $-u+\lambda\psi$ are plurisubharmonic outside $B$
and locally bounded above near $B$. Applying the result of removable singularities in \cite[Chapter~I, 
Theorem~5.24]{smR-Demailly}, we see that $\pm u +\lambda\psi$ extend uniquely to plurisubharmonic functions across $B$. 
The growth assumption on $u$ also implies that $u$ is locally integrable. Since
$dd^c\psi$ is a positive current supported on $B$, we have 
\[
    -\lambda\, dd^c\psi\leq dd^cu\leq\lambda\, dd^c\psi,
\]
in the sense of currents on $X$. Let $\tilde{u}$ be the limit of $u+\lambda\psi$ in $L^1_{\mathrm{loc}}(X)$. Then $dd^c\tilde{u}=0$ as a current on $X$.  
Applying the elliptic regularity in \cite[Theorem 2.35]{FRRO22}, we see that $u$ 
can be extended  pluriharmonically to all of $X$. Since $X$ is a connected
compact complex manifold, the function $u$ must be a constant.
\end{proof}

\begin{corollary} \label{arithmetic-limit-0}
    Let $V/F$ and $\overline{P}$ be as in Proposition \ref{vertical-case-Green-const}. Fix an 
    $x_0\in V(\overline{F})$. Suppose further that the adelic line bundle $\overline{M}^\prime$ 
    associated with $\overline{P}$ satisfies 
    $$
    \mathfrak{h}_{\overline{M}^\prime}(x_0)=0.
    $$
   Then the Green function $g$ of $\overline{P}$ is identically zero, and $P_{\mathcal Y}$ defined in \eqref{arithmetic-limit} has zero coefficient at every prime
divisor horizontal over $\mathbb{Z}$, on every normal projective
arithmetic model $\mathcal Y$.
\end{corollary}
\begin{proof}
    Let $\mathcal{V}\to \mathcal{B}^\circ$ be a quasi-projective arithmetic model of $V/F$.
    By applying Proposition \ref{vertical-case-Green-const} to $\overline{P}$, 
    its Green function $g$ is fiberwise constant on $\mathcal V(\mathbb{C})$.
    By Lemma \ref{zero-g}, this fiberwise constant must be zero for every $w\in \mathcal  B^\circ(\mathbb{C})$ since 
    $\mathfrak{h}_{\overline{M}^\prime}(x_0)=0$. Consequently, we have $g\equiv0$.

Next, we fix a prime divisor $\mathcal{W}$ horizontal over $\Spec \mathbb Z$ on a
normal projective arithmetic model $\mathcal Y$ of $V/F$. Pass to common
models dominating $\mathcal Y$. Let $\{\mathcal{X}_n\}_{n\geq 0}$ be a cofinal sequence of normal 
projective arithmetic models of $\mathcal{V}$ such that $\overline{P}$ is represented by model divisors $\{\overline{\mathcal{P}}_n\}_{n\ge 0}$. We write $\mathcal{W}_n$ for the strict
transform of $\mathcal{W}$ on $\mathcal X_n$. 

 We let $g_n$ and $g_{\overline{\mathcal{E}}_n}$ denote the Green functions of $\overline{\mathcal{P}}_n$ and boundary 
 divisors $\overline{\mathcal{E}}_n$ on $\mathcal{X}_n(\mathbb{C})$ respectively. We have    
 $$
 |g_n|\leq\varepsilon_n \,g_{\overline{\mathcal{E}}_n}
 $$
 on $\mathcal{V}(\mathbb{C})$. If $\mathcal{W}_n$ meets $\mathcal{V}$, then for all $n\ge 0$ we have
 $$
 \operatorname{ord}_{\mathcal{W}_n}(\mathcal P_n)=0.
 $$
 If $\mathcal{W}_n$ is contained in $\mathcal{X}_n\setminus \mathcal{V}$, then
 by comparison of the logarithmic singularities of $g_n$ and $g_{\overline{\mathcal{E}}_n}$ near a general point of $\mathcal{W}_n(\mathbb{C})$, we obtain
\[
    |\operatorname{ord}_{\mathcal{W}_n}(\mathcal P_n)|
       \leq\varepsilon_n \operatorname{ord}_{\mathcal{W}_n}(\mathcal E_n).
\]
Taking limit of push-forwards of $\mathcal P_n$ to $\mathcal{Y}$, we deduce that
$\operatorname{ord}_{\mathcal{W}}(P_{\mathcal Y})=0$, as required.
\end{proof}

 Now we need to deal with the vertical part of 
 $P_\mathcal{X}$.  Let $\mathcal{X} \to \mathcal{B}$ be a normal 
 projective arithmetic model of $V/F$. For any vertical prime divisor 
 $Z\subset \mathcal B$, let $\eta_Z$ be its generic point. Then the local ring 
 $\mathcal{O}_{\mathcal{B},\eta_Z}$ is an excellent discrete valuation ring. Let $K$ be its fraction field and $k$ be its residue field.
 
Since $\mathcal{X}$ is generally non-smooth, a Weil divisor on 
$\mathcal{X}$ does not always define a Cartier divisor. To derive more 
numerical information about $P_\mathcal{X}$, we will apply a resolution result in \cite[Theorem 2.1]{Hau15} to the localization of $\mathcal{X}$. If $\mathscr X$ is
a reduced quasi-excellent scheme, there is a projective birational
morphism $f:\mathscr Y\to\mathscr X$ which is an isomorphism over
$\mathscr X_{\mathrm{reg}}$, and $\mathscr Y$ is regular at every
point whose image has codimension at most two in $\mathscr X$.
When $\mathscr X$ is integral, we
may assume that $\mathscr Y$ is normal by normalization.  

More generally, for a given excellent discrete valuation
ring $R$, let $K=\Frac (R)$ be its fraction field and let $k$ be its residue field.  A vertical divisor on an $R$-scheme
is a divisor supported on its special fiber. Then we have the following lemma.

\begin{lemma}
\label{lem:reduce-to-surface}
Let $R$ be an excellent discrete valuation
ring, $K$ be its fraction field, and $k$ be its residue field. Assume
that $\operatorname{char}K=0$. Let $\mathscr X$ be a normal integral projective
flat $R$-scheme with smooth geometrically integral generic fiber
of dimension $d\geq 2$. Let $A$ be an ample line bundle on
$\mathscr X$. Denote by $\mathscr{X}_s$ the special fiber of 
$\mathscr{X}$. Suppose that $(Q_{\mathscr Z})_{\mathscr Z}$ is the $b$-divisor 
associated with an adelic divisor $\overline{Q}$ on $\mathscr X/\mathbb{Z}$, 
where $\mathscr Z$ runs over birational models of $\mathscr X$ 
over $\Spec R$ with projective birational morphism 
$q_{\mathscr Z}:\mathscr Z\to \mathscr X$. Suppose further that for every $\mathscr Z$, the Weil $\mathbb{R}$-divisor $Q_{\mathscr{Z}}$ is vertical and
\begin{equation}\label{eq:smR-local-orthogonal}
    Q_{\mathscr Z}\cdot \mathcal{T} \cdot(q_{\mathscr Z}^*A)^{d-1}=0
\end{equation}
for every vertical Cartier divisor $\mathcal{T}$ on $\mathscr{Z}$.
Then $ Q_{\mathscr X}$ is a principal $\mathbb R$-divisor, and there exists a $c\in\mathbb R$ such that
\[
    Q_{\mathscr X}=c \cdot \mathscr{X}_s.
\]
\end{lemma}

\begin{proof}
\emph{Step 1:}
We apply \cite[Theorem~2.1]{Hau15} and normalize the resulting scheme to obtain a birational morphism
$f:\mathscr Y\to\mathscr X$, isomorphic over $\mathscr X_{\mathrm{reg}}$,
such that $\mathscr Y$ is regular at every point whose image has
codimension at most two in $\mathscr X$. 
We set
\[
 W:=f(\operatorname{Sing}\mathscr Y),
\]
which is the image of the singular locus $\operatorname{Sing}\mathscr Y$ of $\mathscr Y$.
The set $W$ is a closed subset of $\mathscr X_s$. It has
codimension $\geq 3$ in $\mathscr X$. Hence $\dim_kW\leq d-2$,
and $\mathscr U:=f^{-1}(\mathscr X\setminus W)$ is regular.
Write
\begin{equation}\label{Y-Q}
 \mathscr Y_s=\sum_{i=1}^\ell m_iY_i,\quad \mathrm{and} \quad
 Q_{\mathscr Y}=\sum_{i=1}^\ell c_iY_i,
\end{equation}
where $m_i>0$ is the multiplicity of the irreducible component $Y_i$ of $\mathscr Y_s$  and $c_i\in \mathbb{R}$.

\emph{Step 2:}
For each $i$, let
$\mathcal I_i:=\mathcal O_{\mathscr Y}(-Y_i)\subset\mathcal O_{\mathscr Y}$
be the ideal sheaf defining $Y_i$.
Set $\mathscr Z_0=\mathscr Y$ and $h_0=\operatorname{id}_{\mathscr Y}$.
We can inductively blow up the ideal
\[
 \mathcal J_i:=\mathcal I_i\mathcal O_{\mathscr Z_{i-1}}
\]
to get normalized blowups
$$
 b_i:\mathscr Z_i:=(\operatorname{Bl}_{\mathcal J_i}
                   \mathscr Z_{i-1})^\nu
                  \longrightarrow\mathscr Z_{i-1},
$$
where $(\operatorname{Bl}_{\mathcal J_i} \mathscr Z_{i-1})^\nu$ denotes the normalization of $\operatorname{Bl}_{\mathcal J_i} \mathscr Z_{i-1}$. 
We denote by $h_i$ the resulting birational morphism $\mathscr{Z}_i \to\mathscr{Z}_0$.
Set $\mathscr Z:=\mathscr Z_\ell$, $h:=h_\ell$, and $q:=f\circ h$.
 Since each $\mathcal I_i$
restricts to the unit ideal on the generic fiber, these morphisms $b_i$ are isomorphisms over generic fibers. We see that each ideal sheaf $\mathcal I_i\mathcal O_{\mathscr Z}$ on 
$\mathscr Z$ is
invertible, therefore inducing an effective vertical Cartier divisor
$\mathcal T_i$ by
\[
 \mathcal I_i\mathcal O_{\mathscr Z}
       =\mathcal O_{\mathscr Z}(-\mathcal T_i).
\]
On the regular scheme $\mathscr U$, all $Y_i$ are Cartier. The hypothesis on $Q_{\mathscr Z}$ implies
\begin{equation}\label{local-intsct}
 Q_{\mathscr Z}\cdot\mathcal T_i\cdot(q_{\mathscr Z}^*A)^{d-1}=0,
\end{equation}
for $1\leq i\leq \ell$. 

\emph{Step 3:} Let $\varpi\in R$ be a uniformizer.
We can inductively choose $d-1$ sections $s_1,\dots,s_{d-1} \in H^0(\mathscr{X}_s,A|_{\mathscr{X}_s}^{\otimes N})$ and denote by $V(s_i)$ the zero locus of $s_i$ in $\mathscr{X}_s$,
so that $V(s_i)$ meets properly with each other and
$$
V(s_1,\dots,s_{d-1})=V(s_1)\cap V(s_2) \cap \cdots\cap V(s_{d-1})
$$
contains no irreducible component of any closed non-empty subsets in $\mathscr{X}_s$ that are images of the following:
\begin{enumerate}[$\bullet$]
    \item  $W\subset \mathscr X_s$;
    \item  the irreducible components $Y_i \subset \mathscr Y_s$ for $1\leq i \leq \ell$;
    \item   the irreducible components $Y_i\cap Y_j \subset \mathscr Y_s$ for $i\neq j$;
    \item  the irreducible components of $\mathscr Z_s$.
\end{enumerate}
This is possible because $A^{\otimes N}$ is very ample and the natural morphism
$$H^0(\mathscr{X},A^{\otimes N}) \longrightarrow H^0(\mathscr{X}_s,A|_{\mathscr{X}_s}^{\otimes N})$$
is surjective. 
Since $\dim_k W\leq d-2$ and $\dim_k \mathscr Y_s \leq d$, we have $W \cap V(s_1,\dots,s_{d-1})= \varnothing$, 
and $\dim(Y_i\cap Y_j\cap f^{-1} V(s_1,\dots,s_{d-1}))\leq 0$ for generic choices of $s_1,\dots,s_{d-1}$.

Now we observe that $H^0(\mathscr{X},A^{\otimes N})\otimes_R K=H^0(\mathscr{X}_K,A_K^{\otimes N})$. 
Let $\tilde s_1,\ldots, \tilde s_{d-1}$ be a tuple of lifts of $s_1,\dots,s_{d-1}$.
Then all possible lifts $s_1^\prime,\dots,s_{d-1}^\prime$ of $s_1,\dots,s_{d-1}$ form a coset
$$
(\tilde s_1, \dots,\tilde s_{d-1})+\varpi \cdot  H^0(\mathscr{X},A^{\otimes N})^{d-1},
$$
which is Zariski dense in the $K$-linear space $H^0(\mathscr{X}_K,A_K^{\otimes N})^{d-1}$. Recall that Bertini's 
theorems ensure that there exists a non-empty open subset of the parameter space defining $s_i$ that cut out smooth 
geometrically integral complete intersection curves on $\mathscr X_K$. Therefore, we can choose appropriate lifts so that
$$
V(s_1^\prime,\dots,s_{d-1}^\prime)_K
$$
is a smooth geometrically integral curve over $K$ with the scheme-theoretic intersection structure. Let
$$
\mathscr S:=f^{-1}(V(s_1^\prime,\dots,s_{d-1}^\prime)) \subset \mathscr Y.
$$
Since $f_K:\mathscr Y_K\to\mathscr X_K$ is an isomorphism of generic fibers, we see that $\mathscr S_K$ is a smooth
geometrically integral curve which is contained in $\mathscr U_K$.

\emph{Step 4: } Finally, we can write \eqref{local-intsct} as an intersection product on $\mathscr S$.
Since $ \mathscr S\subset \mathscr U$, we define
$$
S_i:=Y_i|_{\mathscr S}=\mathcal T_i|_{\mathscr S}.
$$
Then we have 
$$
\mathscr S_s=  \sum_{i=1}^\ell m_iS_i,~~~\mathrm{and} ~~~Q_{\mathscr Y}|_{\mathscr S}= \sum_{i=1}^\ell c_iS_i.
$$ 
Each $S_i$ is an effective Cartier divisor contained in $\mathscr S_s$ or zero.
Furthermore, $S_i \neq 0$ whenever $f(Y_i)$ is a component of $\mathscr X_s$, i.e., $Y_i$ is non-exceptional.
Then we have
\[
 Q_{\mathscr Y}|_{\mathscr S} \cdot S_i
  =N^{d-1}\bigl(Q_{\mathscr Z}\cdot\mathcal T_i
                         \cdot(q^*A)^{d-1}\bigr)=0,
\]
by projection formula since $\mathscr S$ is the inverse image of intersections of $d-1$ hypersurfaces 
induced by sections of $A^{\otimes N}$.
Therefore, we have
\[
 Q_{\mathscr Y}|_{\mathscr S}^2=\sum_i c_i(Q_{\mathscr Y}|_{\mathscr S}\cdot S_i)=0.
\]
 Write $m_{ij}:=S_i\cdot S_j\geq 0$.
Also we have $\sum_j m_j(S_i\cdot S_j)=0$, since
$\mathscr S_s=\operatorname{div}_{\mathscr S}(\varpi)$ is principal.
Similarly as the vertical case of curves in \eqref{self-intersection}, we have
\[
 Q_{\mathscr Y}|_{\mathscr S}^2=-\sum_{i<j}m_{ij}\, m_i m_j
 \left(\frac{c_i}{m_i}
             -\frac{c_j}{m_j}\right)^2=0.
\]
Since $\mathscr S_s$ is connected, we see that $\frac{c_1}{m_1}=\cdots=\frac{c_\ell}{m_\ell}=c$ for some $c\in\mathbb R$ and
$D=c\mathscr S_s$.
Because distinct $S_i$ has no common component, we have $c_i=cm_i$ whenever $S_i\neq 0$ in \eqref{Y-Q}, particularly for every
non-exceptional $Y_i$.
Finally from $f_*Q_{\mathscr Y}= Q_{\mathscr X}$ and
$f_*\mathscr Y_s=\mathscr X_s$, we conclude that
$Q_{\mathscr X}=c\mathscr X_s$.
\end{proof}

\begin{remark}\label{rem:smR-finite-tests}
The proof uses vanishing in \eqref{local-intsct} only for the finitely many effective
vertical Cartier divisors $\mathcal T_1,\ldots,\mathcal T_\ell$ on the
single model $q:\mathscr Z\to\mathscr X$.
The models are constructed over $R$, independently of $\overline{Q}$.
\end{remark}

The next lemma relates the intersection product of an adelic divisor $\overline P$ which satisfies orthogonality with a fixed tuple of line bundles to the intersection product of its b-divisor with their underlying line bundles. 

\begin{lemma}
\label{expression-of-arith-intsct}
Let $\overline P$ be an integrable
adelic $\mathbb{R}$-divisor on $V$. Fix a normal projective arithmetic
model $\pi_{\mathcal Y}:\mathcal Y\to\mathcal B$, where
$\mathcal B$ is a normal projective arithmetic model of $F$ and the generic fiber
$\mathcal Y_F$ is a projective compactification of $V$.
Let $\mathcal T$ be a vertical Cartier divisor on $\mathcal Y$ whose image under $\pi_\mathcal{Y}$ is a 
closed point. Denote by $\mathcal O(\overline{\mathcal T})$ the model adelic line bundle obtained by endowing 
$\mathcal O(\mathcal T)$ with the trivial metric for which its
canonical rational section has norm $1$ at infinity.
For any hermitian line bundles $\overline{\mathcal A}$ on $\mathcal Y$ and
$\overline{\mathcal H}$ on $\mathcal B$, we have
\begin{equation}\label{geom-intsct}
\begin{aligned}
    &\mathcal O(\overline P)\cdot\mathcal O(\overline{\mathcal{T}})
       \cdot\overline{\mathcal{A}}^{d-1}\cdot\pi_{\mathcal{Y}}^*\overline{\mathcal{H}}^r\\
    &\qquad=
       P_{\mathcal Y}\cdot \mathcal{T} \cdot\mathcal A^{d-1}
          \cdot\pi_{\mathcal Y}^*\mathcal H^r,
\end{aligned}
\end{equation}
where $\mathcal{A}$ and $\mathcal{H}$ are the underlying line bundles of $\overline{\mathcal A}$ and $\overline{\mathcal H}$ respectively.
\end{lemma}

\begin{proof}
Let 
$\{\overline{\mathcal P}_n\}$ be a sequence of model adelic divisors converging to $\overline{P}$ with respect to a cofinal family of arithmetic models $\mathcal{X}_n$. Let
$p_{n,\mathcal Y}:\mathcal X_n\to\mathcal Y$ be a compatible family of dominant morphisms for 
$n> N$ and $N$ be a sufficiently large integer. Let $p_n:\mathcal{X}_n \to \mathcal{X}_0$ and $p_\mathcal{Y}: \mathcal{Y} \to \mathcal{X}_0$ be the 
dominant morphisms, then by the definition of arithmetic intersection of adelic line bundles, we 
have
\begin{align*}
    \mathcal O(\overline P)\cdot\mathcal O(\overline{\mathcal{T}})
       \cdot \overline{\mathcal{A}}^{d-1}\cdot\pi_{\mathcal{Y}}^*\overline{\mathcal{H}}^r
       &=
       \lim_{n\to \infty} \mathcal{O}(\overline{\mathcal P}_n)\cdot p_{n,\mathcal{Y}}^*\mathcal O(\overline{\mathcal{T}})
       \cdot  p_{n,\mathcal{Y}}^*\overline{\mathcal{A}}^{d-1}\cdot p_{n,\mathcal{Y}}^*\pi_{\mathcal{Y}}^*\overline{\mathcal{H}}^r.
\end{align*}
Since the metric of $\mathcal O(\overline{\mathcal{T}})$ is trivial, the intersection product has no archimedean contribution. Together with the projection formula, it yields
\begin{align*}
    \mathcal{O}(\overline{\mathcal P}_n)\cdot p_{n,\mathcal{Y}}^*\mathcal O(\overline{\mathcal{T}})
       \cdot p_{n,\mathcal{Y}}^*\overline{\mathcal{A}}^{d-1}\cdot p_{n,\mathcal{Y}}^*\pi_{\mathcal{Y}}^*\overline{\mathcal{H}}^r
       &= 
       \mathcal{P}_n \cdot p_{n,\mathcal{Y}}^* \mathcal{T}\cdot  
      p_{n,\mathcal{Y}}^*\mathcal{A}^{d-1}\cdot p_{n,\mathcal{Y}}^*\pi_{\mathcal{Y}}^*\mathcal{H}^r \\
       &=(p_{n,\mathcal{Y}})_*\mathcal{P}_n \cdot \mathcal{T} \cdot \mathcal{A}^{d-1} \cdot \pi_{\mathcal{Y}}^*\mathcal{H}^r. 
\end{align*}
Taking $n\to \infty$, this proves \eqref{geom-intsct}.
\end{proof}

\begin{proposition}
\label{prop:smR-zero-limit}
Let $F$ be a finitely generated 
field of transcendence degree $r$ over $\mathbb{Q}$, and let $V$ be a normal, geometrically integral quasi-projective variety over 
$F$ of dimension $d\ge 2$. Fix an $x_0\in V(\overline{F})$. Let $\overline P$ be an integrable adelic $\mathbb{R}$-divisor on
$V$ associated with an adelic line bundle $\overline{M}^\prime$. Suppose $P_Y=0$ for any normal projective model $Y$ of $V$ and 
$$
\mathfrak{h}_{\overline{M}^\prime}(x_0)=0.
$$
Suppose further that $\overline{P}$ satisfies orthogonality. Then  $\overline P=0$
as an adelic divisor.
\end{proposition}

\begin{proof}
Choose a sequence of model adelic divisors $\overline{\mathcal P}_n$
on $\mathcal X_n$ that converges to $\overline{P}$, 
and a strictly effective boundary divisor
$\overline{\mathcal E}_0$ on $\mathcal X_0$, where $\{\mathcal{X}_n\}$ is a cofinal 
family of projective arithmetic models of $V/F$. Let $\overline{\mathcal E}_n$ be the 
pullback of $\overline{\mathcal E}_0$ to $\mathcal{X}_n$. Let $\mathcal{V} \subset \mathcal{X}_0$ be a normal quasi-projective model of $V/F$. Corollary \ref{arithmetic-limit-0} shows that the Green
function associated with $\overline{P}$ is identically zero on $\mathcal{V}(\mathbb{C})$ and $\mathcal{P}_{\mathcal{Y}}$ 
has zero coefficient at prime divisors horizontal over
$\Spec \mathbb Z$ for every normal projective model of $\mathcal{V}$. It remains to 
consider the coefficients of $\mathcal{P}_{\mathcal{Y}}$ over vertical divisors of 
$\mathcal{Y}$. We write $\mathcal{X}:=\mathcal{X}_0$ and $\pi:\mathcal{X}\to \mathcal{B}$, where $\mathcal{B}$ is an arithmetic model of $\Spec F$.

Now we fix a divisorial valuation $v$ of $F$ such that its residue field $\kappa(v)$ has characteristic
$p>0$. We normalize $v$ by requiring that the value group $v(F^\times)=\mathbb Z$. Choose a normal
projective arithmetic model $\mathcal B$ of $F$ and a prime divisor
$Z$ realizing $v$. Let $\eta_Z$ be its generic point. We set 
\[
    R:=\mathcal O_{\mathcal B,\eta_Z}
\]
as a discrete valuation ring and $ k:=\kappa(Z)$ as its residue field.
We choose a uniformizer $\varpi_v$ of $R$. Let
\[
    \mathscr X:=\mathcal X\times_{\mathcal B}\Spec R.
\]
It is normal, integral, projective and flat over $R$, with smooth
geometrically integral generic fiber. Fix two arithmetically ample
hermitian line bundles $\overline{\mathcal A}$ on $\mathcal X$
and $\overline{\mathcal H}$ on $\mathcal B$. 

We recall the constructions in
Lemma~\ref{lem:reduce-to-surface} on arithmetic models of $\mathscr X$. By applying \cite[Theorem 2.1]{Hau15}
to $\mathcal X$, we obtain a projective birational morphism
$f:\mathcal Y\to\mathcal X$ such that $\mathcal{Y}$ is regular at every 
$y\in \mathcal{Y}$ such that $f(y)$ has codimension $\le 2$ in $\mathcal{X}$. 
The localization over $\eta_Z$ maintains the regularity above the points
of $\mathscr X$ of codimension $\leq 2$. This is because these local rings of $\mathscr{X}$ are local rings of
$\mathcal X$, with the same codimension. The morphism $f$ is an
isomorphism on the smooth generic fiber over $F$.

Let $Y_1,\ldots,Y_{\ell}\subset \mathscr Y$ be the components of the special fiber of
$\mathscr Y:=\mathcal{Y}\times_{\mathcal B}\,\Spec R$. Let $\mathcal{Y}_i$ be their Zariski closure in $\mathcal Y$. Each $\mathcal Y_i$ is a
prime divisor supported above $Z$, hence vertical.
By successively blowing up the pullbacks of the ideals
$\mathcal O_{\mathcal Y}(-\mathcal Y_i)$ and taking normalization, we get a morphism
\[
    h:\mathcal Z\longrightarrow\mathcal Y.
\]
Let $q=f\circ h: \mathcal Z\longrightarrow\mathcal X $ be the composition morphism and 
$\mathcal T_1,\ldots,\mathcal T_\ell$ be the effective Cartier divisors on $\mathcal Z$ determined
by the total transforms of these ideals $\mathcal O_{\mathcal Y}(-\mathcal Y_i)$. All $\mathcal T_i$
are supported above $Z$.

We endow every $\mathcal O(\mathcal T_i)$ with the trivial metric at infinity and write $\mathcal O(\overline{\mathcal T}_i)$ for the resulting adelic line bundle on $\mathcal{Z}$.
Since $\overline{P}$ satisfies orthogonality by assumption, we see that
$$
\mathcal O(\overline P)\cdot \mathcal O(\overline{\mathcal T}_i) \cdot q^* \overline{\mathcal A}^{d-1}
    \cdot\pi_{\mathcal{Z}}^*\overline{\mathcal H}^r=0,
$$
for every $1\le i \le \ell$ by \eqref{eq-orthogonal} and fixed ample $\overline{\mathcal A}$ on $\mathcal X$ and ample $\overline{\mathcal H}$ on $\mathcal{B}$, where 
$\pi_{\mathcal{Z}}=\pi \circ q$ is the dominant morphism $\mathcal{Z}\to \mathcal{B}$. This holds because we can 
approximate $q^* \overline{\mathcal A}$ by arithmetically ample
$q^*\overline{\mathcal A}+\delta\overline{\mathcal B}$
where $\delta\in\mathbb Q_{>0}$ and
a fixed arithmetically ample hermitian line bundle
$\overline{\mathcal B}$ on $\mathcal Z$. By
Lemma~\ref{expression-of-arith-intsct}, we have
\begin{equation}\label{qqq}
P_{\mathcal{Z}} \cdot \mathcal{T}_i \cdot q^*\mathcal A^{d-1}
          \cdot\pi_{\mathcal Z}^*\mathcal H^r=0.
\end{equation}
Here $P_{\mathcal{Z}}$ is defined as in \eqref{arithmetic-limit}. Set $\mathscr Z:=\mathcal{Z}\times_{\mathcal B}\,\Spec R$. We write
$q_R:\mathscr Z\to\mathscr X$ for the base change of $q$ and
$\mathcal T_{i,R}:=\mathcal T_i|_{\mathscr Z}$. Let $P_{\mathscr Z}=P_{\mathcal Z}|_{\Spec R}$ be the localization. 
Set
\[
\gamma_Z:=(\overline{\mathcal H}^{\,r}|Z)>0.
\]
Applying the projection formula to \eqref{qqq}, we have
\begin{equation}\label{q_R}
   \gamma_Z \cdot \big( P_{\mathscr Z}\cdot \mathcal T_{i,R}
                         \cdot(q_R^*\mathcal{A}_R)^{d-1}\big)=0.
\end{equation}
Since $\gamma_Z>0$, this indicates that $P_{\mathscr Z}$ satisfies the hypothesis of Lemma \ref{lem:reduce-to-surface}; see also
Remark~\ref{rem:smR-finite-tests}. Let $\mathscr{X}_s,~\mathscr{X}_{n,s}$ be the special fibers of $\mathscr{X},~\mathscr{X}_n$ respectively. Hence it implies that
\begin{equation*}
    P_{\mathscr X}
       =c_v\,\mathscr{X}_s
\end{equation*}
is a principal divisor for some $c_v\in\mathbb R$. The number $c_v$ is independent of 
the choice of the model $\mathcal{X}$. This also holds on models $\mathscr X$ with 
non-smooth generic fiber.

On $\mathscr{X}_n:=\mathcal{X}_n|_{\Spec R}$, we have $P_{\mathscr X_n}
       =c_v\,\operatorname{div}_{\mathscr X_n}(\varpi_v)$. From \eqref{arith-limit-difference}, we have
\begin{equation}\label{eq:smR-local-model-bound}
    -\varepsilon_n\mathcal E_{n,R}
      \leq\mathcal P_{n,R}
               -c_v\mathscr{X}_{n,s}
      \leq\varepsilon_n\mathcal E_{n,R},
\end{equation}
where $\mathcal E_{n,R}$ and $\mathcal P_{n,R}$ denote the
restrictions to $\mathscr X_n$. 

Let $F_v$ be the completion of $F$ with respect to $v$. We set
$\lambda_v:=-\log|\varpi_v|_v>0$. Denote by $g_{n,v}$, $g_v$ and $g_{E,v}$
the Green functions of $\overline{\mathcal P}_n$,
$\overline P$ and the fixed boundary divisor $\overline{\mathcal{E}}_0$ on 
$V_v^{\mathrm{an}}$ respectively. 
On a normal arithmetic model of $V$, an effective Cartier $\mathbb R$-divisor
has non-negative Green function away from its generic support. Therefore, the inequality ~\eqref{eq:smR-local-model-bound} implies
\[
    |g_{n,v}-\lambda_vc_v|
       \leq\varepsilon_n g_{E,v}
       \]
on $V_v^{\mathrm{an}}$.
Since $g_{E,v}$ is bounded on any compact subset of $V_v^{\mathrm{an}}$, we have  $g_{n,v}\to g_v$ locally uniformly. Hence
$g_v=\lambda_vc_v$ is a constant. By a same normalization step as in Lemma \ref{zero-g}, this constant must be zero. Hence we have
\begin{equation}\label{eq:smR-finite-green-zero}
    g_v=0\quad\text{on }V_v^{\mathrm{an}}.
\end{equation}
    
    Next, we will show that 
    $$
    P_{\mathcal X}=0,
    $$ for every arithmetic model $\mathcal{X}$ of $V$. 
Let $W$ be a prime
divisor on any normal projective arithmetic model $\mathcal{X}$ of $V$.
If it is horizontal over $\mathbb Z$, then the coefficient $\ord_{W}(P_{\mathcal{X}})$ of $P_{\mathcal{X}}$ at $W$ is already
zero by Corollary \ref{arithmetic-limit-0}. Otherwise, W is vertical, and let
$u:=\ord_W$ be its normalized divisorial valuation with value group $\mathbb Z$. The restriction $u|_F$
is nontrivial. We write $u|_F=e\cdot w$, where $e\geq1$ is an integer and
$w(F^\times)=\mathbb Z$. We see that $w$ is a divisorial valuation of 
$F$ by \cite[Lemma 4.5]{HKP21}. Therefore, there is a normal arithmetic model $\widetilde{\mathcal{X}}$ 
dominating $\mathcal{X}$ where the strict transform of $W$ is a component over 
the prime divisor realizing $w$. We have 
$\ord_W(P_{\widetilde{\mathcal{X}}})=0$ by \eqref{eq:smR-finite-green-zero}.
Since the components of the b-divisor associated with $\overline P$ are compatible with projective birational push-forwards, we conclude that
$
P_{\mathcal{X}}=0.
$

Finally, from \eqref{arith-limit-difference}, we know that
\[
    -\varepsilon_n\mathcal E_n
       \leq\mathcal P_n\leq\varepsilon_n\mathcal E_n,
\]
for $n\geq 0$.
The archimedean Green function of $\overline{P}$ is zero by
Corollary \ref{arithmetic-limit-0}. Combining with the geometric part, we have
\[
    -\varepsilon_n\overline{\mathcal E}_n
       \leq\overline{\mathcal P}_n
       \leq\varepsilon_n\overline{\mathcal E}_n.
\]
This means that $\overline{P}=0$ as an adelic divisor.
\end{proof}

Therefore, we have shown the vertical case of Theorem \ref{Rdr}, summarized in the following corollary.
\begin{corollary}[Vertical case]
\label{cor:smR-vertical}
Let $F$ be a finitely generated 
field of transcendence degree $r\ge 0$ over $\mathbb{Q}$, and $V$ be a normal, geometrically integral quasi-projective variety over 
$F$ of dimension $d$. Suppose 
$\overline M\in\widehat{\mathrm{Pic}}(V/\mathbb Z)
_{\mathrm{int},\mathbb Q}$ is vertical. Fix a point $x_0\in V(\overline{F})$. If $\overline M$ satisfies the vanishing hypothesis in
\eqref{hypotheses}, then
\[
    \overline M=\pi_V^*\overline N_0,
\]
for $\overline N_0:=\mathfrak h_{\overline M}(x_0)$.
\end{corollary}

\begin{proof}
Let $\overline D$ be an 
adelic divisor represented by
$\{\overline{\mathcal D}_n\}$ with $\mathcal D_{n,F}=0$ on every normal projective
model $X_n$ of $V$.
Let $\overline{M}^\prime:=\overline M-\pi_V^*\overline N_0 $. And $\overline P:=\overline D-\pi_V^*\overline\Delta_0$
be its associated divisor introduced in $\S$\ref{subsec-pushforward}. $\overline{M}^\prime$ also satisfies the vanishing hypothesis in
\eqref{hypotheses}.
Corollary~\ref{Cor-Orthogonality} implies that $\overline P$ satisfies 
orthogonality. Corollary \ref{arithmetic-limit-0} and Proposition \ref{prop:smR-zero-limit} together imply that we have $\overline P=0$.
Therefore, the associated adelic line bundle $\overline{M}^\prime$=0, which 
proves the assertion.
\end{proof}

\subsection{General case}
\label{subsec-Rdr general}

Let $F/\mathbb Q$ be finitely generated of transcendence degree
$r$, and let $V/F$ be normal and geometrically integral of dimension
$d\geq2$. 
In this subsection, we will show that for an arbitrary adelic line 
bundle $\overline{M}$ on $V$ satisfying \eqref{hypotheses}, must be a 
restriction of some flat adelic line bundle on a projective model of $V$.

We take a smooth projective compactification $X\supset V$
and a point $x_0\in V(\overline F)$. Let $x_0^\prime \in V$ be the associated closed point. Write
$\pi_X:X\to\Spec F$ and $\pi_V:V\to\Spec F$ for the structure
morphisms. 
Let $\overline M$ be an integrable adelic line bundle on $V$ which satisfies \eqref{hypotheses}. Set
$\overline N_0:=\mathfrak h_{\overline M}(x_0)$ and
$\overline M':=\overline M-\pi_V^*\overline N_0$. Then 
$\mathfrak{h}_{\overline{M}^\prime}(x_0)=0$. Let
$\overline P$ be the adelic divisor associated with
$\overline M'$ in \S\ref{subsec-pushforward}. 
Denote the component of the b-divisor $b(\overline P)$ on $X$ by $P_X$, which is introduced in \S \ref{subsec-pushforward}. By Lemma~\ref{lem-limit divisor}, $P_X$ is a numerically
trivial Cartier $\mathbb{R}$-divisor and its pullbacks induces all components of $b(\overline P)$.

We recall the notion of flat metrized line bundles on projective varieties over a number field $K$ in \cite[Appendix A.4]{YZ17}. A numerically
trivial line bundle $L$ on a smooth projective variety $Y/K$ belongs
to $\mathrm{Pic}^0(Y)_{\mathbb Q}$. Then we can endow $L$ with a flat metric $\|\cdot\|$ so that the integrable metrized line bundle $\overline{L}=(L,\|\cdot\|)$ is flat, i.e.
$$
c_1(\overline{L}|_C)=0,
$$
on $C^{\mathrm{an}}$ for any morphism $f:C\to Y$ from a projective curve $C/K$. In \cite[\S4.3]{Car20A}, Carney 
further generalizes the notion of flatness to adelic line bundles on projective varieties over finitely 
generated fields using a Tate limit argument. By \cite[Definition 4.2]{Car20A}, an integrable adelic line bundle 
$\overline{L}$ on a projective variety $X$ over finitely generated $F$ is called \emph{flat} if
$$
\overline{L} \cdot \overline{N}\cdot \overline{L}_1 \cdots \overline{L}_{d-1}\equiv_{\mathrm{num}} 0
$$
for any vertical adelic line bundle $\overline{N}$ on $X$ and any integrable $\overline L_{i}$ on $X$. We call 
an integrable adelic divisor \emph{flat} if its associated adelic line bundle is flat. We 
summarize the form of \cite[Lemma 4.3]{Car20A} that can be applied to our settings as follows.
\begin{lemma}[Flat lifts]
\label{lem:smR-flat-lift}
Let $G$ be a numerically trivial $\mathbb{Q}$-Cartier divisor on a
smooth geometrically integral projective variety $X/F$ such
that $x_0'\notin\operatorname{Supp}(G)$. There is a flat integrable adelic $\mathbb{Q}$-line bundle 
$\mathcal{O}(\overline G)$ with underlying $\mathbb Q$-line bundle 
$\mathcal O(G)$. Moreover, we have
\[
    \mathfrak{h}_{\mathcal{O}(\overline G)}(x_0)=0.
\]
The construction of $\overline{G}$ is preserved by pullbacks to normal projective models of $X$. In particular, let $f:Y\to X$ be a normal projective
birational model, and let $\pi_Y:Y\to\Spec F$ be its structure
morphism. If $\overline T$ is an integrable vertical adelic line
bundle on $Y$, $\overline A_1,\ldots,\overline A_{d-1}$ are
integrable adelic line bundles on $Y$, and
$\overline H_1,\ldots,\overline H_r$ are integrable adelic line
bundles on $\Spec F$, then
\begin{equation}\label{eq:flat-orthogonal}
    f^*\mathcal O(\overline G)\cdot\overline T
       \cdot\overline A_1\cdots\overline A_{d-1}
       \cdot\pi_Y^*\overline H_1\cdots\pi_Y^*\overline H_r=0.
\end{equation}
\end{lemma}

\begin{proof}
Choose an integer $a>0$ for which $aG$ is integral and
algebraically trivial. Then we can find a $\overline{G}$ by endowing $aG$ with a flat metric as 
in \cite[Lemma 4.3]{Car20A} and then dividing by $a$. Let $\mathcal{Y} \to \mathcal{B}$ be an arithmetic model of $V\to \Spec F$. This implies
that \eqref{eq:flat-orthogonal} holds for $\overline{G}$. Therefore, We have
$$
c_1( \mathcal O(\overline G))c_1(\overline{A}_1) \cdots c_1(\overline{A}_{d-1})\pi_Y^*c_1(\overline{H}_1)\cdots \pi_Y^*c_1(\overline{H}_r) =0 
$$
on $\mathcal{Y}^\mathrm{an}$. 

Finally, we can replace $\mathcal{O}(\overline{G})$ by itself tensoring with the inverse of $\pi_V^* \mathfrak{h}_{\mathcal{O}(\overline G)}(x_0)$.
It changes neither the geometric part $\mathcal{O}(G)$ nor the flatness since $\overline{G}$ is flat. Consequently, we have $\mathfrak{h}_{\mathcal{O}(\overline G)}(x_0)=0.$
\end{proof}

\begin{proposition}
\label{prop:smR-flat-reduction}
Let $V,~F,~ \overline P$ be as in \S\ref{subsec-Rdr-ineq}. We set $G:=P_X$ be the 
component of the b-divisor associated with $\overline{P}$ over the normal projective 
model $X$ of $V$. Then there is an adelic $\mathbb{R}$-divisor $\overline G$ as a $\mathbb R$-linear
combination of flat adelic divisors on $X$, with underlying divisor $G$, such that
\begin{equation}\label{eq:smR-flat-reduction}
    \overline P=\overline G|_V.
\end{equation}
\end{proposition}

\begin{proof} First, $G$ is numerically trivial by Lemma \ref{lem-limit divisor}. 
We fix a closed point $x_0^\prime$ of $V$ and its corresponding $x_0\in V(\overline{F})$.
Let $B_1,\ldots,B_t$ be the prime components of the boundary
$X\setminus V$. Choose a $\mathbb Q$-Cartier divisor $G_0$ on $X$
extending the underlying divisor $P$ of $\overline P$ on $V$, with
support avoiding $x_0'$. The construction of the geometric b-divisor shows that
$G$ is an $\mathbb R$-linear composition of $G_0,B_1,\ldots,B_t$.
Let $\left[G_0\right]$ and $\left[B_i\right]$ be their images in the group $N^1(X)_{\mathbb Q}$ of numerically equivalent $\mathbb Q$-classes. We have a linear map
\[
  \varphi:  \mathbb{Q}G_0+  \mathbb{Q}B_1+\cdots+ \mathbb{Q}B_t  \longrightarrow N^1(X)_{\mathbb Q}.
\]
Then $G$ is contained in $\ker(\varphi) \otimes_{\mathbb{Q}}\mathbb{R}$. 
We can choose a $\mathbb{Q}$-basis $G_1,\ldots,G_\rho$ of this kernel,
where $\rho=\dim_{\mathbb{Q}}\ker(\varphi)$. Every $G_j$ is numerically trivial, and its support avoids $x_0'$. And there are
$a_1,\ldots,a_\rho\in\mathbb R$ such that
\[
    G=\sum_{j=1}^{\rho}a_jG_j.
\]
We can choose the flat lifts $\overline G_j$ from
Lemma~\ref{lem:smR-flat-lift}, and set
\[
    \overline G:=\sum_{j=1}^{\rho}a_j\overline G_j,
    \quad \text{and}\quad
    \overline{P}^\prime:=\overline P-\overline G|_V.
\]

By \eqref{eq:smR-trace-pullback}, the geometric b-divisor of
$\overline{P}^\prime$ is zero. We have 
$\mathfrak{h}_{\mathcal{O}(\overline{P}^\prime)}(x_0)=0$ 
by the constructions of these $\overline{G}_j$. 
Also we know that $\mathcal O(\overline P)=\overline M'$ satisfies 
orthogonality by Corollary~\ref{Cor-Orthogonality}. Since the pullbacks of the 
$\overline{G}_j$ to a projective model dominating $X$ are still flat, 
we see that $\mathcal O(\overline G|_V)$ also satisfies the equality \eqref{eq:flat-orthogonal} and so does the difference $\overline P^\prime$.
Consequently, we have $\overline P^\prime=0$ by Proposition~\ref{prop:smR-zero-limit}.
\end{proof}

To reduce the general case to the case of curves, we need the following two lemmas.

\begin{lemma}[Restriction]
\label{lem:smR-flat-restriction}
Let $X/F$ be a smooth geometrically integral projective variety of
dimension $d\geq2$, and let $\overline G$
be a flat adelic $\mathbb{R}$-divisor on $X$. Set
$\overline L:=\mathcal O(\overline G)$.
Choose line bundles $A_1,\ldots,A_{d-1}$ on $X$ with nonzero sections
$s_i\in H^0(X,A_i)$ for $1\leq i\leq d-1$. Suppose that we have a flag
\[
    Y_0:=X,\qquad
    Y_i:=Y_{i-1}\cap\operatorname{div}(s_i)
       \quad(1\leq i\leq d-1)
\]
such that every $Y_i$ is smooth, geometrically integral, and
an effective Cartier divisor of multiplicity $1$ on $Y_{i-1}$. Write
$\overline C:=Y_{d-1}$, and denote by
$\pi_0:X\to\Spec F$ and
$\pi_{\overline C}:\overline C\to\Spec F$ the structure morphisms.
 Then for any model adelic line bundles $\overline A_i$ with underlying line 
 bundles $A_i$, and every integrable $r$-tuple of adelic line bundles
$\overline H_1,\ldots,\overline H_r$ on $\Spec F$, one has
\begin{equation}\label{eq:smR-flat-slicing}
\begin{aligned}
    &(\overline L|_{\overline C})^2
        \cdot\pi_{\overline C}^*\overline H_1\cdots
              \pi_{\overline C}^*\overline H_r\\
    &\qquad=
       \overline L^2\cdot\overline A_1\cdots\overline A_{d-1}
        \cdot\pi_0^*\overline H_1\cdots\pi_0^*\overline H_r.
\end{aligned}
\end{equation}
\end{lemma}
\begin{proof}
Let $\pi_i: Y_i\to \Spec F$ be the structure morphism of $Y_i$ for $1\le i \le d-1$ and $\pi_{d-1}=\pi_{\overline{C}}$. We only need to apply the induction 
formula of arithmetic intersection products successively along the flag.
From $Y_{i-1}$ to $Y_i~(1\le i\le d-1)$, we have
\begin{align*}
           &\overline L|_{Y_{i-1}}^2\cdot\overline A_i|_{Y_{i-1}} \cdots\overline A_{d-1}|_{Y_{i-1}}
        \cdot\pi_{i-1}^*\overline H_1\cdots  \pi_{i-1}^*\overline H_{r}
        \\
        =&
        \overline L|_{Y_{i}}^2\cdot\overline A_{i+1}|_{Y_{i}}\cdots\overline A_{d-1}|_{Y_{i}}
        \cdot\pi_{i}^*\overline H_1\cdots\pi_{i}^*\overline H_{r}
        -\sum_{v\in M_{\mathbb{Q}}} \int_{{\mathcal{Y}}_{i-1,v}^{\mathrm{an}}}
        \log\|s_{i}\|_v \mu_{i-1}  ,\\
       &\mathrm{where}~ \mu_{i-1}:=c_1(\overline L|_{Y_{i-1}})_v^2\prod_{j=i+1}^{d-1} c_1(\overline A_{j}|_{Y_{i-1}})_v  
       \prod_{k=1}^{r}\pi_{i-1}^* c_1(\overline H_{k})_v .
\end{align*}
Here, $\mathcal{Y}_{i-1}$ is an arithmetic model of $Y_{i-1}/F$ and ${\mathcal{Y}}_{i-1,v}^{\mathrm{an}}$ is the Berkovich analytification of $\mathcal{Y}_{i-1}$ with respect to $v$.
Since $\overline{L}$ is flat, the integration terms vanish for all $i$. This proves the lemma. 
\end{proof}

\begin{lemma}[Reducing to a complete intersection curve]
\label{reduction-complete-intsct-curve}
Let $F$, $V$, $X$, and $\overline M$ be as in the beginning of \S\ref{subsec-Rdr general}. Choose a flag on $X$ with a complete intersection curve
$\overline C$ as in Lemma~\ref{lem:smR-flat-restriction}, and assume $C:=\overline C\cap V\ne\varnothing$.
Let $j:C\hookrightarrow V$ be the inclusion and
$\pi_C:C\to\Spec F$ the structure morphism. Then we have
\[
     \pi_{C,*}\langle j^*\overline M,j^*\overline M\rangle
       \equiv_{\mathrm{num}}0,
    \quad \mathrm{and} \quad
    \deg(j^*\widetilde M)=0.
\]
\end{lemma}

\begin{proof} Fix $x_0 \in V(\overline{F})$ corresponding to a regular closed 
point $x_0^\prime$ of $V$. Let 
$\overline{N}_0:=\mathfrak h_{\overline{M}}(x_0)$. We set $\overline M':=\overline M-\pi_V^*\overline N_0$, and the 
corresponding divisor $\overline P$ as in \S\ref{subsec-pushforward}. Then  
 $\overline P=\overline G|_V$
for a flat adelic divisor on $X$
 from Proposition~\ref{prop:smR-flat-reduction}.
Let
$\overline A_1,\ldots,\overline A_{d-1}$ be the model adelic
bundles that give rise to this flag. Since $\overline P$ satisfies
\eqref{hypotheses}, we have
\[
    \mathcal O(\overline G)^2\cdot\overline A_1\cdots\overline A_{d-1}
       \cdot\pi_X^*\overline H_1\cdots\pi_X^*\overline H_r=0
\]
for every integrable base tuple
$\overline H_1,\ldots,\overline H_r$. The above lemma therefore shows that
\[
    (j^*\overline M')^2
       \cdot\pi_C^*\overline H_1\cdots\pi_C^*\overline H_r=0.
\]
The geometric part $j^*\widetilde M$ of $j^*\overline M'$ is represented by the
restriction of $\mathcal O_X(G)$ to
$\overline C$. Its degree is zero because $G$ is numerically
trivial on $X$. The geometric part of a base pullback is zero,
so $\deg(j^*\widetilde M)=0$ as well. Finally, since we have
\[
    j^*\overline M=j^*\overline M'+\pi_C^*\overline N_0,
\]
the projection formula ensures that
\[
    (j^*\overline M)^2
       \cdot\pi_C^*\overline H_1\cdots\pi_C^*\overline H_r=0
\]
for every integrable $r$-tuple $\overline{H}_1,\dots,\overline{H}_r$ as required.
\end{proof}

We can now apply the already established curve theorem $R(1,r)$ and the height 
theorem $H(d,r)$ to $j^*\overline{M}=\overline{M}|_{C}$ and $\overline{M}$.

\begin{proof}[Proof of Theorem \ref{chap1-Rdr}]
We fix two points $x_0,~x\in V(\overline F)$ associated with regular closed points $x_0^\prime,~x^\prime$ of $V$. and put
$\overline N_0:=\mathfrak h_{\overline M}(x_0)$. We then choose a smooth complete intersection curve $\overline C\subset V_{\mathrm{reg}}$ through 
the two closed points. Here, the regular locus ${V_{\mathrm{reg}}}$ is a non-empty open subvariety of $V$ Applying Lemma \ref{reduction-complete-intsct-curve} and the assertion $R(1,r)$ to the restriction of $\ol M$ on $C= \ol C\cap V_{\mathrm{reg}}$, we have
\[
    \overline M|_C=\pi_C^*\overline N_0.
\]
This shows that $\mathfrak h_{\overline M}(x)=\mathfrak h_{\overline M}(x_0)=\ol N_0$. Varying $x\in V_{\mathrm{reg}}(\ol F)$, we see that
$$
\mathfrak h_{\overline M}(x)=\ol N_0,
$$
for every $x\in V_{\mathrm{reg}}(\ol F)$.
Therefore, we can apply $H(d,r)$ to $(\overline M-\pi_V^*\overline N_0)|_{V_{\mathrm{reg}}}$ to see that 
$$ 
(\overline M-\pi_V^*\overline N_0)|_{V_{\mathrm{reg}}}= \pi_V^*\overline N_1|_{V_{\mathrm{reg}}},
$$
for some adelic $\mathbb Q$-line bundle $\ol N_1$ on $\Spec F$. Moreover, by \cite[Corollary 3.4.2]{YZ26}, the restriction map
$$
\Pic(V/\mathbb Z) \longrightarrow \Pic(V_{\mathrm{reg}}/\mathbb Z)
$$
induced by the open immersion $V_{\mathrm{reg}} \hookrightarrow V$ is injective. Consequently, we have 
$$
\overline M=\pi_V^*(\overline N_0+\overline N_1)
$$ 
on $V$, which concludes the proof of Theorem \ref{chap1-Rdr}.
\end{proof}

\subsection{Proof of Theorem \ref{num-triv}} \label{Hodge-implies-a-to-c}
Now we can deduce the equivalence of three conditions of Theorem~\ref{num-triv} from the preceding height results and the Hodge index type theorem as follows.

\begin{proof}[Proof of Theorem~\ref{num-triv}]
We first prove $(a)\Rightarrow(c)$.
For integrable $\overline{A}_1,\ldots,\overline{A}_{d-1}$ on $V$,
put
$\delta=\widetilde{M}\cdot\widetilde{A}_1\cdots
\widetilde{A}_{d-1}$.
Let $\overline{H}_0,\dots,\overline{H}_r$ be $r+1$ nef adelic line bundles on $\Spec F$ with positive intersection product. Then we have
$$
0=\pi_{V,*}\langle\overline{M},\pi_V^*\overline{H}_0,
\overline{A}_1,\ldots,\overline{A}_{d-1}\rangle \cdot\overline{H}_1\cdots\overline{H}_r= \delta\,
\overline{H}_0\cdot\overline{H}_1\cdots\overline{H}_r,
$$
by \cite[Lemma 4.6.1(3)]{YZ26}.
Hence we have $\delta=0$. Taking one of the test bundles in $(a)$ to be
$\overline{M}$ provides that
\[
\pi_{V,*}\langle\overline{M},\overline{M},
\overline{A}_1,\ldots,\overline{A}_{d-1}\rangle
\equiv_{\mathrm{num}}0.
\]
Applying Theorem~\ref{Rdr} to $\overline{M}$ shows that
$\overline{M}=\pi_V^*\overline{N}$ when $d\geq2$;
for $d=1$, use the corresponding curve equality theorem $R(1,r)$.
Choose $d$ nef adelic line bundles
$\overline{B}_1,\ldots,\overline{B}_d$ on $V$ whose geometric
intersection $b=\widetilde{B}_1\cdots\widetilde{B}_d>0$.
The projection formula and $(a)$ imply that
\[
0= \pi_{V,*}\langle\pi_V^*\overline{N},
\overline{B}_1,\ldots,\overline{B}_{d}\rangle\cdot \overline{H}_1\cdots\overline{H}_r =b\,(\overline{N}\cdot\overline{H}_1\cdots\overline{H}_r)
\]
for every integrable tuple $\overline{H}_i$. Since $b>0$, we have
$\overline{N}\equiv_{\mathrm{num}}0$.

For $(b)\Rightarrow(c)$, Theorem~\ref{H(d,r)} gives
$\overline{M}=\pi_V^*\overline{N}$.
Lemma~\ref{height2} identifies
$\mathfrak{h}_{\overline{M}}(x)$ with $\overline{N}$, so
$\overline{N}\equiv_{\mathrm{num}}0$ by $(b)$.

Conversely, $(c)\Rightarrow(b)$ follows immediately from
Lemma~\ref{height2}. Finally, for arbitrary integrable
$\overline{L}_1,\ldots,\overline{L}_d$ on $V$ and
$\overline{H}_1,\ldots,\overline{H}_r$ on $F$, the projection
formula gives
\begin{align*}
&\pi_{V,*}\langle\pi_V^*\overline{N},
\overline{L}_1,\ldots,\overline{L}_d\rangle
\cdot\overline{H}_1\cdots\overline{H}_r\\
&\qquad=
(\widetilde{L}_1\cdots\widetilde{L}_d)
(\overline{N}\cdot\overline{H}_1\cdots\overline{H}_r)=0.
\end{align*}
Hence $\pi_{V,*}\langle\pi_V^*\overline{N},
\overline{L}_1,\ldots,\overline{L}_d\rangle$ is numerically trivial. This proves $(c)\Rightarrow(a)$ and completes the proof of Theorem \ref{num-triv}.
\end{proof}

\newpage

\end{document}